\documentclass[11pt,reqno]{amsart}
\usepackage[cp1251]{inputenc}
\usepackage[english]{babel}
\usepackage{amsmath}
\usepackage{amsfonts}
\usepackage{amssymb}
\usepackage[dvipsnames]{xcolor}
\usepackage{mathrsfs} 
\usepackage{tikz}
\usepackage{pgfplots}
\usepackage{float}
\floatstyle{boxed}
\restylefloat{figure}
\usepackage{amsthm}
\theoremstyle{theorem}
\newtheorem{theorem}{Theorem}[section]
\newtheorem{proposition}[theorem]{Proposition}
\newtheorem{lemma}[theorem]{Lemma}
\newtheorem{condition}[theorem]{Condition}
\newtheorem{remark}[theorem]{Remark}

\newtheorem{example}[theorem]{Example}

\usepackage[paperwidth=20.5cm,paperheight=29.7cm,width=16cm,height=25.5cm,top=24mm,left=19mm,headsep=11mm]{geometry}

\numberwithin{equation}{section}

\DeclareMathOperator*{\esssup}{ess\,sup}

\begin{document}

\title[Approximation for solutions of parabolic systems]{Homogenization of periodic parabolic
systems \\in the $L_2(\mathbb{R}^d)$-norm \\with the corrector taken into account}

\author{Yu. M. Meshkova}

\address{St. Petersburg 
State University\\
Chebyshev Laboratory\\
199178, 14th Line V.~O.~29B\\St.~Petersburg,
Russia; \textit{Current address of the author: }Department of Mathematics and Statistics, University of Helsinki,
P.~O. Box 68 (Gustaf
H\"allstr\"omin katu 2), FI-00014 Helsinki, Finland; 
 \textit{Current address of Chebyshev Laboratory: }199178, 14th Line V.~O.~29\\St.~Petersburg,
Russia.}

\email{y.meshkova@spbu.ru, iuliia.meshkova@helsinki.fi}



\thanks{The research was supported by the Russian Science Foundation, grant  no.~14-21-00035.}

\begin{abstract}
In $L_2(\mathbb{R}^d;\mathbb{C}^n)$, consider a self-adjoint matrix second order elliptic  differential operator $\mathcal{B}_\varepsilon$, $0<\varepsilon \leqslant 1$. The principal part of the operator is given in a factorised form, the operator contains first and zero order terms. The operator $\mathcal{B}_\varepsilon$ is positive definite, its coefficients are periodic and depend on $\mathbf{x}/\varepsilon$. We study the behaviour in the small period limit of the operator exponential $e^{-\mathcal{B}_\varepsilon t}$, $t\geqslant 0$. The approximation in the $(L_2\rightarrow L_2)$-operator norm with error estimate of order  $O(\varepsilon ^2)$ is obtained. The corrector is taken into account in this approximation. The results are applied to homogenization of the solutions for the Cauchy problem for parabolic systems.
\end{abstract}

\maketitle

\section*{Introduction}

The research is devoted to homogenization of periodic differential operators (DO's). A broad literature is devoted to homogenization problems; see, for example, the books
 \cite{BaPa,BeLP,ZhKO,Sa}. We rely on the spectral approach to homogenization
problems based on the Floquet-Bloch theory and the analytic perturbation theory. This approach was developed in the series of papers  \cite{BSu0,BSu,BSu05-1,BSu05,BSu06} by M.~Sh.~Birman and T.~A.~Suslina. A great deal of considerations  is carried out in abstract operator-theoretic terms.

\subsection{Problem setting} We study the behavior in the small period limit of the solution for the parabolic system 
\begin{equation}
\label{Introduction Cauchy problem}
\begin{cases}
G(\mathbf{x}/\varepsilon )\partial _s\mathbf{u}_\varepsilon (\mathbf{x},s)=-\mathcal{B}_\varepsilon \mathbf{u}_\varepsilon (\mathbf{x},s),\quad\mathbf{x}\in\mathbb{R}^d,\;s>0;\\
G(\mathbf{x}/\varepsilon)\mathbf{u}_\varepsilon (\mathbf{x},0)=\boldsymbol{\phi}(\mathbf{x}),\quad\mathbf{x}\in\mathbb{R}^d .
\end{cases}
\end{equation}
Here $\boldsymbol{\phi}\in L_2(\mathbb{R}^d;\mathbb{C}^n)$, $\mathcal{B}_\varepsilon$ is a matrix elliptic second order DO acting in $L_2(\mathbb{R}^d;\mathbb{C}^n)$. The measurable $(n\times n)$-matrix-valued function $G(\mathbf{x})$ is assumed to be bounded, uniformly positive definite, and periodic with respect to some lattice $\Gamma \subset \mathbb{R}^d$. The coefficients of the operator $\mathcal{B}_\varepsilon$ are $\Gamma$-periodic functions depending on $\mathbf{x}/\varepsilon$. For any measurable $\Gamma$-periodic function $\varphi (\mathbf{x})$, $\mathbf{x}\in\mathbb{R}^d$, we use the notation $\varphi ^\varepsilon (\mathbf{x})=\varphi (\mathbf{x}/\varepsilon)$.

The principal part $\mathcal{A}_\varepsilon$ of the operator $\mathcal{B}_\varepsilon$ is given in a factorized form
\begin{equation}
\label{Introduction A_eps}
\mathcal{A}_\varepsilon =b(\mathbf{D})^* g^\varepsilon (\mathbf{x})b(\mathbf{D}),
\end{equation}
where $b(\mathbf{D})$ is a matrix homogeneous first order DO, $g(\mathbf{x})$ is $\Gamma$\nobreakdash-periodic bounded and positive definite matrix-valued function in~$\mathbb{R}^d$. (The precise assumptions on the coefficients of the operator  $\mathcal{B}_\varepsilon$ are given below in~\S\ref{Section Pervaritel'nye svedeniya}.) The operator $\mathcal{B}_\varepsilon$ contains first and zero order terms
\begin{equation}
\label{Introduction hat-B-eps}
\mathcal{B}_\varepsilon \mathbf{u} 
=\mathcal{A}_\varepsilon \mathbf{u}
+\sum \limits _{j=1}^d\left(a_j^\varepsilon(\mathbf{x})D_j\mathbf{u}
+D_j(a_j^\varepsilon (\mathbf{x}))^*\mathbf{u}\right)
+\mathcal{Q}^\varepsilon (\mathbf{x})\mathbf{u}+\lambda \mathbf{u}.
\end{equation}
Here $a_j(\mathbf{x})$, $j=1,\dots,d$, are $\Gamma$-periodic $(n\times n)$-matrix-valued functions belonging to a suitable $L_\varrho$-space on the cell $\Omega$ of the lattice $\Gamma$. In general, the potential $\mathcal{Q}^\varepsilon (\mathbf{x})$ is a distribution generated by a rapidly oscillating measure with values in the Hermitian matrices. (Since the coefficients are not assumed to be bounded, the precise definition of the operator  \eqref{Introduction hat-B-eps} is given via the quadratic form.) The constant $\lambda$ is chosen such that the operator $\mathcal{B}_\varepsilon$ is positive definite. The coefficients of the operator $\mathcal{B}_\varepsilon$ oscillate rapidly as $\varepsilon \rightarrow 0$. 

The principal term of approximation for the solution of the problem \eqref{Introduction Cauchy problem} for small $\varepsilon $ was obtained in \cite{M_AA}:
\begin{equation}
\label{Introduction principal tterm of approx u-eps}
\Vert \mathbf{u}_\varepsilon (\,\cdot\, ,s)-\mathbf{u}_0(\,\cdot\, ,s)\Vert _{L_2(\mathbb{R}^d)}\leqslant C_1\varepsilon (s+\varepsilon ^2)^{-1/2}e^{-C_2 s}\Vert \boldsymbol{\phi}\Vert _{L_2(\mathbb{R}^d)},\quad s\geqslant 0.
\end{equation} 
Here $\mathbf{u}_0$ is the solution of the ``homogenized'' problem
\begin{equation*}
\begin{cases}
\overline{G}\partial _s \mathbf{u}_0(\mathbf{x},s)=-\mathcal{B}^0\mathbf{u}_0(\mathbf{x},s),& \mathbf{x}\in\mathbb{R}^d,\quad s>0;
\\
\overline{G}\mathbf{u}_0(\mathbf{x},0)=\boldsymbol{\phi}(\mathbf{x}),&\mathbf{x}\in\mathbb{R}^d,
\end{cases}
\end{equation*}
where $\mathcal{B}^0$ is the \textit{effective operator }with constant coefficients, $\overline{G}$ is the mean value of the matrix-valued function $G$ over the cell of periodicity: $\overline{G}=\vert \Omega\vert ^{-1}\int _\Omega G(\mathbf{x})\,d\mathbf{x}$. In the case when $G(\mathbf{x})=\mathbf{1}_n$, estimate \eqref{Introduction principal tterm of approx u-eps} admits an obvious formulation in operator terms 
\begin{equation*}
\Vert e^{-\mathcal{B}_\varepsilon s}-e^{-\mathcal{B}^0 s}\Vert _{L_2(\mathbb{R}^d)\rightarrow L_2(\mathbb{R}^d)}
\leqslant  C_1\varepsilon (s+\varepsilon ^2)^{-1/2}e^{-C_2 s}.
\end{equation*}
Thus, results of such type are called \textit{operator error estimates} in homogenization theory. 
(In the general case, we deal with the ``bordered'' exponential $f^\varepsilon e^{-(f^\varepsilon )^*\mathcal{B}_\varepsilon f^\varepsilon s}(f^\varepsilon )^*$, where ${G^{-1}=ff^*}$ and the matrix-valued function $f$ is assumed to be $\Gamma$-periodic.) 

\textit{Our aim} is to refine approximation \eqref{Introduction principal tterm of approx u-eps} by taking the \textit{corrector} into account. In other words, for a fixed time  $s>0$, to approximate the solution  $\mathbf{u}_\varepsilon (\,\cdot\, ,s)$ of problem \eqref{Introduction Cauchy problem} in the $L_2(\mathbb{R}^d;\mathbb{C}^n)$-norm with the error estimate of the order $O(\varepsilon ^2)$. 

\subsection{Main results} In Introduction, we discuss only the case when  $G(\mathbf{x})=\mathbf{1}_n$. 
\textit{Our main result} is the estimate
\begin{equation}
\label{Introduction Main Result}
\Vert e^{-\mathcal{B}_\varepsilon s}-e^{-\mathcal{B}^0s}-\varepsilon \mathcal{K}(\varepsilon ,s)\Vert _{L_2(\mathbb{R}^d)\rightarrow L_2(\mathbb{R}^d)}
\leqslant
C_3\varepsilon ^2 (s+\varepsilon ^2)^{-1}e^{-C_2 s},\quad s>0.
\end{equation}
Here $\mathcal{K}(\varepsilon ,s)$ is the corrector. 
It turns out that it is equal to sum of three terms and the two terms of the corrector are mutually conjugate and contain rapidly oscillating factors. The third summand has constant coefficients.  

\subsection{Survey on operator error estimates} 
An interest to operator error estimates was sparked by the paper \cite{BSu0} of M.~Sh.~Birman and T.~A.~Suslina, where for the resolvent of the operator  \eqref{Introduction A_eps} it was shown that
\begin{equation}
\label{Introduction res A-eps principal}
\Vert (\mathcal{A}_\varepsilon +I)^{-1}-(\mathcal{A}^0 +I)^{-1}\Vert _{L_2(\mathbb{R}^d)\rightarrow L_2(\mathbb{R}^d)}
\leqslant C\varepsilon .
\end{equation}
Here $\mathcal{A}^0=b(\mathbf{D})^*g^0b(\mathbf{D})$ is the effective operator, $g^0$ is the constant effective matrix. Approximation for the operator $(\mathcal{A}_\varepsilon +I)^{-1}$ with the corrector taken into account was obtained in \cite{BSu05}:
\begin{equation}
\label{Introduction res A-eps L2-corrector}
\Vert (\mathcal{A}_\varepsilon +I)^{-1}-(\mathcal{A}^0 +I)^{-1} -\varepsilon K(\varepsilon)\Vert _{L_2(\mathbb{R}^d)\rightarrow L_2(\mathbb{R}^d)}
\leqslant C\varepsilon ^2.
\end{equation}
Approximation for the resolvent $(\mathcal{A}_\varepsilon +I)^{-1}$ in the energy norm was achieved in \cite{BSu06}:
\begin{equation}
\label{Introduction res A-eps H1-corrector}
\Vert (\mathcal{A}_\varepsilon +I)^{-1}-(\mathcal{A}^0 +I)^{-1} -\varepsilon K_1(\varepsilon)\Vert _{L_2(\mathbb{R}^d)\rightarrow H^1(\mathbb{R}^d)}
\leqslant C\varepsilon .
\end{equation}
Here $K_1(\varepsilon)$ corresponds to the traditional corrector in homogenization theory. This operator involves rapidly oscillating factors, so  $\varepsilon\Vert K_1(\varepsilon)\Vert _{L_2\rightarrow H^1}=O(1)$. Note that the corrector in  \eqref{Introduction res A-eps L2-corrector} has the form $K(\varepsilon)=K_1(\varepsilon )+K_1(\varepsilon )^* +K_3$, moreover $K_3$ does not depend on $\varepsilon$. Later, estimates
 \eqref{Introduction res A-eps principal}--\eqref{Introduction res A-eps H1-corrector} were generalized by T.~A.~Suslina \cite{Su10,Su14} to a wider class of operators \eqref{Introduction hat-B-eps}.

The spectral approach was applied to parabolic systems  
in the papers \cite{Su04,Su07} by T.~A.~Suslina, where the principal term of approximation for the operator $e^{-\mathcal{A}_\varepsilon s}$ was obtained, in the paper  \cite{V} by E.~S.~Vasilevskaya, where approximation in the $L_2(\mathbb{R}^d;\mathbb{C}^n)$-operator norm with the corrector taken into account was achieved, and also in  \cite{Su_MMNP}, where the approximation for the operator exponential in the energy norm was obtained:
\begin{align}
\label{Introduction exp-A-eps principal}
\Vert e^{-\mathcal{A}_\varepsilon s}\!-\!e^{-\mathcal{A}^0s}\Vert _{L_2(\mathbb{R}^d)\rightarrow L_2(\mathbb{R}^d)}&\!\leqslant\! C\varepsilon (s\!+\!\varepsilon ^2)^{-1/2},\quad s\!\geqslant \!0;
\\
\label{Introduction Vasilevskaya}
\Vert e^{-\mathcal{A}_\varepsilon s}\!-\!e^{-\mathcal{A}^0s}-\varepsilon \mathcal{K}(\varepsilon ,s)\Vert _{L_2(\mathbb{R}^d)\rightarrow L_2(\mathbb{R}^d)}&\!\leqslant\! C\varepsilon ^2 (s\!+\!\varepsilon ^2)^{-1},\quad s\!\geqslant\! 0;
\\
\label{Introduction exp-A-eps H1}
\Vert e^{-\mathcal{A}_\varepsilon s}\!-\!e^{-\mathcal{A}^0s}\!-\!\varepsilon \mathcal{K}_1(\varepsilon ,s)\Vert _{L_2(\mathbb{R}^d)\rightarrow H^1(\mathbb{R}^d)}&\!\leqslant \!C\varepsilon (s^{-1/2}\!+\!s^{-1}),\quad s\!\geqslant\! \varepsilon ^2.
\end{align}
There are no exponentially decaying factors in these estimates because the point zero is the lower edge of the spectra for the operators $\mathcal{A}_\varepsilon$ and $\mathcal{A}^0$. Generalization of estimates~\eqref{Introduction exp-A-eps principal} and \eqref{Introduction exp-A-eps H1} to the operator \eqref{Introduction hat-B-eps} was achieved in the paper \cite{M_AA} by the author. Generalization of estimate \eqref{Introduction Vasilevskaya} is the main result of the present paper.

\textit{Another approach} to proving operator error estimates in homogenization theory was suggested by V.~V.~Zhikov  \cite{Zh1} and developed by him together with S.~E.~Pastukhova  \cite{ZhPas,ZhPAs_parabol}. The authors of the method called it the ``\textit{shift method}'' or the ``\textit{modified method of the first approximation.}'' The method 
is based on introduction of  an additional parameter, specifically, a shift by a vector from~$\mathbb{R}^d$, on  careful analysis of the first order approximation to the solution, and  subsequent integration over the shift parameter. An important role is played by the Steklov smoothing. In papers \cite{Zh1,ZhPas}, analogs of estimates  \eqref{Introduction res A-eps principal} and \eqref{Introduction res A-eps H1-corrector} were proved for the operators of acoustics and elasticity theory. In  \cite{ZhPAs_parabol}, the shift method was applied  
to parabolic problems  and certain analogues of inequalities 
\eqref{Introduction exp-A-eps principal} and \eqref{Introduction exp-A-eps H1} were obtained. 
A summary of further results of V.~V.~Zhikov, S.~E.~Pastukhova, and their students can be found in the review~\cite{ZhPasUMN}.

So far, we discussed operator error estimates for the operators acting in the whole space $\mathbb{R}^d$. For completeness, note that it is more traditional for homogenization theory to study operators acting in a bounded domain $\mathcal{O}\subset \mathbb{R}^d$. Operator error estimates for such problems were studied by many authors. We highlight, in particular, the works \cite{MoV,Gr1,Gr2,ZhPas,KeLiS,SuMatematika,Xu3,GeSh}. 
(A more detailed survey of results on operator error estimates can be found in introduction to the paper \cite{MSu}.)

\subsection{The method of investigation} Consider the case $G=\mathbf{1}_n$. By using the scaling transformation, we reduce the proof of estimate \eqref{Introduction Main Result} to obtaining a similar estimate for the operator exponential $e^{-\mathcal{B}(\varepsilon) \varepsilon ^{-2} s}$, where $\mathcal{B}(\varepsilon)$ is the operator acting in $L_2(\mathbb{R}^d;\mathbb{C}^n)$ and given by the differential expression 
\begin{equation*}
\mathcal{B}(\varepsilon)=b(\mathbf{D})^*g (\mathbf{x}) b(\mathbf{D})
+\varepsilon\sum _{j=1}^d (a_j(\mathbf{x})D_j+D_ja_j(\mathbf{x})^*)+\varepsilon ^2\mathcal{Q}(\mathbf{x})+\varepsilon ^2 \lambda I.
\end{equation*}
Thus, from the study of the exponential of the operator with rapidly oscillating coefficients, we switch to study of the large time behavior  for the exponential of the operator with coefficients depending on $\mathbf{x}$ (not on $\mathbf{x}/\varepsilon$).

In accordance with~the Floquet-Bloch theory, the operator 
 $\mathcal{B}(\varepsilon)$ decomposes into the direct integral of operators  $\mathcal{B}(\mathbf{k},\varepsilon)$, acting in $L_2(\Omega ;\mathbb{C}^n)$ and depending on the parameter $\mathbf{k}\in\mathbb{R}^d$ called the quasi-momentum. The operator $\mathcal{B}(\mathbf{k},\varepsilon)$ is formally given by the expression
\begin{equation*}
\begin{split}
\mathcal{B}(\mathbf{k},\varepsilon) 
=b(\mathbf{D}+\mathbf{k})^*g(\mathbf{x}) b(\mathbf{D}+\mathbf{k})
&+\varepsilon\sum _{j=1}^d \big( a_j(\mathbf{x}) (D_j+k_j)+(D_j+k_j)a_j(\mathbf{x})^*\big)
\\
&+\varepsilon ^2\mathcal{Q}(\mathbf{x})+\varepsilon ^2\lambda I
\end{split}
\end{equation*}
with periodic boundary conditions. The spectrum of the operator $\mathcal{B}(\mathbf{k},\varepsilon)$ is discrete. 
In accordance with~\cite{Su10,Su14}, we distinguish the one-dimensional parameter $\tau =(\vert \mathbf{k}\vert ^2 +\varepsilon ^2)^{1/2}$ and study the family $\mathcal{B}(\mathbf{k},\varepsilon)$ by using the methods of analytic perturbation theory with respect to  $\tau$.

\subsection{The structure of the paper} The paper consists of Introduction and three chapters. In Chapter~1 (\S\S\ref{Section Predvaritelnye svedenia}--\ref{Section bordered exp abstract}), the abstract operator-theoretic scheme is  expounded. In Chapter~2 (\S\S\ref{Section Pervaritel'nye svedeniya}--\ref{Section f exp B(eps) f*}), periodic differential operators are studied. The approximation for the ``bordered'' operator exponential is obtained (see \S\ref{Section f exp B(eps) f*}). Chapter~3 (\S\S\ref{Section main results in operator terms}--\ref{Section Example}) is devoted to homogenization for solutions of the Cauchy problem for parabolic systems. In \S\ref{Section main results in operator terms}, by using the scaling transformation, from results of \S\ref{Section f exp B(eps) f*} we derive the main result of the paper, that is, estimate  \eqref{Introduction Main Result}. In \S\ref{Section Homogenization of parabolic systems}, 
results in operator terms are applied to homogenization of solutions for parabolic systems. In~\S\ref{Section Example}, we consider the scalar elliptic operator as an example.

\subsection{Notation}
Let $\mathfrak{H}$, $\mathfrak{H}_*$ be complex separable Hilbert spaces. The symbols $(\,\cdot\, ,\,\cdot\,)_\mathfrak{H}$ and $\Vert \,\cdot\,\Vert _\mathfrak{H}$ denote the scalar product and the norm in  $\mathfrak{H}$; the symbol $\Vert \,\cdot\,\Vert _{\mathfrak{H}\rightarrow\mathfrak{H}_*}$ denotes the norm of a linear continuous operator from $\mathfrak{H}$ to $\mathfrak{H}_*$.
 
By $\langle \,\cdot\, ,\,\cdot\,\rangle$, and $\vert \,\cdot\,\vert$ we denote the scalar product and the norm in  $\mathbb{C}^n$, respectively, by $\mathbf{1}_n$ we denote the unit $(n\times n)$-matrix. If $a$ is an $(m\times n)$\nobreakdash-matrix, then the symbol $\vert a\vert$ denotes the norm of the matrix $a$ as an operator form $\mathbb{C}^n$ to~$\mathbb{C}^m$. We use the notation
 $\mathbf{x}=(x_1,\dots , x_d)\in\mathbb{R}^d$, $iD_j=\partial _j =\partial /\partial x_j$, $j=1,\dots,d$, $\mathbf{D}=-i\nabla=(D_1,\dots ,D_d)$. The class $L_p$ of vector-valued functions in a domain  $\mathcal{O}\subset\mathbb{R}^d$ with values in $\mathbb{C}^n$ is denoted by $L_p(\mathcal{O};\mathbb{C}^n)$, $1\leqslant p\leqslant \infty$. The Sobolev class of $\mathbb{C}^n$-valued functions in a domain $\mathcal{O}\subset\mathbb{R}^d$ is denoted by $H^s(\mathcal{O};\mathbb{C}^n)$.  For $n=1$, we simply write  $L_p(\mathcal{O})$, $H^s(\mathcal{O})$, and so on, but (if this does
not lead to confusion) we use such short notation also for the spaces of vector-valued or
matrix-valued functions. The symbol $L_p((0,T);\mathfrak{H})$, $1\leqslant p\leqslant\infty$, means the  $L_p$-space of $\mathfrak{H}$-valued functions on the interval $(0,T)$.

Various constants in estimates are denoted by
 ${c, C,  \mathfrak{C}, \mathscr{C}}$ 
(possibly, with indices and marks).

\subsection{} The author is grateful to T.~A.~Suslina for the problem setting and carefully reading of the manuscript and to the administration of the Chebychev Laboratory for stimulating of completion of the present research.

\section*{Chapter 1. Abstract scheme}

\section{Preliminaries}

\label{Section Predvaritelnye svedenia}

The content of Sections \ref{Subsubsection operator pencils}--\ref{Subsection threshold approximations}  is borrowed from \cite{Su10,Su13}.

\subsection{Operators $X(t)$ and $A(t)$} 
\label{Subsubsection operator pencils}
Let $\mathfrak{H}$ and $\mathfrak{H}_*$ be complex separable Hilbert spaces. Let $X_0 \colon\mathfrak{H}\rightarrow \mathfrak{H}_*$ be a densely defined closed operator, and let the operator $X_1 \colon\mathfrak{H}\rightarrow\mathfrak{H}_*$ be bounded. On the domain $\mathrm{Dom}\,X(t)=\mathrm{Dom}\,X_0$, we define the operator
\begin{equation}
\label{X(t) in S1}
X(t):=X_0+tX_1\colon \mathfrak{H}\rightarrow\mathfrak{H}_*,\quad t\in\mathbb{R}.
\end{equation}
In $\mathfrak{H}$, consider a family of self-adjoint non-negative operators
\begin{equation}
\label{A(t)=}
A(t):=X(t)^*X(t),\quad t\in\mathbb{R},
\end{equation}
corresponding to closed  quadratic forms 
$$
\Vert X(t)u\Vert ^2_{\mathfrak{H}_*},\quad {u\in\mathrm{Dom}\,X_0},
$$ 
in $\mathfrak{H}$. 
Denote $A(0)=X_0^*X_0=:A_0$. Set
\begin{equation*}
\mathfrak{N}:=\mathrm{Ker}\,A_0=\mathrm{Ker}\,X_0,\quad \mathfrak{N}_*:=\mathrm{Ker}\,X_0^*.
\end{equation*}
Let $P$ be the orthogonal projection of the space $\mathfrak{H}$ onto $\mathfrak{N}$, and let $P_*$ be the orthogonal projection of $\mathfrak{H}_*$ onto $\mathfrak{N}_*$.

We impose the following condition.
\begin{condition}
\label{Condition lambda0}
The point $\lambda _0=0$ is isolated in the spectrum of the operator $A_0$, and
\begin{equation*}
0<n:=\mathrm{dim}\,\mathfrak{N}<\infty,\quad n\leqslant n_*:=\mathrm{dim}\,\mathfrak{N}_*\leqslant\infty .
\end{equation*}
\end{condition}

\textit{We denote by $d^0$  the distance between the point zero and the rest of the spectrum of the operator $A_0$.}

\subsection{The operators $Y(t)$ and $Y_2$}

\label{Subsection Y(t) i Y2}

Let $\widetilde{\mathfrak{H}}$ be yet another separable Hilbert space. Let $Y_0: \mathfrak{H}\rightarrow\widetilde{\mathfrak{H}}$ be a densely defined linear operator such that $\mathrm{Dom}\,X_0\subset \mathrm{Dom}\,Y_0$, and let $Y_1 :\mathfrak{H}\rightarrow\widetilde{\mathfrak{H}}$ be a bounded linear operator. Denote $Y(t):=Y_0+tY_1$, $\mathrm{Dom}\,Y(t)=\mathrm{Dom}\,Y_0$. The following condition is assumed to be satisfied. 

\begin{condition}
\label{Condition Y(t)}
There exists a constant $c_1>0$ such that
\begin{equation}
\label{Cond Y(t)u<=}
\Vert Y(t)u\Vert _{\widetilde{\mathfrak{H}}}\leqslant c_1\Vert X(t)u\Vert _{\mathfrak{H}_*},\quad u\in\mathrm{Dom}\,X_0,\quad t\in\mathbb{R}.
\end{equation}
\end{condition}

From estimate \eqref{Cond Y(t)u<=} with $t=0$ it follows that $\mathrm{Ker}\,X_0\subset\mathrm{Ker}\,Y_0$, i.~e., $Y_0P=0$.

Let $Y_2\colon\mathfrak{H}\rightarrow \widetilde{\mathfrak{H}}$ be a densely defined linear operator, and $\mathrm{Dom}\,X_0\subset\mathrm{Dom}\,Y_2$. The following condition is assumed to be satisfied.

\begin{condition}
\label{Condition Y2}
For any $\nu >0$ there exists a constant $C(\nu)>0$ such that
\begin{equation*}
\Vert Y_2 u\Vert ^2 _{\widetilde{\mathfrak{H}}}\leqslant \nu \Vert X(t)u\Vert ^2_{\mathfrak{H}_*}
+C(\nu)\Vert u\Vert ^2 _\mathfrak{H},\quad u\in\mathrm{Dom}\,X_0,\quad t\in\mathbb{R}.
\end{equation*}
\end{condition}

\subsection{The form $\mathfrak{q}$}

\label{Subsection form g abstract scheme}

In the space $\mathfrak{H}$, let $\mathfrak{q}[u,v]$ be a densely defined Hermitian sesquilinear form, and let $\mathrm{Dom}\,X_0\subset \mathrm{Dom}\,\mathfrak{q}$. We impose the following condition.

\begin{condition}
\label{Condition q}
$1^\circ$. There exist constants $c_2\geqslant 0$ and $c_3\geqslant 0$ such that
\begin{equation*}
\begin{split}
\vert \mathfrak{q}[u,v]\vert\leqslant (c_2\Vert X(t)u\Vert ^2_{\mathfrak{H}_*}+c_3\Vert u\Vert ^2_\mathfrak{H})^{1/2}(c_2\Vert X(t)v\Vert ^2_{\mathfrak{H}_*}+c_3\Vert v\Vert ^2_\mathfrak{H})^{1/2},
\\
u,v\in\mathrm{Dom}\,X_0,\quad t\in\mathbb{R}.
\end{split}
\end{equation*}

$2^\circ$. There exist constants $0<\kappa\leqslant 1$ and $c_0\in\mathbb{R}$ such that
\begin{equation*}
\mathfrak{q}[u,u]\geqslant -(1-\kappa)\Vert X(t)u\Vert ^2 _{\mathfrak{H}_*}-c_0\Vert u\Vert ^2_\mathfrak{H},\quad u\in\mathrm{Dom}\,X_0,\quad t\in\mathbb{R}.
\end{equation*}
\end{condition}

\subsection{The operator $B(t,\varepsilon)$}

\label{Subsection B(t,eps) abstract}

Let $\varepsilon\in (0,1]$. 
In $\mathfrak{H}$, we consider the Hermitian sesquilinear form
\begin{equation}
\label{b}
\begin{split}
\mathfrak{b}(t,\varepsilon)[u,v]
=(X(t)u,X(t)v)_{\mathfrak{H}_*}&+\varepsilon\left( (Y(t)u,Y_2v)_{\widetilde{\mathfrak{H}}}+(Y_2u,Y(t)v)_{\widetilde{\mathfrak{H}}}\right)
\\
&+\varepsilon ^2\mathfrak{q}[u,v],\quad u,v\in\mathrm{Dom}\,X_0.
\end{split}
\end{equation}

By using Conditions \ref{Condition Y(t)}, \ref{Condition Y2}, and \ref{Condition q}, it is easy to see that
\begin{align}
\mathfrak{b}(t,\varepsilon)[u,u]&\leqslant (2+c_1^2+c_2)\Vert X(t)u\Vert ^2 _{\mathfrak{H}_*}+(C(1)+c_3)\varepsilon ^2\Vert u\Vert ^2_\mathfrak{H},\quad u\in\mathrm{Dom}\,X_0,
\nonumber
\\
\label{frac b(t,e)>=}
\mathfrak{b}(t,\varepsilon)[u,u] &\geqslant \frac{\kappa}{2} \Vert X(t)u\Vert ^2_{\mathfrak{H}_*}-(c_0+c_4)\varepsilon ^2\Vert u\Vert ^2 _\mathfrak{H},\quad u\in\mathrm{Dom}\,X_0.
\end{align}
Here
\begin{equation}
\label{c4}
c_4=4\kappa ^{-1}c_1^2C(\nu)\quad\mbox{for}\quad\nu=\kappa ^2(16 c_1^2)^{-1}.
\end{equation}
A detailed proof of these inequalities can be found in \cite[Subsection~1.4]{Su10}. Thus, the form $\mathfrak{b}(t,\varepsilon)$ is closed and lower semi-bounded. 

The self-adjoint operator acting in the space $\mathfrak{H}$ and corresponding to the form \eqref{b} is denoted by $\mathfrak{B}(t,\varepsilon)$. Formally,
\begin{equation*}
\mathfrak{B}(t,\varepsilon)=A(t)+\varepsilon (Y_2^*Y(t)+Y(t)^*Y_2)+\varepsilon ^2Q.
\end{equation*}
Here $Q$ is the formal object that corresponds to the form $\mathfrak{q}$.

Let $Q_0\colon \mathfrak{H}\rightarrow\mathfrak{H}$ be a~bounded positive definite operator. 
 We add the term $\lambda Q_0$ to $Q$ in such a way that the operator
\begin{equation}
\label{B(t,eps) ne mathfrak}
B(t,\varepsilon):=\mathfrak{B}(t,\varepsilon)+\lambda \varepsilon ^2 Q_0 =A(t)+\varepsilon (Y_2^*Y(t)+Y(t)^*Y_2)+\varepsilon ^2(Q+\lambda  Q_0)
\end{equation}
corresponding to the form
\begin{equation}
\label{b(t,eps) forma which >0}
b(t,\varepsilon)[u,v]:=\mathfrak{b}(t,\varepsilon)[u,v]+\lambda\varepsilon^2 (Q_0 u,v)_\mathfrak{H},\quad u,v\in\mathrm{Dom}\,X_0,
\end{equation}
is positive definite. To guarantee this, we impose the following restriction on $\lambda$ 
\begin{equation}
\label{lambda condition}
\begin{split}
\lambda &>\Vert Q_0^{-1}\Vert (c_0+c_4),\quad\mbox{if }\lambda\geqslant 0,
\\
\lambda &>\Vert Q_0\Vert ^{-1}(c_0+c_4),\quad\mbox{if }\lambda<0\quad (\mbox{and }c_0+c_4<0).
\end{split}
\end{equation}
Condition \eqref{lambda condition} leads to the inequality
\begin{equation}
\label{lambda Q0 form >=}
\lambda (Q_0u,u)_{\mathfrak{H}}\geqslant (c_0+c_4+\beta)\Vert u\Vert ^2 _{\mathfrak{H}},\quad u\in\mathfrak{H},
\end{equation}
where $\beta >0$ is defined in terms of $\lambda$ by the rule
\begin{equation}
\label{beta condition}
\begin{split}
\beta &=\lambda \Vert Q_0^{-1}\Vert ^{-1}-c_0-c_4,\quad\mbox{if }\lambda\geqslant 0,
\\
\beta &=\lambda \Vert Q_0\Vert -c_0-c_4,\quad\mbox{if }\lambda<0 \quad (\mbox{and }c_0+c_4<0).
\end{split}
\end{equation}
From \eqref{frac b(t,e)>=} and \eqref{lambda Q0 form >=} it follows that a lower semi-bound for the form \eqref{b(t,eps) forma which >0} looks like this
\begin{equation}
\label{b (t,eps)>=}
b(t,\varepsilon)[u,u]\geqslant \frac{\kappa}{2}\Vert X(t)u\Vert ^2 _{\mathfrak{H}_*}+\beta \varepsilon ^2\Vert u\Vert ^2_\mathfrak{H},\quad u\in\mathrm{Dom}\,X_0.
\end{equation}
Thus, $B(t,\varepsilon)>0$. The operator $B(t,\varepsilon)$ is the main object of investigation in Chapter~1.

\subsection{Passage to the parameters $\tau$, $\vartheta$}

\label{Subsection to parameters tau and vartheta}

The operator family $B(t,\varepsilon)$ is an analytic with respect to the parameters $t$ and $\varepsilon$. If $t=\varepsilon=0$, then the operator $B(0,0)$ coincides with $A_0$ and, in accordance with Condition~\ref{Condition lambda0}, has an isolated eigenvalue $\lambda _0=0$ of multiplicity $n$. It seems natural to apply the analytic perturbation theory. But, for $n>1$, we deal with a multiple eigenvalue and a two-dimensional parameter. Analytic perturbation theory is not directly applicable. Thus we distinct a one-dimensional parameter $\tau =(t^2+\varepsilon ^2)^{1/2}$ and trace the dependence on the additional parameters $\vartheta _1=t\tau ^{-1}$ and $\vartheta _2 =\varepsilon\tau ^{-1}$. Note that the vector $\vartheta =(\vartheta _1,\vartheta _2)$ belongs to the unit circle. In that follows, the operator $B(t,\varepsilon)$ is denoted by $B(\tau;\vartheta)$ and the corresponding form $b(t,\varepsilon)$ is denoted by $b(\tau;\vartheta)$. According to \eqref{b} and \eqref{b(t,eps) forma which >0}, this form can be written as
\begin{equation*}
\begin{split}
b(\tau;\vartheta)[u,v]&=
(X_0u,X_0v)_{\mathfrak{H}_*}+\tau\vartheta _1 \left( (X_0u,X_1v)_{\mathfrak{H}_*}+(X_1u,X_0v)_{\mathfrak{H}_*}\right)
\\
&\quad+\tau ^2 \vartheta _1 ^2 (X_1 u,X_1v)_{\mathfrak{H}_*}+\tau\vartheta _2 \left( (Y_0u,Y_2v)_{\widetilde{\mathfrak{H}}}+(Y_2u,Y_0v)_{\widetilde{\mathfrak{H}}}\right)
\\
&\quad+\tau ^2\vartheta _1\vartheta _2\left((Y_1u,Y_2v)_{\widetilde{\mathfrak{H}}}+(Y_2u,Y_1v)_{\widetilde{\mathfrak{H}}}\right)
\\
&\quad+\tau ^2\vartheta _2 ^2 \left(\mathfrak{q}[u,v]+\lambda (Q_0u,v)_{\mathfrak{H}}\right),
\quad u,v\in\mathrm{Dom}\,X_0.
\end{split}
\end{equation*}
Formally,
\begin{equation*}
\begin{split}
B(\tau;\vartheta)=X_0^*X_0&+\tau\vartheta _1(X_0^*X_1+X_1^*X_0)+\tau ^2 \vartheta _1 ^2 X_1^*X_1
+\tau\vartheta _2(Y_2^*Y_0+Y_0^*Y_2)\\
&+\tau ^2 \vartheta _1\vartheta _2(Y_2^*Y_1+Y_1^*Y_2)
+\tau ^2\vartheta _2^2(Q+\lambda Q_0).
\end{split}
\end{equation*}

Let $F(\tau;\vartheta;\mu)$ be the spectral projector of the operator $B(\tau;\vartheta)$ corresponding to the interval $[0,\mu]$. Fix a constant $\delta\in (0,\kappa d^0/13)$ and choose a number $\tau _0>0$ such that
\begin{equation}
\label{tau 0 in abstract scheme}
\tau _0\leqslant \delta ^{1/2}\left(
(2+c_1^2+c_2)\Vert X_1 \Vert ^2 +C(1)+c_3+\vert \lambda\vert \Vert Q_0\Vert
\right)^{-1/2}.
\end{equation}
As was shown in \cite[Proposition 1.5]{Su10}, for $\vert\tau\vert \leqslant\tau _0$, we have 
\begin{equation*}
F(\tau;\vartheta;\delta)=F(\tau;\vartheta;3\delta),\quad\mathrm{rank}\,F(\tau;\vartheta
;\delta)=n,\quad\vert\tau\vert \leqslant\tau _0.
\end{equation*}
In that follows we often write $F(\tau;\vartheta)$ instead of $F(\tau;\vartheta ;\delta)$.

In Chapter~1, we trace the dependence of constants in estimates on the following  ``data'':
\begin{equation}
\label{data for abstract scheme}
\kappa ^{1/2}, \kappa ^{-1/2},  \delta, \delta ^{-1/2},  \tau _0,  c_1, c_2^{1/2}, c_3^{1/2}, C(1)^{1/2}, \vert \lambda\vert , \Vert X_1\Vert , \Vert Y_1\Vert , \Vert Q_0\Vert 
\end{equation}
(and also on the constant $\check{c}_*^{-1}$ introduced in Subsection~\ref{Subsection exp abstract approx} below). It is important for application of the results of the present chapter to differential operators that the constants in estimates (possibly after completion) depend polynomially on these quantities, and the coefficients of polynomials are positive numbers.

\subsection{The operators $Z$ and $\widetilde{Z}$}

\label{Subsection Z and tilde-Z-abstract}

In the present and next subsections, we define operators appearing in considerations inspired by the analytic perturbation theory.

Set $\mathcal{D}:=\mathrm{Dom}\,X_0\cap \mathfrak{N}^\perp$. Since the point $\lambda _0=0$ is isolated in the spectrum of $A_0$, the form $(X_0\phi,X_0\psi)_{\mathfrak{H}_*}$, $\phi,\psi\in\mathcal{D}$ defines a scalar product in $\mathcal{D}$, so $\mathcal{D}$ becomes a Hilbert space. 
Let $\omega\in\mathfrak{N}$. Consider an equation for the element $\phi\in\mathcal{D}$ (cf. \cite[Chapter~1, (1.7)]{BSu}):
\begin{equation}
\label{first eq. for Z}
X_0^*(X_0\phi +X_1\omega)=0,
\end{equation}
which is understood in a weak sense. In other words, $\phi\in\mathcal{D}$ satisfies the identity
\begin{equation}
\label{def Z phi tozd}
(X_0\phi ,X_0\zeta )_{\mathfrak{H}_*}=-(X_1\omega ,X_0\zeta)_{\mathfrak{H}_*},\quad \forall \zeta\in\mathcal{D}.
\end{equation}
Thus the right-hand side in \eqref{def Z phi tozd} is a continuous anti-linear functional on $\zeta\in\mathcal{D}$, from the Riesz theorem it follows that equation \eqref{first eq. for Z} has a unique solution $\phi (\omega)$. Obviously, $\Vert X_0\phi\Vert_{\mathfrak{H}_*}\leqslant\Vert X_1\omega \Vert _{\mathfrak{H}_*}$.  Define a~bounded operator $Z:\mathfrak{H}\rightarrow\mathfrak{H}$ by identities
\begin{equation*}
Z\omega =\phi (\omega),\quad\omega\in\mathfrak{N};\quad Zx=0,\quad x\in\mathfrak{N}^\perp .
\end{equation*}
Obviously,
\begin{equation}
\label{ZP=Z, PZ=0}
ZP=Z,\quad PZ=0,\quad PZ^*=Z^*,\quad Z^*P=0.
\end{equation}
It is easy to check that (see \cite[(1.21)]{Su13})
\begin{equation}
\label{Z<=}
\Vert Z\Vert _{\mathfrak{H}\rightarrow\mathfrak{H}}\leqslant\kappa^{1/2}(13\delta)^{-1/2}\Vert X_1\Vert.
\end{equation}

Similarly, for given $\omega\in\mathfrak{N}$ consider  the equation
\begin{equation}
\label{eq tilde Z abstract}
X_0^*X_0\psi+Y_0^*Y_2\omega =0
\end{equation}
for an element $\psi\in\mathcal{D}$. This equation is understood as the identity
\begin{equation*}
(X_0\psi ,X_0\zeta )_{\mathfrak{H}_*}=-(Y_2\omega,Y_0\zeta)_{\widetilde{\mathfrak{H}}},\quad\zeta\in\mathcal{D}.
\end{equation*}
By Condition~\ref{Condition Y(t)}, the right-hand side is a continuous anti-linear  functional of  $\zeta\in\mathcal{D}$. By the Riesz theorem, there exists a unique solution $\psi (\omega)$. Define a bounded operator $\widetilde{Z}:\mathfrak{H}\rightarrow \mathfrak{H}$ by the formulas
\begin{equation*}
\widetilde{Z}\omega = \psi (\omega),\quad\omega\in\mathfrak{N};\quad\widetilde{Z}x=0,\quad x\in \mathfrak{N}^\perp .
\end{equation*}
Note that $\widetilde{Z}$ maps $\mathfrak{N}$ onto $\mathfrak{N}^\perp$, and $\mathfrak{N}^\perp $ onto $\lbrace 0\rbrace$. So,
\begin{equation}
\label{tilde Z P}
\widetilde{Z}P=\widetilde{Z},\quad
P\widetilde{Z}=0,\quad P\widetilde{Z}^*=\widetilde{Z}^*, \quad \widetilde{Z}^*P=0.
\end{equation}
We need the estimate (see \cite[(1.25)]{Su13}):
\begin{equation}
\label{tilde Z<=}
\Vert \widetilde{Z}\Vert _{\mathfrak{H}\rightarrow \mathfrak{H}}\leqslant c_1 \kappa^{1/2}C(1)^{1/2}(13\delta )^{-1/2}.
\end{equation}

\subsection{The operators $R$ and $S$}

\label{Subsection R and S in abstract scheme}

Define an operator $R\colon \mathfrak{N}\rightarrow \mathfrak{N}_*$ (see \cite[Chapter~1, Subsection~1.2]{BSu}) as follows
$ R\omega =X_0\psi(\omega)+X_1\omega\in\mathfrak{N}_*$. Another definition of the operator $R$ is given by the formula $R=P_*X_1\vert _\mathfrak{N}$.

By the \textit{spectral germ} of the operator family \eqref{A(t)=} for ${t\!=\!0}$~we call\! (see\! \cite[\!Chapter~1,\! Subsection~1.3]{BSu})\! {the self-adjoint operator~${S\!=\!R^*R \colon\! \mathfrak{N}\!\rightarrow \!\mathfrak{N}}$,} which also satisfies the identity $${S=PX_1^*P_*X_1\vert _\mathfrak{N}}.$$ 

\subsection{Analytic branches of eigenvalues and eigenvectors of the operator $B(\tau;\vartheta)$}

According to the~analytic perturbation theory, (see \cite{K}), for $\vert \tau\vert\leqslant \tau_0$ there exist real-analytic functions $\lambda _l(\tau;\vartheta)$ and real-analytic (in $\tau$) $\mathfrak{H}$-valued functions $\phi _l(\tau ;\vartheta)$ such that
\begin{equation}
\label{B phi =lambda phi}
B(\tau;\vartheta)\phi _l(\tau;\vartheta)=\lambda _l(\tau;\vartheta)\phi _l(\tau;\vartheta),\quad l=1,\dots,n,\quad\vert \tau\vert \leqslant \tau _0,
\end{equation}
moreover, $\phi _l(\tau;\vartheta)$, $l=1,\dots,n$, form an orthogonal basis in the eigenspace $ F(\tau;\vartheta)\mathfrak{H}$. For sufficiently small $\tau_* (\vartheta)$ ($\leqslant \tau _0$) and $\vert \tau\vert \leqslant \tau_* (\vartheta)$ we have the following convergent power series expansions 
\begin{align}
\label{lambda_l(t)= series}
&\lambda _l(\tau;\vartheta)=\gamma _l (\vartheta) \tau ^2 +\mu _l (\vartheta) \tau ^3+\dots ,\quad \gamma _l(\vartheta)\geqslant 0,
\quad l=1,\dots,n;\\
&\phi _l(\tau;\vartheta)=\omega _l (\vartheta) +\tau \phi _l^{(1)}(\vartheta)+\tau ^2\phi _l^{(2)}(\vartheta)+\dots,\quad l=1,\dots,n.
\nonumber
\end{align}
The elements $\omega _l (\vartheta)$, $l=1,\dots ,n$, form an orthonormal basis in $\mathfrak{N}$.

In \cite[(1.32), (1.33)]{Su10}, it was shown that the elements
\begin{equation*}
\widetilde{\omega}_l(\vartheta):=\phi _l^{(1)}(\vartheta)-\vartheta _1 Z\omega _l(\vartheta)-\vartheta _2 \widetilde{Z}\omega _l (\vartheta),\quad l=1,\dots,n,
\end{equation*}
satisfy the identities
\begin{equation}
\label{omega i tilde omega tozd}
\left(\widetilde{\omega}_k(\vartheta),\omega _l(\vartheta)\right)_\mathfrak{H}
+\left(\omega _k(\vartheta),\widetilde{\omega}_l(\vartheta)\right)_\mathfrak{H}=0,\quad k,l=1,\dots,n.
\end{equation}

In accordance with \cite[Subsection~1.8]{Su10}, the {\it spectral germ of the operator family $B(\tau;\vartheta)$ for $\tau =0$} is the operator acting in $\mathfrak{N}$ and given by the expression
\begin{equation}
\label{S(theta)=}
\begin{split}
S(\vartheta)\!=\!\vartheta _1^2S&\!+\!\vartheta_1\vartheta_2(
-(X_0Z)^*X_0\widetilde{Z}\!-\!(X_0\widetilde{Z})^*X_0Z\!+\!P(Y_2^*Y_1\!+\!Y_1^*Y_2)
)\Big\vert _{\mathfrak{N}}
\\
&+\vartheta _2^2\left(
-(X_0\widetilde{Z})^*X_0\widetilde{Z}\vert _\mathfrak{N}+Q_\mathfrak{N}+\lambda Q_{0\mathfrak{N}}\right).
\end{split}
\end{equation}
Here $Q_\mathfrak{N}$ is the self-adjoint operator in $\mathfrak{N}$, generated by the form $\mathfrak{q}[\omega,\omega]$ on $\omega\in\mathfrak{N}$, and $Q_{0\mathfrak{N}}:=PQ_0\vert _\mathfrak{N}$.

In \cite[Proposition~1.6]{Su10}, it was proved that the numbers $\gamma _l(\vartheta)$ and elements $\omega _l(\vartheta)$, $l=1,\dots,n$, are, respectively, eigenvalues and eigenvectors for the operator $S(\vartheta)$:
\begin{equation}
\label{S omega _l=}
S(\vartheta)\omega_l(\vartheta)=\gamma_l(\vartheta)\omega_l(\vartheta),\quad l=1,\dots, n.
\end{equation}
The quantities $\gamma _l(\vartheta)$ and $\omega _l(\vartheta)$, $l=1,\dots ,n$, are called the \textit{threshold characteristics at the bottom of the spectrum} for the operator family $B(\tau;\vartheta)$.

\subsection{Threshold approximations} 

\label{Subsection threshold approximations}

As was shown in \cite[Theorem~3.1]{Su13},
\begin{equation}
\label{F-P}
F(\tau;\vartheta)-P=\Phi(\tau;\vartheta),\quad \Vert \Phi (\tau;\vartheta)\Vert _{\mathfrak{H}\rightarrow\mathfrak{H}} \leqslant C_1 \vert \tau\vert,\quad \vert \tau\vert\leqslant \tau_0 .
\end{equation}
Along with \eqref{F-P}, we need a sharper approximation for the spectral projection obtained in \cite[Theorem~3.1]{Su13}:
\begin{align}
\label{F(t)=P+tF_1+F_2(t)}
F(\tau;\vartheta )&=P+\tau F_1 (\vartheta)+F_2(\tau;\vartheta),\quad \Vert F_2(\tau;\vartheta)\Vert\leqslant C_2 \tau ^2,\quad \vert \tau\vert\leqslant \tau _0.
\end{align}
By \cite[(1.41)]{Su13}, the operator $F_1(\vartheta)$ has the form
\begin{equation}
\label{F_1=}
F_1(\vartheta)=\vartheta _1 (Z+Z^*)+\vartheta _2(\widetilde{Z}+\widetilde{Z}^*).
\end{equation}
From \eqref{ZP=Z, PZ=0}, \eqref{tilde Z P}, and \eqref{F_1=} it follows that
\begin{equation}
\label{F1P=ZP}
F_1(\vartheta) P=\vartheta_1 Z +\vartheta _2\widetilde{Z},
\quad
PF_1 (\vartheta)=\vartheta _1Z^*+\vartheta _2\widetilde{Z}^*.
\end{equation}

From \eqref{Z<=}, \eqref{tilde Z<=}, and \eqref{F_1=}, we derive the estimate
\begin{equation}
\label{F_1<=}
\Vert F_1(\vartheta)\Vert _{\mathfrak{H}\rightarrow\mathfrak{H}}\leqslant C_{F_1}=
2\kappa ^{1/2}(13\delta)^{-1/2}(\Vert X_1\Vert +c_1C(1)^{1/2}).
\end{equation}

We need an approximation obtained in \cite[Theorem~2.2]{Su10}:
\begin{equation}
\begin{split}
\label{BF-SP<=}
B(\tau;\vartheta)F(\tau;\vartheta)-\tau ^2S(\vartheta)P&=\Phi _1(\tau;\vartheta),
\\
\Vert \Phi _1(\tau;\vartheta)\Vert _{\mathfrak{H}\rightarrow\mathfrak{H}}&\leqslant C_3\vert \tau\vert ^3,\quad\vert\tau\vert \leqslant \tau _0.
\end{split}
\end{equation}

A refinement of approximation \eqref{BF-SP<=} was achieved in \cite[Theorem~3.3]{Su13}:
\begin{align}
\label{BF =tau^2S+tau3K+Phi2}
B(\tau;\vartheta)F(\tau;\vartheta)&=\tau ^2S(\vartheta)P+\tau ^3 K(\vartheta)+\Phi_2(\tau ;\vartheta),
\\
\label{Phi 2<=}
\Vert \Phi _2(\tau ;\vartheta)\Vert &\leqslant C_4\tau ^4,\quad \vert \tau\vert \leqslant \tau _0.
\end{align}

By \cite[(3.18)--(3.20)]{Su13}, we have
$
K(\vartheta)=K_0(\vartheta)+N(\vartheta)$, 
where
\begin{align}
\begin{split}
K_0(\vartheta)&=\sum _{l=1}^n \gamma _l(\vartheta)\Bigl(
(\,\cdot\, ,\omega _l(\vartheta))_\mathfrak{H}(\vartheta _1 Z\omega _l(\vartheta)+\vartheta _2\widetilde{Z}\omega _l(\vartheta))
\\
&\quad+(\,\cdot\, ,\vartheta _1 Z\omega _l(\vartheta)+\vartheta _2\widetilde{Z}\omega _l(\vartheta))_\mathfrak{H}\omega _l(\vartheta)\Bigr),
\end{split}
\nonumber
\\
\label{N=N0+N*}
N(\vartheta)&=N_0(\vartheta)+N_*(\vartheta),
\\
N_0(\vartheta)&=\sum _{l=1}^n \mu _l(\vartheta)(\,\cdot\,,\omega _l(\vartheta))_\mathfrak{H}\omega _l(\vartheta),
\nonumber
\\
N_*(\vartheta)&=\sum _{l=1}^n \gamma _l(\vartheta)\bigl((\,\cdot\, ,\omega _l(\vartheta))_\mathfrak{H}\widetilde{\omega}_l(\vartheta)+(\,\cdot\, ,\widetilde{\omega}_l(\vartheta))_\mathfrak{H}\omega _l(\vartheta)\bigr).
\nonumber
\end{align}
In the basis $\lbrace \omega _l(\vartheta)\rbrace _{l=1}^n$, the operators $N(\vartheta)$, $N_0(\vartheta)$, and $N_*(\vartheta)$ (after restriction to~$\mathfrak{N}$) are represented by the $(n\times n)$-matrices. The operator $N_0(\vartheta)$ is diagonal: $(N_0(\vartheta)\omega _j (\vartheta),\omega _l(\vartheta))=\mu _j (\vartheta)\delta _{jl}$, $j,l=1,\dots,n$. By \eqref{omega i tilde omega tozd}, the matrix elements of the operator $N_*(\vartheta)$ have the form
\begin{equation}
\label{N* matrix elements}
\begin{split}
(N_*(\vartheta)\omega _j (\vartheta),\omega _l(\vartheta))&=\gamma _l(\vartheta)\left(\omega _j(\vartheta),\widetilde{\omega}_l(\vartheta)\right)
+\gamma _j(\vartheta)\left(\widetilde{\omega}_j(\vartheta),\omega _l(\vartheta)\right)
\\
&=\left(\gamma _j(\vartheta)-\gamma _l(\vartheta)\right)\left(\widetilde{\omega}_j(\vartheta),\omega _l(\vartheta)\right),
\end{split}
\end{equation}
$j,l=1,\dots,n$. Thus, the diagonal elements $N_*(\vartheta)$ are equal to zero. Moreover, $$\left(N_*(\vartheta)\omega _j(\vartheta),\omega _l(\vartheta)\right)=0,$$ if $\gamma _j(\vartheta)=\gamma _l(\vartheta)$. In the case when $n=1$, we have $N_*(\vartheta)=0$, i.~e., $N(\vartheta)=N_0(\vartheta)$.

Invariant representations for the operators $K_0(\vartheta)$ and $N(\vartheta)$ 
were found in \cite[(3.21), (3.30)--(3.34)]{Su13}:
\begin{align}
\label{K0(theta)=}
K_0(\vartheta)&=\vartheta _1\left(ZS(\vartheta)P+S(\vartheta)PZ^*\right)
+\vartheta _2\left(\widetilde{Z}S(\vartheta)P+S(\vartheta)P\widetilde{Z}^*\right),
\\
\label{N(theta)=}
N(\vartheta)&=\vartheta _1^3N_{11}+\vartheta _1^2\vartheta _2N_{12}+\vartheta _1\vartheta _2^2N_{21}+\vartheta _2^3N_{22}.
\end{align}
Here
\begin{align}
\label{N11 abstract}
N_{11}&=(X_1Z)^*RP+(RP)^*X_1Z,
\\
\label{N12 abstract}
\begin{split}
N_{12}&=(X_1\widetilde{Z})^*RP+(RP)^*X_1\widetilde{Z}+(X_1Z)^*X_0\widetilde{Z}
\\
&\quad+(X_0\widetilde{Z})^*X_1Z+(Y_2Z)^*Y_0Z+(Y_0Z)^*Y_2Z
\\
&\quad+(Y_2Z)^*Y_1P+(Y_1P)^*Y_2Z+(Y_2P)^*Y_1Z+(Y_1Z)^*Y_2P,
\end{split}
\\
\label{N21 abstract}
\begin{split}
N_{21}&=(X_0\widetilde{Z})^*X_1\widetilde{Z}+(X_1\widetilde{Z})^*X_0\widetilde{Z}+(Y_2Z)^*Y_0\widetilde{Z}
\\
&\quad+(Y_0\widetilde{Z})^*Y_2Z+(Y_2\widetilde{Z})^*Y_0Z+(Y_0Z)^*Y_2\widetilde{Z}+(Y_2\widetilde{Z})^*Y_1P
\\
&\quad+(Y_1P)^*Y_2\widetilde{Z}+(Y_1\widetilde{Z})^*Y_2P+(Y_2P)^*Y_1\widetilde{Z}
\\
&\quad+Z^*QP+PQZ+\lambda ({Z}^*Q_0P+PQ_0Z),
\end{split}
\\
\label{N22 abstract}
\begin{split}
N_{22}&=(Y_0\widetilde{Z})^*Y_2\widetilde{Z}+(Y_2\widetilde{Z})^*Y_0\widetilde{Z}+\widetilde{Z}^*QP+PQ\widetilde{Z}
\\
&\quad+\lambda(\widetilde{Z}^*Q_0P+PQ_0\widetilde{Z}).
\end{split}
\end{align}
We clarify that, in \eqref{N21 abstract}, the formal notation $Z^*QP+PQZ$ is used for the bounded self-adjoint operator acting in $\mathfrak{H}$ and corresponding to the form $\mathfrak{q}[Pu,Zu]+\mathfrak{q}[Zu,Pu]$, $u\in\mathfrak{H}$. 
Similarly, in \eqref{N22 abstract}, $\widetilde{Z}^*QP+PQ\widetilde{Z}$ is understood as the~bounded self-adjoint operator in $\mathfrak{H}$ corresponding to the form $\mathfrak{q}[Pu,\widetilde{Z}u]+\mathfrak{q}[\widetilde{Z}u,Pu]$, $u\in\mathfrak{H}$.

By \eqref{ZP=Z, PZ=0}, \eqref{tilde Z P}, and \eqref{K0(theta)=}--\eqref{N22 abstract}, for the operator $
K(\vartheta)=K_0(\vartheta)+N(\vartheta)$ we have 
\begin{equation}
\label{PKP}
PK(\vartheta)P=PN(\vartheta)P=N(\vartheta).
\end{equation}

We need the following estimates obtained in \cite[(3.35), (3.36), (3.44)--(3.49)]{Su13}:
\begin{align}
\label{N(theta)<=}
&\Vert N(\vartheta)\Vert _{\mathfrak{H}\rightarrow\mathfrak{H}}\leqslant C_N ,
\\
\label{K<=}
&\Vert K(\vartheta)\Vert _{\mathfrak{H}\rightarrow\mathfrak{H}}\leqslant C_K.
\end{align}

\begin{remark}
\label{Remark constants Sec.1}
The constants $C_1$, $C_2$, $C_3$, $C_4$, and $C_{F_1}$ introduced above depend on the~data  \eqref{data for abstract scheme} polynomially. The corresponding polynomials have positive numerical  coefficients. The constants $C_N$ and $C_K$ can be estimated in terms of polynomials with numeric coefficients and variables \eqref{data for abstract scheme}. The dependence of constants $C_1$, $C_2$, $C_3$, $C_4$,  $C_N$, and $C_K$ on the data  \eqref{data for abstract scheme} was traced in \cite{Su10,Su13}.
\end{remark}

\section{Approximation of the operator $e^{-B(t,\varepsilon)s}$ with the corrector taken into account}

\subsection{Approximation of the operator $e^{-B(\tau;\vartheta)s}$}

\label{Subsection exp abstract approx}

{\it Assume that for some $c_*>0$ we have}
\begin{equation}
\label{A(t)>= abstract}
A(t)\geqslant c_*t^2 I,\quad c_*>0,\quad\vert t\vert \leqslant \tau _0.
\end{equation}
Then, by \eqref{b (t,eps)>=},
\begin{equation}
\label{B(tau,theta)>=}
B(\tau;\vartheta)\geqslant \check{c}_*\tau ^2 I,\quad \vert \tau\vert \leqslant\tau _0;
\quad \check{c}_*=\frac{1}{2}\min\lbrace\kappa c_*;2\beta\rbrace.
\end{equation}
Together with \eqref{B phi =lambda phi} this implies that the eigenvalues $\lambda _l(\tau;\vartheta)$ of the operator $B(\tau;\vartheta)$ satisfy the estimates
\begin{equation}
\label{lambda l (tau,theta)>=}
\lambda _l(\tau;\vartheta)\geqslant \check{c}_*\tau ^2,\quad l=1,\dots,n,\quad \vert \tau\vert\leqslant\tau _0.
\end{equation}
Bringing \eqref{lambda_l(t)= series} and \eqref{lambda l (tau,theta)>=} together, we see that $\gamma _l(\vartheta)\geqslant \check{c}_*$, $l=1,\dots,n$, and so (see \eqref{S omega _l=})
\begin{equation}
\label{S(theta)>=}
S(\vartheta)\geqslant\check{c}_* I_\mathfrak{N}.
\end{equation}

\textit{Let} $\vert \tau\vert\leqslant \tau _0$. We rewrite the operator $e^{-B(\tau;\vartheta)s}$, $s\geqslant 0$, as
\begin{equation}
\label{expB=}
e^{-B(\tau;\vartheta)s}=e^{-B(\tau;\vartheta)s}F(\tau;\vartheta)^\perp +e^{-B(\tau;\vartheta)s} F(\tau;\vartheta).
\end{equation}
Let $s>0$. By \eqref{B(tau,theta)>=}, the first summand on the right-hand side \eqref{expB=} satisfies the estimate
\begin{equation}
\label{expB first summand est}
\Vert e^{-B(\tau;\vartheta)s}F(\tau;\vartheta)^\perp \Vert _{\mathfrak{H}\rightarrow\mathfrak{H}}
\leqslant e^{-\delta s/2}e^{-\check{c}_*\tau ^2 s/2}\leqslant 2(\delta s)^{-1}e^{-\check{c}_*\tau ^2 s/2}.
\end{equation}
The second summand can be represented as
\begin{equation}
\label{expB second summand tozd}
e^{-B(\tau;\vartheta)s} F(\tau;\vartheta)=Pe^{-B(\tau;\vartheta)s} F(\tau;\vartheta)+P^\perp e^{-B(\tau;\vartheta)s} F(\tau;\vartheta).
\end{equation}
By \eqref{F(t)=P+tF_1+F_2(t)},
\begin{equation*}
\begin{split}
&P^\perp e^{-B(\tau;\vartheta)s} F(\tau;\vartheta)
=(I-P)e^{-B(\tau;\vartheta)s} F(\tau;\vartheta)
\\
&=(\tau F_1(\vartheta)+F_2(\tau;\vartheta))e^{-B(\tau;\vartheta)s} F(\tau;\vartheta)
\\
&=(\tau F_1(\vartheta)+F_2(\tau;\vartheta))
\left( (e^{-B(\tau;\vartheta)s} F(\tau;\vartheta)-e^{-\tau ^2 S(\vartheta)s}P)+e^{-\tau ^2 S(\vartheta)s}P\right).
\end{split}
\end{equation*}
Set
\begin{equation}
\label{Pi(tau,theta,s)=}
\Pi (\tau;\vartheta ;s):=e^{-B(\tau;\vartheta)s}F(\tau;\vartheta)-e^{-\tau ^2 S(\vartheta)s}P.
\end{equation}
Then
\begin{equation}
\label{P perp exp B F=}
\begin{split}
P^\perp e^{-B(\tau;\vartheta)s} F(\tau;\vartheta)= \tau F_1(\vartheta)e^{-\tau ^2 S(\vartheta)s}P
&+\tau F_1(\vartheta)\Pi (\tau;\vartheta ;s)
\\
&+F_2(\tau;\vartheta)e^{-B(\tau;\vartheta)s} F(\tau;\vartheta).
\end{split}
\end{equation}
We need the estimate obtained in \cite[Subsection~2.2]{M_AA}:
\begin{equation}
\label{Pi(tau,theta,s)<=}
\Vert \Pi (\tau;\vartheta ;s)\Vert _{\mathfrak{H}\rightarrow\mathfrak{H}}
\leqslant
(2C_1\vert\tau\vert +C_3\vert\tau\vert ^3 s)e^{-\check{c}_*\tau ^2 s}.
\end{equation}
Together with \eqref{F_1<=}, this implies
\begin{equation*}
\Vert \tau F_1(\vartheta)\Pi (\tau;\vartheta ;s)\Vert _{\mathfrak{H}\rightarrow\mathfrak{H}}
\leqslant C_{F_1}(2C_1\tau ^2+C_3\tau ^4s)e^{-\check{c}_*\tau ^2 s}.
\end{equation*}
Introducing the parameter $\alpha =\tau ^2 s$, we arrive at the inequality
\begin{equation}
\label{F1Pi<=}
\Vert \tau F_1(\vartheta)\Pi (\tau;\vartheta ;s)\Vert _{\mathfrak{H}\rightarrow\mathfrak{H}}
\leqslant C_{F_1}(2C_1\alpha +C_3\alpha ^2)e^{-\check{c}_*\alpha}s^{-1}.
\end{equation}
The third summand on the right-hand side of \eqref{P perp exp B F=} is estimated with the help of \eqref{F(t)=P+tF_1+F_2(t)} and \eqref{B(tau,theta)>=}:
\begin{equation}
\label{F2expBF<=}
\Vert F_2(\tau;\vartheta)e^{-B(\tau;\vartheta)s} F(\tau;\vartheta)\Vert  _{\mathfrak{H}\rightarrow\mathfrak{H}}
\leqslant 
C_2\tau ^2 e^{-\check{c}_*\tau ^2 s}
=C_2\alpha e^{-\check{c}_*\alpha}s^{-1}.
\end{equation}
Bringing  \eqref{P perp exp B F=}, \eqref{F1Pi<=}, and \eqref{F2expBF<=} together, we obtain
\begin{equation}
\label{P perp exp B F<=}
\Vert P^\perp e^{-B(\tau;\vartheta)s} F(\tau;\vartheta)-\tau F_1(\vartheta)e^{-\tau ^2 S(\vartheta)s}P
\Vert  _{\mathfrak{H}\rightarrow\mathfrak{H}}
\leqslant 
\phi (\alpha) s^{-1},\quad s>0,
\end{equation}
where
\begin{equation}
\label{phi(alpha)=}
\phi (\alpha)=C_{F_1}(2C_1\alpha +C_3\alpha ^2)e^{-\check{c}_*\alpha}+C_2\alpha e^{-\check{c}_*\alpha}.
\end{equation}

Now, consider the first summand on the right-hand side of \eqref{expB second summand tozd}:
\begin{equation}
\label{P exp B F=}
P e^{-B(\tau;\vartheta)s} F(\tau;\vartheta)=e^{-\tau ^2 S(\vartheta) s}P+\Upsilon (\tau;\vartheta ;s),
\end{equation}
where
$
\Upsilon (\tau;\vartheta ;s):=P e^{-B(\tau;\vartheta)s} F(\tau;\vartheta)-e^{-\tau ^2 S(\vartheta) s}P$. 
In \cite[Subsection~2.2]{M_AA}, it was shown that 
\begin{align*}
\Upsilon (\tau;\vartheta ;s)&=e^{-\tau ^2 S(\vartheta) s}P\Phi (\tau;\vartheta)
-\mathcal{J}(\tau;\vartheta ;s),
\\
\mathcal{J}(\tau;\vartheta ;s):&=
\int\limits_0^s e^{-\tau ^2 S(\vartheta)P(s-\widetilde{s})}P\Phi _1(\tau ;\vartheta)e^{-B(\tau;\vartheta)\widetilde{s}}F(\tau;\vartheta)\,d\widetilde{s}.
\end{align*}
Here the operators $\Phi (\tau;\vartheta)$ and $\Phi _1(\tau;\vartheta)$ are defined in \eqref{F-P} and \eqref{BF-SP<=}. Taking~\eqref{F-P} and \eqref{F(t)=P+tF_1+F_2(t)} into account, we represent the operator $\Upsilon (\tau;\vartheta ;s)$ as
\begin{equation}
\label{Upsilon tozd}
\begin{split}
\Upsilon (\tau;\vartheta ;s)=e^{-\tau ^2 S(\vartheta)Ps}P\tau F_1(\vartheta)
+e^{-\tau ^2 S(\vartheta)Ps}P F_2(\tau;\vartheta)
-\mathcal{J}(\tau;\vartheta ;s).
\end{split}
\end{equation}
The second summand on the right-hand side of \eqref{Upsilon tozd} can be estimated with the help of \eqref{F(t)=P+tF_1+F_2(t)} and \eqref{S(theta)>=}:
\begin{equation}
\Vert e^{-\tau ^2 S(\vartheta)Ps}P F_2(\tau;\vartheta)\Vert _{\mathfrak{H}\rightarrow\mathfrak{H}}\leqslant
C_2\tau ^2 e^{-\tau ^2\check{c}_*s}=C_2\alpha e^{-\check{c}_*\alpha}s^{-1}.
\end{equation}
By \eqref{BF-SP<=}, \eqref{BF =tau^2S+tau3K+Phi2}, \eqref{PKP}, and \eqref{Pi(tau,theta,s)=}, 
the third summand on the right-hand side of \eqref{Upsilon tozd} can be rewritten as
\begin{equation}
\label{tr summ in Upsilon tozd}
\begin{split}
&\mathcal{J}(\tau;\vartheta ;s)
\\
&\!=\!\!\int\limits_0^s \!e^{-\tau ^2S(\vartheta)P(s\!-\!\widetilde{s})}P (\tau ^3 K(\vartheta)\!+\!\Phi_2(\tau;\vartheta))e^{-B(\tau;\vartheta)\widetilde{s}}F(\tau;\vartheta)\,d\widetilde{s}
\\
&\!=\!\!\int\limits_0^s \!e^{\!-\tau ^2S(\vartheta)P(s\!-\!\widetilde{s})}P (\tau ^3 K(\vartheta)\!+\!\Phi_2(\tau;\vartheta))
\!\Big(\Pi (\tau;\vartheta;\widetilde{s})\!+\!e^{\!-\tau ^2 S(\vartheta)\widetilde{s}}P\Big)\,d\widetilde{s}
\\
&\!=\!\mathcal{J}_1(\tau;\vartheta ;s)\!+\!\mathcal{J}_2(\tau;\vartheta ;s)\!+\!\mathcal{J}_3(\tau;\vartheta ;s).
\end{split}
\end{equation}
Here
\begin{align*}
\mathcal{J}_1(\tau;\vartheta ;s):&=\tau ^3 \int\limits_0^s  e^{-\tau ^2S(\vartheta)P(s-\widetilde{s})}P  N(\vartheta)Pe^{-\tau ^2 S(\vartheta)\widetilde{s}}P\,d\widetilde{s},
\\
\mathcal{J}_2(\tau;\vartheta ;s):&=\tau ^3\int\limits_0^s  e^{-\tau ^2S(\vartheta)P(s-\widetilde{s})}P K(\vartheta)\Pi (\tau;\vartheta;\widetilde{s})\,d\widetilde{s},
\\
\mathcal{J}_3(\tau;\vartheta ;s):&=\int\limits_0^s  e^{-\tau ^2S(\vartheta)P(s-\widetilde{s})}P \Phi _2(\tau;\vartheta)e^{-B(\tau;\vartheta)\widetilde{s}}F(\tau;\vartheta)\,d\widetilde{s}.
\end{align*}
We estimate the term $\mathcal{J}_2(\tau;\vartheta ;s)$ with the help of \eqref{K<=}, \eqref{S(theta)>=}, and \eqref{Pi(tau,theta,s)<=}:
\begin{equation}
\begin{split}
\Bigl\Vert \mathcal{J}_2(\tau;\vartheta ;s)\Bigr\Vert _{\mathfrak{H}\rightarrow\mathfrak{H}}
&\leqslant
\vert \tau \vert ^3\int\limits_0^s e^{-\check{c}_*\tau ^2 (s-\widetilde{s})}C_K(2C_1\vert \tau\vert +C_3\vert \tau\vert ^3 \widetilde{s}) e^{-\check{c}_*\tau ^2 \widetilde{s}}\,d\widetilde{s}
\\
&=
\vert \tau\vert ^3 e^{-\check{c}_*\tau ^2 s}C_K \int\limits_0^s (2C_1\vert \tau\vert +C_3\vert \tau\vert ^3 \widetilde{s})\,d\widetilde{s}
\\
&=e^{-\check{c}_*\tau ^2 s}C_K (2C_1\tau ^4 s +2^{-1}C_3\tau ^6 s^2)
\\
&= C_Ke^{-\check{c}_*\alpha}(2C_1\alpha ^2+2^{-1}C_3\alpha ^3)s^{-1}.
\end{split}
\end{equation}
Using \eqref{Phi 2<=},  \eqref{B(tau,theta)>=}, and \eqref{S(theta)>=}, we estimate
$\mathcal{J}_3(\tau;\vartheta ;s)$:
\begin{equation}
\label{tr summ in tr summ in Upsilon tozd}
\begin{split}
\Bigl\Vert \mathcal{J}_3(\tau;\vartheta ;s)\Bigr\Vert _{\mathfrak{H}\rightarrow\mathfrak{H}}
&\leqslant 
\int\limits_0^s e^{-\check{c}_*\tau ^2 (s-\widetilde{s})}C_4\tau ^4 e^{-\check{c}_*\tau ^2\widetilde{s}}\,d\widetilde{s}
\\
&=e^{-\check{c}_*\tau ^2 s}C_4\tau ^4 s =C_4e^{-\check{c}_*\alpha}\alpha ^2 s^{-1}.
\end{split}
\end{equation}

We summarise the results. By \eqref{P exp B F=}--\eqref{tr summ in tr summ in Upsilon tozd},
\begin{equation}
\label{Neraw}
\begin{split}
\Big\Vert  Pe^{-B(\tau;\vartheta)s}\!F(\tau;\vartheta)\!-\!e^{\!-\tau ^2 S(\vartheta)s}P\!\!-\!\tau e^{\!-\tau ^2 S(\vartheta)s}PF_1(\vartheta)\!+\!\tau ^3\!\mathcal{M}(\tau;\vartheta ;s)\Big\Vert _{\mathfrak{H}\!\rightarrow\!\mathfrak{H}}
\\
\!\leqslant\!\!
\big( C_2\alpha \!+\!2C_1C_K\alpha ^2\!+\!2^{-1}C_3C_K\alpha ^3\!+\!C_4\alpha ^2\big) e^{-\check{c}_*\alpha}s^{-1},
\end{split}
\end{equation}
where
\begin{equation}
\label{M=}
\mathcal{M}(\tau;\vartheta ;s):=\int\limits_0^s e^{-\tau ^2 S(\vartheta)P(s-\widetilde{s})}PN(\vartheta)Pe^{-\tau ^2 S(\vartheta)P\widetilde{s}}\,d\widetilde{s}.
\end{equation}
Now from  \eqref{F1P=ZP}, \eqref{expB=}--\eqref{expB second summand tozd}, \eqref{P perp exp B F<=}, \eqref{phi(alpha)=}, \eqref{Neraw}, and \eqref{M=}, we derive the identity
\begin{equation}
\label{Th exp abstract theta tau}
\begin{split}
e^{-B(\tau;\vartheta)s}&=e^{-\tau ^2 S(\vartheta)s}P
+\tau (\vartheta _1 Z+\vartheta _2\widetilde{Z})e^{-\tau ^2 S(\vartheta)s}P
\\
&\quad+e^{-\tau ^2 S(\vartheta)s}P\tau(\vartheta _1 Z^*+\vartheta _2\widetilde{Z}^*)
-\tau ^3 \mathcal{M}(\tau;\vartheta;s) 
+\mathcal{R}(\tau;\vartheta ;s),
\end{split}
\end{equation}
where the remainder term $\mathcal{R}(\tau;\vartheta ;s)$ is subject to the estimate
\begin{align}
\Vert \mathcal{R}(\tau;\vartheta ;s)\Vert _{\mathfrak{H}\rightarrow\mathfrak{H}}
&\leqslant C_5 s^{-1}e^{-\check{c}_*\tau ^2 s/2},\quad s>0;
\nonumber
\\
\label{D0}
 C_5&=2\delta ^{-1}+\max _{\alpha\geqslant 0}\check{\phi}(\alpha)e^{-\check{c}_*\alpha/2}.
\end{align}
Here we introduced the notation 
\begin{equation*}
\begin{split}
\check{\phi}(\alpha)&=
2(C_1C_{F_1}+C_2)\alpha
+(C_3C_{F_1}+2C_1C_K+C_4)\alpha ^2
+2^{-1}C_3C_K\alpha ^3.
\end{split}
\end{equation*}
Changing the variables $\check{c}_*\alpha/2=:\alpha$, we can show that the constant $C_5$ is controlled by the polynomial of $\check{c}_*^{-1}$ and the data \eqref{data for abstract scheme}. 

Obviously, for $s\geqslant 1$ one has $s^{-1}\leqslant 2(s+1)^{-1}$, and, consequently, $\Vert \mathcal{R}(\tau;\vartheta ;s)\Vert \leqslant 2C_5(s+1)^{-1}e^{-\check{c}_*\tau ^2 s/2}$. For $0\leqslant s<1$, from \eqref{N(theta)<=}, \eqref{M=} it follows that for $\vert\tau\vert\leqslant\tau _0$ we have the inequality $\vert \tau \vert ^3 \Vert \mathcal{M}(\tau;\vartheta ;s)\Vert \leqslant \tau ^3_0C_N e^{-\check{c}_*\tau ^2 s}$. Thus, using \eqref{Z<=}, \eqref{tilde Z<=}, \eqref{B(tau,theta)>=}, and \eqref{S(theta)>=}, from \eqref{Th exp abstract theta tau} we derive
\begin{equation*}
\begin{split}
\Vert \mathcal{R}(\tau;\vartheta ;s)\Vert 
\leqslant C_{6}(1+s)^{-1}e^{-\check{c}_*\tau ^2 s/2},\quad 0\leqslant s<1,
\end{split}
\end{equation*}
where 
$
C_{6}=4 +4\kappa ^{1/2}(13\delta )^{-1/2}\tau _0(\Vert X_1\Vert +c_1C(1)^{1/2})+2\tau ^3_0C_N$.

Taking Remark~\ref{Remark constants Sec.1} into account, we arrive at the following result.

\begin{theorem}
\label{Theorem exp abstract theta tau}
For $\vert \tau\vert\leqslant\tau _0$, $s\geqslant 0$, identity \eqref{Th exp abstract theta tau} makes sense, 
where the operator $ \mathcal{M}(\tau;\vartheta;s)$ is defined in \eqref{M=}, and the operator $\mathcal{R}(\tau;\vartheta ;s)$ is subject to the estimates
\begin{align*}
\Vert \mathcal{R}(\tau;\vartheta ;s)\Vert &\leqslant C_7(s+1)^{-1}e^{-\check{c}_*\tau ^2 s/2},\quad s\geqslant 0;
\\
\Vert \mathcal{R}(\tau;\vartheta ;s)\Vert &\leqslant C_5s^{-1}e^{-\check{c}_*\tau ^2 s/2},\quad s>0.
\end{align*}
The constant $C_5$ was defined in \eqref{D0}. The constants $C_5$ and $C_7=\max\lbrace 2C_5;C_{6}\rbrace $ are controlled by polynomials of $\check{c}_*^{-1}$ and the data \eqref{data for abstract scheme}.
\end{theorem}

\subsection{Enumeration of the matrix elements of the operator $\mathcal{M}(\tau;\vartheta ;s)$} 

Let us look at the operator  \eqref{M=} in more detail. By \eqref{N=N0+N*},
\begin{equation}
\label{M=M0+M*}
\mathcal{M}(\tau;\vartheta ;s)= \mathcal{M}_0(\tau;\vartheta ;s)+\mathcal{M}_*(\tau;\vartheta ;s),
\end{equation}
where
\begin{align*}
\mathcal{M}_0(\tau;\vartheta ;s):&=\int\limits_0^s e^{-\tau ^2 S(\vartheta )P(s-\widetilde{s})}PN_0(\vartheta)e^{-\tau ^2 S(\vartheta )P\widetilde{s}}\,d\widetilde{s},
\\
\mathcal{M}_*(\tau;\vartheta ;s):&=\int\limits_0^s e^{-\tau ^2 S(\vartheta )P(s-\widetilde{s})}PN_*(\vartheta)e^{-\tau ^2 S(\vartheta )P\widetilde{s}}\,d\widetilde{s}.
\end{align*}
The operators $S(\vartheta)$ and $N_0(\vartheta)$ are diagonal in the basis $\lbrace \omega _l(\vartheta)\rbrace _{l=1}^n$, thus $N_0(\vartheta)S(\vartheta)=S(\vartheta )N_0(\vartheta)$ and
\begin{equation}
\label{M0=}
\mathcal{M}_0(\tau;\vartheta ;s)=N_0(\vartheta )e^{-\tau ^2 S(\vartheta)Ps}Ps.
\end{equation}
We calculate the matrix elements of the operator $\mathcal{M}_*(\tau;\vartheta ;s)$ in the basis  $\lbrace \omega _l(\vartheta)\rbrace _{l=1}^n$. According to \eqref{S omega _l=} we have $e^{-\tau ^2 S(\vartheta)Ps}\omega _j(\vartheta)=e^{-\tau ^2\gamma _j(\vartheta)s}\omega _j(\vartheta)$. Thus, the matrix elements of the operator $\mathcal{M}_*(\tau;\vartheta ;s)$ have the form
\begin{equation*}
\begin{split}
\Bigl(& \mathcal{M}_*(\tau;\vartheta ;s) \omega _j(\vartheta),\omega _k(\vartheta)\Bigr)
\\
&=\int\limits_0^s 
\left(
e^{-\tau ^2 S(\vartheta)P(s-\widetilde{s})}N_*(\vartheta)e^{-\tau ^2 S(\vartheta)P\widetilde{s}}\omega _j(\vartheta),\omega _k(\vartheta)\right)
\,d\widetilde{s}
\\
&=\bigl( N_*(\vartheta )\omega _j(\vartheta),\omega _k(\vartheta)\bigr)
e^{-\tau ^2 \gamma _k(\vartheta)s}\int\limits_0^s e^{-\tau ^2 (\gamma _j(\vartheta)-\gamma _k(\vartheta))\widetilde{s}}\,d\widetilde{s}.
\end{split}
\end{equation*}
Together with \eqref{N* matrix elements}, this implies that the functions $\gamma _j(\vartheta)=\gamma _k(\vartheta)$ satisfy 
$$
\big( \mathcal{M}_*(\tau;\vartheta ;s) \omega _j(\vartheta),\omega _k(\vartheta)\big)=0.
$$
Now, let $\gamma _j(\vartheta)\neq \gamma _k(\vartheta)$. By \eqref{omega i tilde omega tozd} and \eqref{N* matrix elements},
\begin{equation*}
\begin{split}
\Bigl( \mathcal{M}_*&(\tau;\vartheta ;s) \omega _j(\vartheta),\omega _k(\vartheta)\Bigr)
\\
&=\left( N_*(\vartheta)\omega _j(\vartheta),\omega _k(\vartheta)\right)
e^{-\tau ^2 \gamma _k(\vartheta)s}\frac{e^{-\tau ^2 (\gamma _j(\vartheta)-\gamma _k(\vartheta))s}-1}{(\gamma _k(\vartheta )-\gamma_j(\vartheta))\tau ^2}
\\
&=
\frac{e^{-\tau ^2\gamma _j(\vartheta)s}-e^{-\tau ^2 \gamma _k(\vartheta)s}}{\tau ^2}\left(\omega _j(\vartheta),\widetilde{\omega}_k(\vartheta)\right)
\\
&=\frac{e^{-\tau ^2 \gamma _j(\vartheta)s}}{\tau ^2}\left(\omega _j(\vartheta),\widetilde{\omega}_k(\vartheta)\right)
+\frac{e^{-\tau ^2 \gamma _k(\vartheta)s}}{\tau ^2}\left(\widetilde{\omega}_j(\vartheta),\omega _k(\vartheta)\right).
\end{split}
\end{equation*}
Thus,
\begin{equation}
\label{M* matrix element almost final}
\begin{split}
\Bigl(\mathcal{M}_*&(\tau;\vartheta ;s) \omega _j(\vartheta),\omega _k(\vartheta)\Bigr)
\\
&=\tau ^{-2}\Bigl(
\bigl(e^{-\tau ^2 \gamma _j(\vartheta)s}\omega _j(\vartheta),\widetilde{\omega}_k(\vartheta)\bigr)
+\bigr(\widetilde{\omega}_j(\vartheta),e^{-\tau ^2 \gamma _k(\vartheta)s}\omega _k(\vartheta)\bigr)
\Bigr).
\end{split}
\end{equation}
By \eqref{omega i tilde omega tozd}, for $\gamma _j(\vartheta)=\gamma _k(\vartheta)$ the right-hand side in \eqref{M* matrix element almost final} is equal to zero, so equality \eqref{M* matrix element almost final} holds true in this case. From \eqref{M* matrix element almost final} with the help of identities \eqref{omega i tilde omega tozd} we get the representation
\begin{equation}
\label{M*=}
\mathcal{M} _*(\tau;\vartheta ;s)=-\tau ^{-2}\sum _{l=1}^n e^{-\tau ^2 \gamma _l(\vartheta)s}
\Bigl(
\bigl(\,\cdot\, ,\widetilde{\omega}_l(\vartheta)\bigr)\omega _l(\vartheta)+\bigl( \,\cdot\, ,\omega _l(\vartheta)\bigr)\widetilde{\omega}_l(\vartheta)\Bigr).
\end{equation}
Combining \eqref{M=M0+M*}, \eqref{M0=}, and \eqref{M*=}, we find
\begin{equation*}
\begin{split}
\mathcal{M}(\tau;\vartheta ;s)&=N_0(\vartheta )e^{-\tau ^2 S(\vartheta)Ps}Ps
\\
&\quad-\tau ^{-2}\sum _{l=1}^n e^{-\tau ^2 \gamma _l(\vartheta)s}
\Bigl(
\bigl(\,\cdot\, ,\widetilde{\omega}_l(\vartheta)\bigr)\omega _l(\vartheta)+\bigl( \,\cdot\, ,\omega _l(\vartheta)\bigr)\widetilde{\omega}_l(\vartheta)\Bigr).
\end{split}
\end{equation*}

Finally, note that formula \eqref{Th exp abstract theta tau} simplifies if $N(\vartheta)=0$ or $N_*(\vartheta)=0$.

\begin{proposition}
Let $N(\vartheta)=0$. Then, under conditions of Theorem~\textnormal{\ref{Theorem exp abstract theta tau}} we have
\begin{equation*}
\begin{split}
e^{-B(\tau;\vartheta)s}=e^{-\tau ^2 S(\vartheta)s}P
&+\tau (\vartheta _1 Z+\vartheta _2\widetilde{Z})e^{-\tau ^2 S(\vartheta)s}P
\\
&+e^{-\tau ^2 S(\vartheta)s}P\tau(\vartheta _1 Z^*+\vartheta _2\widetilde{Z}^*)
+\mathcal{R}(\tau;\vartheta ;s).
\end{split}
\end{equation*}
\end{proposition}

\begin{proposition}
\label{Proposition N*(theta)=0}
Under conditions of Theorem~\textnormal{\ref{Theorem exp abstract theta tau}}, let $N_*(\vartheta)=0$. Then
\begin{equation*}
\begin{split}
e^{-B(\tau;\vartheta)s}\!=\!e^{-\tau ^2 S(\vartheta)s}P
\!+\!\tau (\vartheta _1 Z\!+\!\vartheta _2\widetilde{Z})e^{-\tau ^2 S(\vartheta)s}P
&+\!e^{-\tau ^2 S(\vartheta)s}P\tau(\vartheta _1 Z^*\!+\!\vartheta _2\widetilde{Z}^*)
\\
-\tau ^3 N_0(\vartheta )e^{-\tau ^2 S(\vartheta)Ps}Ps
&+\mathcal{R}(\tau;\vartheta ;s).
\end{split}
\end{equation*}
\end{proposition}

\subsection{Returning to parameters $t$ and $\varepsilon$}

Let us return back to the original parameters $t,\varepsilon$, recalling that $t=\tau\vartheta _1$, $\varepsilon =\tau\vartheta _2$. By \eqref{S(theta)=}, the operator $\tau ^2 S(\vartheta)=:L(t,\varepsilon)$ is given by the expression
\begin{equation}
\label{L(t,eps)}
\begin{split}
L(t,\varepsilon)
\!=\!t^2 S &\!+\!t\varepsilon \Big(\! -(X_0Z)^*X_0\widetilde{Z}
\!-\!(X_0\widetilde{Z})^*X_0Z
\!+\!P(Y_2^*Y_1\!+\!Y_1^*Y_2)\Big)\Big\vert _\mathfrak{N}
\\
&\!+\!\varepsilon ^2 \Big(-(X_0\widetilde{Z})^*X_0\widetilde{Z}\vert _\mathfrak{N}+Q_\mathfrak{N}+\lambda Q_{0\mathfrak{N}}\Big).
\end{split}
\end{equation}
Note an estimate, which follows from \eqref{S(theta)>=}:
\begin{equation}
\label{L(t,eps)>=}
L(t,\varepsilon)\geqslant\check{c}_*(t^2+\varepsilon ^2)I_\mathfrak{N},\quad t\in\mathbb{R},\quad 0\leqslant\varepsilon\leqslant 1.
\end{equation}

By \eqref{N(theta)=}, the operator $\tau ^3 N(\vartheta)=:N(t,\varepsilon)$ has the form
\begin{equation}
\label{N(t,eps)=}
N(t,\varepsilon)=t^3 N_{11}+t^2\varepsilon N_{12}+t\varepsilon ^2 N_{21}+\varepsilon ^3 N_{22}.
\end{equation}
By \eqref{N(theta)<=}, we have
\begin{equation}
\label{N(t,eps)<=}
\Vert N(t,\varepsilon)\Vert \leqslant C_N(t^2+\varepsilon ^2 )^{3/2},\quad t\in\mathbb{R},\quad 0\leqslant\varepsilon\leqslant 1.
\end{equation}

Rewriting \eqref{Th exp abstract theta tau} in terms of original parameters, we give an equivalent formulation of Theorem~\ref{Theorem exp abstract theta tau} that is convenient for further applications to differential operators. By the \textit{corrector} we mean the operator
\begin{equation}
\label{K(t,eps,s)}
\begin{split}
K(t,\varepsilon ,s)\!:=\!\big(tZ\!+\!\varepsilon \widetilde{Z}\big)e^{\!-L(t,\varepsilon)s}P
&\!+\!e^{\!-L(t,\varepsilon)s}P(tZ^*\!\!+\!\varepsilon \widetilde{Z}^*)
\\
&\!-\!\!\int\limits_0^s\! e^{-L(t,\varepsilon)(s-\widetilde{s})}PN(t,\varepsilon)e^{-L(t,\varepsilon)\widetilde{s}}P\,d\widetilde{s}.
\end{split}
\end{equation}

By \eqref{L(t,eps)>=} and \eqref{N(t,eps)<=},
\begin{equation}
\label{K(t,eps)<= start}
\begin{split}
\Vert K(t,\varepsilon ,s)\Vert 
\leqslant
2\max\lbrace \Vert Z\Vert ;\Vert \widetilde{Z}\Vert \rbrace (t&+\varepsilon) e^{-\check{c}_*(t^2+\varepsilon ^2)s}
\\
\quad+C_Ns(t^2&+\varepsilon ^2)^{3/2}e^{-\check{c}_*(t^2+\varepsilon ^2)s}.
\end{split}
\end{equation}
Combining \eqref{Z<=}, \eqref{tilde Z<=}, and \eqref{K(t,eps)<= start} and using the~elementary inequalities
$e^{-\alpha }\leqslant \alpha ^{-1}e^{-\alpha /2}$ and $e^{-\alpha}\leqslant 3\alpha ^{-2}e^{-\alpha /2}$, $\alpha >0$, we get
\begin{align}
\label{K(t,eps)<= for s>0}
\Vert  K(t,\varepsilon ,s)\Vert 
&\leqslant C_8 s^{-1}(t^2+\varepsilon ^2)^{-1/2}e^{-\check{c}_*(t^2+\varepsilon ^2)s/2},\quad s>0;
\\
C_8&=2^{3/2}\check{c}_*^{-1}\kappa ^{1/2}(13\delta )^{-1/2}\max\lbrace\Vert X_1\Vert ;c_1 C(1)^{1/2}\rbrace + 3\check{c}_*^{-2}C_N.
\nonumber
\end{align}
For $s\geqslant 0$, the estimate is more bulky: now for the exponentials we use inequalities of the form  
$e^{-\alpha}\leqslant 2(1+\alpha )^{-1}e^{-\alpha /2}$ and $e^{-\alpha}\alpha \leqslant 4(1+\alpha )^{-1}e^{-\alpha /2}$, where $\alpha\geqslant 0$. As a result, we obtain
\begin{align}
\label{K(t,eps)<= for s>=0}
\Vert  K(t,\varepsilon ,s)\Vert 
&\!\leqslant\! C_9 (t^2\!+\!\varepsilon ^2)^{1/2}(1\!+\!\check{c}_*(t^2\!+\!\varepsilon ^2)s)^{-1}e^{-\check{c}_*(t^2+\varepsilon ^2)s/2},\quad s\!\geqslant\! 0;
\\
C_9&=2^{5/2}\kappa ^{1/2}(13\delta )^{-1/2}\max\lbrace\Vert X_1\Vert ;c_1 C(1)^{1/2}\rbrace + 4\check{c}_*^{-1}C_N .
\nonumber
\end{align}

\begin{theorem}
\label{Theorem exp abstract}
We have
\begin{equation}
\label{2.40a}
e^{-B(t,\varepsilon)s}=e^{-L(t,\varepsilon)s}P+K(t,\varepsilon,s)+\mathcal{R}(t,\varepsilon,s),
\end{equation}
where the operators $B(t,\varepsilon)$, $L(t,\varepsilon)$, and $K(t,\varepsilon ,s)$ are defined in \eqref{B(t,eps) ne mathfrak}, \eqref{L(t,eps)}, and \eqref{K(t,eps,s)}, respectively. For the operator $K(t,\varepsilon ,s)$ we have estimates \eqref{K(t,eps)<= for s>0} and~\eqref{K(t,eps)<= for s>=0}. The operator $\mathcal{R}(t,\varepsilon,s):= \mathcal{R}(\tau;\vartheta ;s)$ for $t^2+\varepsilon ^2\leqslant\tau _0^2$ is subject to the  inequalities\textnormal{:}
\begin{align}
\label{R2(t,eps,s) estimates}
\Vert \mathcal{R}(t,\varepsilon ,s)\Vert &\leqslant C_7(s+1)^{-1}e^{-\check{c}_*(t^2+\varepsilon ^2) s/2},\quad s\geqslant 0;
\\
\label{R2(t,eps,s) estimates s>0}
\Vert \mathcal{R}(t,\varepsilon , s)\Vert &\leqslant C_5s^{-1}e^{-\check{c}_*(t^2+\varepsilon ^2) s/2},\quad s>0.
\end{align}
The constants $C_5$, $C_7$, $C_8$, and $C_9$ are controlled in terms of polynomials whose coefficients are positive numbers and whose variables are  $\check{c}_*^{-1}$ and the  data set \eqref{data for abstract scheme}.
\end{theorem}

\section{Approximation for the bordered operator exponential}

\label{Section bordered exp abstract}

\subsection{The operator family $A(t)=M^*\widehat{A}(t)M$}
\label{Subsection A(t) and hat-A(t)}

Let $\widehat{\mathfrak{H}}$ be yet another Hilbert space, and let $\widehat{X}(t)=\widehat{X}_0 +t\widehat{X}_1 \colon \widehat{\mathfrak{H}}\rightarrow \mathfrak{H}_*$ be the family of the form \eqref{X(t) in S1} satisfying conditions of Subsection~\ref{Subsubsection operator pencils}. We emphasize that the space  $\mathfrak{H}_*$ is the same as before. \textit{All the objects corresponding to $\widehat{X}(t)$ will be denoted by the mark} ``$\widehat{\phantom{a}}\,$''. 
Let $M:\mathfrak{H}\rightarrow\widehat{\mathfrak{H}}$ be an isomorphism
\begin{equation}
\label{MDomX0=}
M\mathrm{Dom}\,X_0=\mathrm{Dom}\,\widehat{X}_0,
\end{equation}
$X(t)=\widehat{X}(t)M\colon \mathfrak{H}\rightarrow\mathfrak{H}_*$; $X_0=\widehat{X}_0M$, $X_1=\widehat{X}_1M$. Then $A(t)=M^*\widehat{A}(t)M$, where $\widehat{A}(t)=\widehat{X}(t)^*\widehat{X}(t)$. Note that $\widehat{\mathfrak{N}}=M\mathfrak{N}$, 
$\widehat{n}=n$, and $\widehat{\mathfrak{N}}_*=\mathfrak{N}_*$, $\widehat{n}_*=n_*$, $\widehat{P}_*=P_*$. Set
\begin{equation}
\label{G= abstract scheme}
G=(MM^*)^{-1}\colon \widehat{\mathfrak{H}}\rightarrow\widehat{\mathfrak{H}}.
\end{equation}
Let $G_{\widehat{\mathfrak{N}}}$ be the block of the operator $G$ in the subspace $\widehat{\mathfrak{N}}$: $
G_{\widehat{\mathfrak{N}}}=\widehat{P}G\vert _{\widehat{\mathfrak{N}}} \colon \widehat{\mathfrak{N}}\rightarrow\widehat{\mathfrak{N}}$. 
Obviously, $G_{\widehat{\mathfrak{N}}}$ is an isomorphism in $\widehat{\mathfrak{N}}$. 

Let $\widehat{S}=\widehat{R}^*\widehat{R}\colon \widehat{\mathfrak{N}}\rightarrow\widehat{\mathfrak{N}}$ be the spectral germ of the operator family $\widehat{A}(t)$ with $t=0$. According to \cite[Chapter~1, Subsection~1.5]{BSu}, we have
\begin{equation}
\label{R and hat-R}
R=\widehat{R}M\vert _\mathfrak{N},\quad\mathrm{rank}\,R=\mathrm{rank}\,\widehat{R},
\end{equation}
and $S=PM^*\widehat{S}M\vert _\mathfrak{N}$.

\subsection{The operator family $B(t,\varepsilon)=M^*\widehat{B}(t,\varepsilon)M$}

\label{Subsection B(t,epa) and hat-B(t,eps)}

Assume that the operator $\widehat{Y}_0 \colon \widehat{\mathfrak{H}}\rightarrow\widetilde{\mathfrak{H}}$ satisfies the conditions of Subsection~\ref{Subsection Y(t) i Y2}. Note that the space $\widetilde{\mathfrak{H}}$ does not change. Set $Y_0:=\widehat{Y}_0M$, $M\mathrm{Dom}\,Y_0=\mathrm{Dom}\,\widehat{Y}_0$. By \eqref{MDomX0=} and condition $\mathrm{Dom}\,\widehat{X}_0\subset \mathrm{Dom}\,\widehat{Y}_0$, the inclusion $\mathrm{Dom}\,X_0\subset \mathrm{Dom}\,Y_0$ makes sense. Let $\widehat{Y}_1\colon  \widehat{\mathfrak{H}}\rightarrow \widetilde{\mathfrak{H}}$ be a bounded operator and let $Y_1=\widehat{Y}_1 M :\mathfrak{H}\rightarrow\widetilde{\mathfrak{H}}$. Set $\widehat{Y}(t):=\widehat{Y}_0+t\widehat{Y}_1 \colon \widehat{\mathfrak{H}}\rightarrow \widetilde{\mathfrak{H}}$, $\mathrm{Dom}\,\widehat{Y}(t)=\mathrm{Dom}\,\widehat{Y}_0$. Let $Y(t)=\widehat{Y}(t)M=Y_0+tY_1 : \mathfrak{H}\rightarrow \widetilde{\mathfrak{H}}$, $\mathrm{Dom}\,Y(t)=\mathrm{Dom}\,Y_0$. Assume that the operators $\widehat{X}(t)$ and $\widehat{Y}(t)$ are subject to Condition~\ref{Condition Y(t)} with some constant $\widehat{c}_1$. Then we automatically have $\Vert Y(t)u\Vert _{\widetilde{\mathfrak{H}}}\leqslant c_1 \Vert X(t)u\Vert _{\mathfrak{H}_*}$, $u\in\mathrm{Dom}\,X_0$, where $c_1=\widehat{c}_1$.

Let the operator $\widehat{Y}_2\colon \widehat{\mathfrak{H}}\rightarrow\widetilde{\mathfrak{H}}$ satisfy conditions of Subsection~\ref{Subsection Y(t) i Y2}. Put $Y_2:=\widehat{Y}_2M \colon \mathfrak{H}\rightarrow\widetilde{\mathfrak{H}}$, $M\mathrm{Dom}\,Y_2=\mathrm{Dom}\,\widehat{Y}_2$. Since $M$ is an isomorphism, and the operator $\widehat{Y}_2$ is densely defined, the operator $Y_2$ is also densely defined. By \eqref{MDomX0=}, the inclusion $\mathrm{Dom}\,X_0\subset\mathrm{Dom}\,Y_2$ holds. Assume that the operators $\widehat{X}(t)$ and $\widehat{Y}_2$ are subject to Condition~\ref{Condition Y2} with some constants $\widehat{C}(\nu)>0$. Then, automatically, for any  $\nu >0$   where exists a constant $C(\nu)=\widehat{C}(\nu)\Vert M\Vert ^2 >0$ such that for $u\in\mathrm{Dom}\,X_0$, $t\in\mathbb{R}$, we have $\Vert Y_2u\Vert ^2_{\widetilde{\mathfrak{H}}}\leqslant\nu \Vert X(t)u\Vert ^2_{\mathfrak{H}_*}+C(\nu)\Vert u\Vert ^2 _{\mathfrak{H}}$.

Put $Q_0:=M^*M$. Then $Q_0$ is a~bounded positive definite operator in $\mathfrak{H}$. (The role of $\widehat{Q}_0$ is played by the identity operator in $\widehat{\mathfrak{H}}$.)

In the space $\widehat{\mathfrak{H}}$, consider the form $\widehat{\mathfrak{q}}$ satisfying the conditions of Subsection~\ref{Subsection form g abstract scheme}. Define a~form $\mathfrak{q}$ acting by the rule $\mathfrak{q}[u,v]=\widehat{\mathfrak{q}}[Mu,Mv]$, $u,v\in\mathrm{Dom}\,\mathfrak{q}$, $M\mathrm{Dom}\,\mathfrak{q}=\mathrm{Dom}\,\widehat{\mathfrak{q}}$. Formally, $Q=M^*\widehat{Q}M$. Assume that the operator $\widehat{X}(t)$ and the form $\widehat{\mathfrak{q}}$ satisfy Condition~\ref{Condition q} with some constants $\kappa$, $\widehat{c}_0$, $\widehat{c}_2$, and $\widehat{c}_3$. Taking \eqref{MDomX0=} into account, one can check that the operator $X(t)=\widehat{X}(t)M$ and the form $\mathfrak{q}$ also satisfy Condition~\ref{Condition q} with the constants
\begin{equation}
\label{c0 relations}
c_0=\Vert M\Vert ^2\widehat{c}_0,\quad\mbox{if }\widehat{c}_0\geqslant 0,\quad
c_0=\Vert M^{-1}\Vert ^{-2}\widehat{c}_0,\quad\mbox{if }\widehat{c}_0<0,
\end{equation}
$c_2=\widehat{c}_2$, $c_3=\Vert M\Vert ^2\widehat{c}_3$, and the same constant $\kappa$. From \eqref{c4} it follows that the constants $c_4$ and $\widehat{c}_4=4\kappa ^{-1}\widehat{c}_1^{\,2}\widehat{C}(\nu)$ for $\nu=\kappa ^2(16\widehat{c}_1^{\,2})^{-1}$ are related by the identity
\begin{equation}
\label{c4 and hat c4}
c_4=\Vert M\Vert ^2 \widehat{c}_4.
\end{equation}

Under the above assumptions, the operator pencil
\begin{equation}
\label{hat B(t,eps)}
\widehat{B}(t,\varepsilon)=\widehat{A}(t)+\varepsilon(\widehat{Y}_2^*\widehat{Y}(t)+\widehat{Y}(t)^*\widehat{Y}_2)+
\varepsilon ^2\widehat{Q}+\lambda\varepsilon ^2 I
\end{equation}
is related to the pencil~\eqref{B(t,eps) ne mathfrak} by the identity 
\begin{equation}
\label{B=M*hat-B M abstract}
B(t,\varepsilon)=M^*\widehat{B}(t,\varepsilon)M. 
\end{equation}
The constant $\lambda$ is chosen from condition~\eqref{lambda condition} for the operator \eqref{B(t,eps) ne mathfrak}. By using relations \eqref{c0 relations}, \eqref{c4 and hat c4}, and the equality $Q_0=M^*M$, we find that, under our choice of $\lambda$, condition~\eqref{lambda condition} for the operator~\eqref{hat B(t,eps)} is also satisfied.

Note that for the operator \eqref{hat B(t,eps)} both variants of relation \eqref{beta condition} has the form $\widehat{\beta}=\lambda-\widehat{c}_0-\widehat{c}_4$. Together with \eqref{beta condition}, \eqref{c0 relations}, and \eqref{c4 and hat c4}, this implies that 
\begin{equation}
\label{beta<=||M^-1||^-2hat beta}
\beta \leqslant \Vert M^{-1}\Vert ^{-2}\widehat{\beta}.
\end{equation}

\subsection{The relationship between the spectral germs $S(\vartheta)$ and $\widehat{S}(\vartheta)$}

In expressions for the operators $B(t,\varepsilon)$ and $\widehat{B}(t,\varepsilon)$, we now switch to the parameters $\tau$, $\vartheta$. 
Consider the spectral germ \eqref{S(theta)=} and a similar germ for the family \eqref{hat B(t,eps)}:
\begin{equation}
\label{hat S(vartheta)=}
\begin{split}
\widehat{S}(\vartheta)&\!=\!
\vartheta _1 ^2\widehat{S}
\!-\!\vartheta _1\vartheta _2(\widehat{X}_0\widehat{Z})^*(\widehat{X}_0\widehat{\widetilde{Z}})\vert _{\widehat{\mathfrak{N}}}\!-\!\vartheta _1\vartheta _2 (\widehat{X}_0\widehat{\widetilde{Z}})^*(\widehat{X}_0\widehat{Z})\vert _{\widehat{\mathfrak{N}}}
\\
&\quad\!-\!\vartheta _2^2 (\widehat{X}_0\widehat{\widetilde{Z}})^*(\widehat{X}_0\widehat{\widetilde{Z}})\vert _{\widehat{\mathfrak{N}}}
\!+\!\vartheta _1\vartheta _2 \widehat{P}(\widehat{Y}_2^*\widehat{Y}_1+\widehat{Y}_1^*\widehat{Y}_2)\vert _{\widehat{\mathfrak{N}}}
\!+\!\vartheta _2 ^2(\widehat{Q}_{\widehat{\mathfrak{N}}}\!+\!\lambda I_{\widehat{\mathfrak{N}}}).
\end{split}
\end{equation}
In \cite[Proposition~1.8]{M_AA}, the following relation between operators \eqref{S(theta)=} and \eqref{hat S(vartheta)=} was proved:
\begin{equation}
\label{S(theta) and hat-S(theta)}
S(\vartheta)=PM^*\widehat{S}(\vartheta)M\vert _\mathfrak{N}.
\end{equation}

Let's go back to the parameters $t$ and $\varepsilon$. 
The operator $\tau ^2 \widehat{S}(\vartheta)=:\widehat{L}(t,\varepsilon)$ is given by the expression
\begin{equation}
\label{hat L(t,eps)}
\begin{split}
\widehat{L}(t,\varepsilon)
=t^2 \widehat{S} &+t\varepsilon \Big(\!\! -(\widehat{X}_0\widehat{Z})^*\widehat{X}_0\widehat{\widetilde{Z}}
\!-\!(\widehat{X}_0\widehat{\widetilde{Z}})^*\widehat{X}_0\widehat{Z}
\!+\!\widehat{P}(\widehat{Y}_2^*\widehat{Y}_1+\widehat{Y}_1^*\widehat{Y}_2)\Big)\Big\vert _{\widehat{\mathfrak{N}}}
\\
&+\varepsilon ^2 \Big(\!\!-(\widehat{X}_0\widehat{\widetilde{Z}})^*\widehat{X}_0\widehat{\widetilde{Z}}\vert _{\widehat{\mathfrak{N}}}+\widehat{Q}_{\widehat{\mathfrak{N}}}+\lambda I_{\widehat{\mathfrak{N}}}\Big).
\end{split}
\end{equation}
By \eqref{S(theta) and hat-S(theta)}, $L(t,\varepsilon)=PM^*\widehat{L}(t,\varepsilon)M\vert _{\mathfrak{N}}$. 
Here $L(t,\varepsilon)$ is the operator~\eqref{L(t,eps)}.

\subsection{The operators $\widehat{Z}_G$ and $\widehat{\widetilde{Z}}_G$}
\label{Subsection ZG and tilde ZG abstract scheme}

Let $\widehat{Z}_G$ be the operator in $\widehat{\mathfrak{H}}$ mapping an element $\widehat{u}\in\widehat{\mathfrak{H}}$ into the (unique) solution $\widehat{\phi}_G$ of the equation
\begin{equation*}
\widehat{X}_0^*(\widehat{X}_0\widehat{\phi}_G+\widehat{X}_1\widehat{\omega})=0,\quad G\widehat{\phi}_G\perp \widehat{\mathfrak{N}},
\end{equation*}
where $\widehat{\omega}=\widehat{P}\widehat{u}$. The equation is understood in a weak sense (cf. \eqref{def Z phi tozd}). Thus, as was shown in \cite[Lemma~6.1]{BSu05-1},
\begin{equation}
\label{hat_Z-G and Z}
\widehat{Z}_G=MZM^{-1}\widehat{P}.
\end{equation}

Similarly, let $\widehat{\widetilde{Z}}_G$ be the operator in $\widehat{\mathfrak{H}}$ mapping an element $\widehat{u}\in\widehat{\mathfrak{H}}$ into the (unique) solution $\widehat{\psi}_G$ of the equation
\begin{equation*}
\widehat{X}_0^*\widehat{X}_0\widehat{\psi}_G+\widehat{Y}_0^*\widehat{Y}_2\widehat{\omega}=0,\quad G\widehat{\psi}_G\perp \widehat{\mathfrak{N}},
\end{equation*}
where $\widehat{\omega}=\widehat{P}\widehat{u}$. The equation is understood in a weak sense. By recalculation in equation \eqref{eq tilde Z abstract} and the identities $M\mathfrak{N}=\widehat{\mathfrak{N}}$, \eqref{MDomX0=}, and \eqref{G= abstract scheme}, we get
\begin{equation}
\label{hat_tilde-Z-G and tilde Z}
\widehat{\widetilde{Z}}_G=M\widetilde{Z}M^{-1}\widehat{P}.
\end{equation}

\subsection{The operator $\widehat{N}_G(t,\varepsilon)$}

Using relationship between operators $X_0$, $X_1$, $Y_0$, $Y_1$, $Y_2$, and $Q$ and the corresponding operators marked by the ``hat'', applying the identities $\widehat{\mathfrak{N}}=M\mathfrak{N}$, $Q_0=M^*M$, \eqref{N11 abstract}--\eqref{N22 abstract}, \eqref{N(t,eps)=},  \eqref{R and hat-R}, \eqref{hat_Z-G and Z}, and \eqref{hat_tilde-Z-G and tilde Z}, we conclude that the operator
\begin{equation}
\label{NG def in abstr scheme}
\widehat{N}_G(t,\varepsilon):=\widehat{P}(M^*)^{-1}N(t,\varepsilon)M^{-1}\widehat{P}
\end{equation}
has the form
\begin{equation}
\label{NG(t,eps):=}
\widehat{N}_G(t,\varepsilon)=t^3 \widehat{N}_{G,11}+t^2\varepsilon \widehat{N}_{G,12}
+t\varepsilon ^2 \widehat{N}_{G,21}+\varepsilon ^3 \widehat{N}_{G,22},
\end{equation}
where
\begin{align}
\label{NG11 abstract}
\widehat{N}_{G,11}&=(\widehat{X}_1\widehat{Z}_G)^*\widehat{R}\widehat{P}+(\widehat{R}\widehat{P})^*\widehat{X}_1\widehat{Z}_G,
\\
\begin{split}
\widehat{N}_{G,12}&=(\widehat{X}_1\widehat{\widetilde{Z}}_G)^*\widehat{R}\widehat{P}+(\widehat{R}\widehat{P})^*\widehat{X}_1\widehat{\widetilde{Z}}_G+(\widehat{X}_1\widehat{Z}_G)^*\widehat{X}_0\widehat{\widetilde{Z}}_G
\\
&+(\widehat{X}_0\widehat{\widetilde{Z}}_G)^*\widehat{X}_1\widehat{Z}_G
+(\widehat{Y}_2\widehat{Z}_G)^*\widehat{Y}_0\widehat{Z}_G+(\widehat{Y}_0\widehat{Z}_G)^*\widehat{Y}_2\widehat{Z}_G
\\
&+(\widehat{Y}_2\widehat{Z}_G)^*\widehat{Y}_1\widehat{P}+(\widehat{Y}_1\widehat{P})^*\widehat{Y}_2\widehat{Z}_G
+(\widehat{Y}_2\widehat{P})^*\widehat{Y}_1\widehat{Z}_G+(\widehat{Y}_1\widehat{Z}_G)^*\widehat{Y}_2\widehat{P},
\end{split}
\\
\begin{split}
\widehat{N}_{G,21}&=(\widehat{X}_0\widehat{\widetilde{Z}}_G)^*\widehat{X}_1\widehat{\widetilde{Z}}_G
+(\widehat{X}_1\widehat{\widetilde{Z}}_G)^*\widehat{X}_0\widehat{\widetilde{Z}}_G
+(\widehat{Y}_2\widehat{Z}_G)^*\widehat{Y}_0\widehat{\widetilde{Z}}_G
\\
&+(\widehat{Y}_0\widehat{\widetilde{Z}}_G)^*\widehat{Y}_2\widehat{Z}_G
+(\widehat{Y}_2\widehat{\widetilde{Z}}_G)^*\widehat{Y}_0\widehat{Z}_G
+(\widehat{Y}_0\widehat{Z}_G)^*\widehat{Y}_2\widehat{\widetilde{Z}}_G
+(\widehat{Y}_2\widehat{\widetilde{Z}}_G)^*\widehat{Y}_1\widehat{P}
\\
&+(\widehat{Y}_1\widehat{P})^*\widehat{Y}_2\widehat{\widetilde{Z}}_G
+(\widehat{Y}_1\widehat{\widetilde{Z}}_G)^*\widehat{Y}_2\widehat{P}
+(\widehat{Y}_2\widehat{P})^*\widehat{Y}_1\widehat{\widetilde{Z}}_G
\\
&+\widehat{Z}_G^*\widehat{Q}\widehat{P}+\widehat{P}\widehat{Q}\widehat{Z}_G
+\lambda (\widehat{Z}_G^*\widehat{P}+\widehat{P}\widehat{Z}_G),
\end{split}
\\
\label{N_G22 abstract}
\begin{split}
\widehat{N}_{G,22}&=(\widehat{Y}_0\widehat{\widetilde{Z}}_G)^*\widehat{Y}_2\widehat{\widetilde{Z}}_G
+(\widehat{Y}_2\widehat{\widetilde{Z}}_G)^*\widehat{Y}_0\widehat{\widetilde{Z}}_G
+\bigl(\widehat{\widetilde{Z}}_G\bigr)^*\widehat{Q}\widehat{P}
+\widehat{P}\widehat{Q}\widehat{\widetilde{Z}}_G
\\
&+\lambda\Bigl(\bigl(\widehat{\widetilde{Z}}_G\bigr)^*\widehat{P}
+\widehat{P}\widehat{\widetilde{Z}}_G\Bigr)
.
\end{split}
\end{align}
Note the estimate that follows from \eqref{N(t,eps)<=} and \eqref{NG def in abstr scheme}:
\begin{equation}
\label{NG<= abstract scheme}
\Vert \widehat{N}_G(t,\varepsilon)\Vert \leqslant C_G (t^2+\varepsilon ^2)^{3/2};\quad C_G= C_N\Vert M^{-1}\Vert ^2 .
\end{equation}
By Remark~\ref{Remark constants Sec.1}, the following observation holds true.

\begin{remark}
\label{Remark CG}
The constant $C_G$ is controlled in terms of a polynomial with positive numerical coefficients and variables  $\Vert M^{-1}\Vert $ and \eqref{data for abstract scheme}.
\end{remark}

\subsection{Approximation for the bordered operator exponential}

In the present subsection, our goal is to approximate the operator exponential $e^{-B(t,\varepsilon)s}$ generated by the operator 
\eqref{B=M*hat-B M abstract} in terms of the isomorphism $M$ and the threshold characteristics of the operator $\widehat{B}(t,\varepsilon)$. 

Assume that the operator ${A}(t)$ is subject to inequality \eqref{A(t)>= abstract} with some constant ${c}_*>0$: 
$
{A}(t)\geqslant {c}_*t^2I$, $\vert t\vert \leqslant {\tau}_0$. 
Here $\tau _0$ satisfies condition \eqref{tau 0 in abstract scheme} for the operator $B(t,\varepsilon)$. Then estimate \eqref{B(tau,theta)>=} is fulfilled for $B(t,\varepsilon)$. 

Therefore, the operator $B(t,\varepsilon)$ of the form \eqref{B=M*hat-B M abstract} satisfies all the assumptions of Theorem~\ref{Theorem exp abstract}. Applying that theorem and multiplying \eqref{2.40a} by the operator $M$ from the left and by $M^*$ from the right, we obtain
\begin{equation}
\label{M Th exp M*}
\begin{split}
Me^{-B(t,\varepsilon)s}M^*&=Me^{-L(t,\varepsilon)s}PM^*
+MK(t,\varepsilon,s)M^*
+M\mathcal{R}(t,\varepsilon,s)M^*.
\end{split}
\end{equation}
By \eqref{K(t,eps,s)},
\begin{equation}
\label{K(t,eps,s) bordered in sec 3}
\begin{split}
MK(t,\varepsilon,s)M^*
&=M(tZ+\varepsilon \widetilde{Z})e^{-L(t,\varepsilon)s}PM^*
\\
&+Me^{-L(t,\varepsilon)s}P(tZ^*+\varepsilon \widetilde{Z}^*)M^*
\\
&-\int\limits_0^s M e^{-L(t,\varepsilon)(s-\widetilde{s})}PN(t,\varepsilon)Pe^{-L(t,\varepsilon)\widetilde{s}}PM^*\,d\widetilde{s}.
\end{split}
\end{equation}

Now we invoke the identity obtained in \cite[Proposition 3.1]{M_AA}:
\begin{equation}
\label{exp eff M tozd}
Me^{-L(t,\varepsilon)s}PM^*=M_0e^{-M_0\widehat{L}(t,\varepsilon)M_0 s}M_0\widehat{P}.
\end{equation}
Here $L(t,\varepsilon)$ and $\widehat{L}(t,\varepsilon)$ are the operators \eqref{L(t,eps)} and \eqref{hat L(t,eps)}, respectively, and
\begin{equation}
\label{M0}
M_0:=(G_{\widehat{\mathfrak{N}}})^{-1/2}.
\end{equation}

Denote $K_G(t,\varepsilon,s):= MK(t,\varepsilon,s)M^*$. We note at once two estimates following from \eqref{K(t,eps)<= for s>0} and \eqref{K(t,eps)<= for s>=0}:
\begin{align}
\label{KG (t,eps)<= for s>0}
\Vert & K _G(t,\varepsilon ,s)\Vert 
\!\leqslant\! C_8\Vert M\Vert ^2 s^{-1}(t^2\!+\!\varepsilon ^2)^{-1/2}e^{-\check{c}_*(t^2+\varepsilon ^2)s/2},\quad s\!>\!0;
\\
\label{KG (t,eps)<= for s>=0}
\begin{split}
\Vert & K_G(t,\varepsilon ,s)\Vert 
\!\leqslant\! C_9\Vert M\Vert ^2 (t^2+\varepsilon ^2)^{1/2}\big(1\!+\!\check{c}_*(t^2\!\!+\!\varepsilon ^2)s\big)^{-1}e^{-\check{c}_*(t^2\!+\varepsilon ^2)s/2},\\
&\qquad\qquad\qquad\qquad\qquad\qquad\qquad\qquad\qquad\qquad\qquad\qquad\qquad s\geqslant 0.
\end{split}
\end{align}
By \eqref{hat_Z-G and Z}--\eqref{NG def in abstr scheme}, \eqref{K(t,eps,s) bordered in sec 3}, and \eqref{exp eff M tozd}, the operator $K_G(t,\varepsilon,s)$ can be represented as
\begin{equation}
\label{KG(t,eps,s) tozd Sec 3}
\begin{split}
K_G(t,\varepsilon,s)
&\!=\!\bigl(t\widehat{Z}_G\!+\!\varepsilon\widehat{\widetilde{Z}}_G\bigr)M_0e^{-M_0\widehat{L}(t,\varepsilon)M_0 s}\!M_0\widehat{P}
\\
&\!+\!M_0e^{-M_0\widehat{L}(t,\varepsilon)M_0 s}\!M_0\widehat{P}\Big(t\widehat{Z}_G^*+\varepsilon\bigl(\widehat{\widetilde{Z}}_G\bigr)^*\Big)
\\
&\!-\!\!\int\limits_0^s\! M_0e^{-M_0\widehat{L}(t,\varepsilon)M_0 (s-\widetilde{s})}\!M_0\widehat{N}_G(t,\varepsilon)M_0e^{-M_0\widehat{L}(t,\varepsilon)M_0 \widetilde{s}}\!M_0\widehat{P}\,d\widetilde{s}.
\end{split}
\end{equation}

Combining \eqref{R2(t,eps,s) estimates}, \eqref{R2(t,eps,s) estimates s>0}, \eqref{M Th exp M*}, \eqref{exp eff M tozd}, and \eqref{KG(t,eps,s) tozd Sec 3},  we arrive at the following result.

\begin{theorem}
\label{Theorem bordered exp final in abstrach scheme}
Let the assumptions of Subsections~\textnormal{\ref{Subsection A(t) and hat-A(t)}} and \textnormal{\ref{Subsection B(t,epa) and hat-B(t,eps)}} be satisfied. Suppose that the operator ${A}(t)$ is subject to inequality \eqref{A(t)>= abstract}. Let the operator $\widehat{L}(t,\varepsilon)$ be defined in \eqref{hat L(t,eps)}. Let $M_0$ be the operator \eqref{M0}. Then 
\begin{equation*}
\begin{split}
Me^{-B(t,\varepsilon)s}M^*=M_0e^{-M_0\widehat{L}(t,\varepsilon)M_0 s}M_0\widehat{P}
+K_G(t,\varepsilon,s)
+M\mathcal{R}(t,\varepsilon,s)M^*.
\end{split}
\end{equation*}
Here $K_G(t,\varepsilon,s)$ is the operator \eqref{KG(t,eps,s) tozd Sec 3} that satisfies estimates \eqref{KG (t,eps)<= for s>0} and \eqref{KG (t,eps)<= for s>=0}, 
and the remainder is subject to the inequalities 
\begin{align*}
\Vert M\mathcal{R}(t,\varepsilon,s)M^*\Vert &\leqslant C_7\Vert M\Vert ^2
(s+1)^{-1}e^{-\check{c}_*(t^2+\varepsilon ^2) s/2},\quad s\geqslant 0;
\\
\Vert M\mathcal{R}(t,\varepsilon,s)M^*\Vert &\leqslant C_5\Vert M\Vert ^2
s^{-1}e^{-\check{c}_*(t^2+\varepsilon ^2) s/2},\quad s>0.
\end{align*}
The constants $C_5$, $C_7$, $C_8$, and $C_9$ are controlled in terms of polynomials with positive numerical coefficients and with the variables $\check{c}_*^{-1}$ and \eqref{data for abstract scheme}.
\end{theorem}

\section*{Chapter 2. Periodic differential operators in $L_2(\mathbb{R}^d;\mathbb{C}^n)$}

\section{Preliminaries}

\label{Section Pervaritel'nye svedeniya}

\subsection{The lattices $\Gamma $ and $\widetilde{\Gamma }$} 
Let $\Gamma$~be a lattice in $\mathbb{R}^d$ generated by a basis
$\mathbf{a}_1,\dots,\mathbf{a}_d$: $\Gamma=\lbrace \mathbf{a}\in 
\mathbb{R}^d:\; \mathbf{a}=\sum _{j=1}^d n^j\mathbf{a}_j,\; n^j\in 
\mathbb{Z}\rbrace$. By $\Omega$ we denote the elementary cell of the lattice 
$\Gamma$: $\Omega=\lbrace\mathbf{x}\in \mathbb{R}^d:\; \mathbf{x}=\sum 
_{j=1}^d\xi ^j\mathbf{a}_j,\; 0<\xi ^j<1\rbrace$. The basis 
$\mathbf{b}^1,\dots,\mathbf{b}^d$, dual to
$\mathbf{a}_1,\dots,\mathbf{a}_d$, is defined by the relations
$\langle\mathbf{b}^l, \mathbf{a}_j\rangle=2\pi \delta ^l_j$. This basis generates a lattice 
 $\widetilde{\Gamma}$ dual to the lattice~$\Gamma$. 
By $\widetilde{\Omega}$ we denote the first \textit{Brillouin zone} of the lattice 
$\widetilde{\Gamma}$: 
$$\widetilde{\Omega}=\lbrace \mathbf{k}\in 
\mathbb{R}^d:\; \vert \mathbf{k}\vert<\vert \mathbf{k}-\mathbf{b}\vert ,\; 
0\neq \mathbf{b}\in \widetilde{\Gamma}\rbrace.$$ 
The domain  
$\widetilde{\Omega}$ is fundamental for $\widetilde{\Gamma}$. 
We use the notation $\vert \Omega\vert 
=\textrm{meas}\,\Omega,$ $\vert \widetilde{\Omega}\vert 
=\textrm{meas}\,\widetilde{\Omega}$. Let $r_0$~be the radius of the ball inscribed in 
$\textrm{clos}\;\widetilde{\Omega}$ and let $2r_1 = 
\text{diam}\, \widetilde{\Omega}$.

If $\Phi (\mathbf{x})$~is a~$\Gamma$-periodic matrix-valued function in~$\mathbb{R}^d$, denote
\begin{equation}
\label{underline and overline definition}
\overline{\Phi}:=\vert \Omega\vert ^{-1}\int\limits_\Omega \Phi(\mathbf{x})\,d\mathbf{x},\qquad
 \underline{\Phi}:=\bigg(\vert \Omega\vert ^{-1}\int\limits_\Omega \Phi(\mathbf{x})^{-1}\,d\mathbf{x}\bigg)^{-1}. 
\end{equation}
Here in the definition of $\overline{\Phi}$ it is assumed that $\Phi\in L_{1,\mathrm{loc}}(\mathbb{R}^d)$ and in the definition of $\underline{\Phi}$ it is assumed that $\Phi$ is squared non-singular matrix-valued function and $\Phi^{-1}\in L_{1,\mathrm{loc}}(\mathbb{R}^d)$. 

\subsection{The factorized second-order operator families} 
\label{Subsection Factorized DOs}
(See \cite{BSu}.)
Let $b(\mathbf{D})= \sum_{l=1}^d b_l D_l  :  
L_2(\mathbb{R}^d;\mathbb{C}^n)\rightarrow 
L_2(\mathbb{R}^d;\mathbb{C}^m)$~be a first order DO. 
Here $b_l$~are $(m\times n)$-matrices with constant entries. 
\textit{Suppose that} $m\geqslant n$.
The symbol $b(\boldsymbol{\xi})= \sum_{l=1}^d b_l \xi_l$ 
\textit{is assumed to satisfy} $\textrm{rank}\,b(\boldsymbol{\xi})=n$, $0\neq 
\boldsymbol{\xi}\in \mathbb{R}^d$. Then for some $\alpha _0,\, 
\alpha _1 >0$ we have
\begin{equation}
\label{<b^*b<}
\alpha _0\mathbf{1}_n \leqslant b(\boldsymbol{\theta})^*b(\boldsymbol{\theta}) 
\leqslant \alpha _1\mathbf{1}_n,\quad 
\boldsymbol{\theta}\in \mathbb{S}^{d-1},\quad 
0<\alpha _0\leqslant \alpha _1<\infty.
\end{equation}
Let an $(n\times n)$-matrix-valued function $f(\mathbf{x})$ and an $(m\times 
m)$-matrix-valued function $h(\mathbf{x})$, $\mathbf{x} \in \mathbb{R}^d$,~be bounded together with the inverses:
\begin{equation}
\label{f,h}
f,\,f^{-1}\in L_\infty(\mathbb{R}^d);\quad h,\;h^{-1}\in L_\infty(\mathbb{R}^d).
\end{equation}
The functions $f$ and $h$ are $\Gamma$-periodic. Consider the DO
\begin{align}
\label{DO X}
&\mathcal{X}:=hb(\mathbf{D})f\,\colon \,L_2(\mathbb{R}^d;\mathbb{C}^n)\rightarrow L_2(\mathbb{R}^d;\mathbb{C}^m),
\\
\label{Dom DO X}
&\textrm{Dom}\,\mathcal{X}:=\lbrace \mathbf{u}\in L_2(\mathbb{R}^d;\mathbb{C}^n)\; \colon \; f\mathbf{u}\in H^1(\mathbb{R}^d;\mathbb{C}^n)\rbrace.
\end{align}
The operator \eqref{DO X} is closed on the domain \eqref{Dom DO X}. 
In~$L_2(\mathbb{R}^d;\mathbb{C}^n)$, consider a self-adjoint operator  
$\mathcal{A}:=\mathcal{X^*X}$ corresponding to the quadratic form  
\begin{equation}
\label{frak a[u,u]}
\mathfrak{a}[\mathbf{u},\mathbf{u}]=\Vert\mathcal{X}\mathbf{u}\Vert^2_{L_2(\mathbb{R}^d)},\quad 
\mathbf{u}\in \textrm{Dom}\,\mathcal{X}.
\end{equation} 
Formally,
$\mathcal{A}=f^*b(\mathbf{D})^*gb(\mathbf{D})f$, where  $g=h^*h$. By using the Fourier transform and 
\eqref{<b^*b<}, and \eqref{f,h}, it is easily seen that 
for $\mathbf{u}\in \textrm{Dom}\,\mathcal{X}$ we have
\begin{equation}
\label{<a[u,u]<}
\alpha _0\Vert g^{-1}\Vert _{L_\infty}^{-1}\Vert \mathbf{D}(f\mathbf{u})\Vert ^2_{L_2(\mathbb{R}^d)}\leqslant \mathfrak{a}[\mathbf{u},\mathbf{u}]\leqslant \alpha _1\Vert g \Vert _{L_\infty}\Vert \mathbf{D}(f\mathbf{u})\Vert ^2_{L_2(\mathbb{R}^d)}.
\end{equation}

\subsection{The operators $\mathcal{Y}$ and $\mathcal{Y}_2$} 
Now we introduce lower order terms. Let the operator 
$\mathcal{Y}:L_2(\mathbb{R}^d;\mathbb{C}^{n})\rightarrow L_2(\mathbb{R}^d;\mathbb{C}^{dn})$ act by the rule
\begin{equation}
\label{mathcal Y}
\mathcal{Y}\textbf{u}=\mathbf{D}(f\mathbf{u})
=\textrm{col}\,\lbrace D_1(f\mathbf{u}),\dots ,D_d(f\mathbf{u})\rbrace ,\quad 
\textrm{Dom}\,\mathcal{Y}=\textrm{Dom}\,\mathcal{X}.
\end{equation}
The lower estimate \eqref{<a[u,u]<} means that
\begin{equation}
\label{Yu<=}
\Vert \mathcal{Y}\mathbf{u}\Vert _{L_2(\mathbb{R}^d)}\leqslant c_1\Vert \mathcal{X}\mathbf{u}\Vert _{L_2(\mathbb{R}^d)},\quad \mathbf{u}\in \textrm{Dom}\,\mathcal{X};\quad
c_1=\alpha _0^{-1/2}\Vert g ^{-1}\Vert _{L_\infty }^{1/2}.
\end{equation}

In $\mathbb{R}^d$, let $\Gamma $-periodic $(n\times 
n)$-matrix-valued functions $a_j(\mathbf{x})$,  $j=1,\dots ,d$, be such that
\begin{equation}
\label{a_j}
a_j\in L_\varrho (\Omega),\quad \varrho =2\;\mbox{for}\;d=1,\quad 
\varrho >d\enskip\text{for}\enskip d\geqslant 2;\quad j=1,\dots ,d.
\end{equation}
Let the operator $\mathcal{Y}_2\colon L_2(\mathbb{R}^d;\mathbb{C}^n)\rightarrow L_2(\mathbb{R}^d;\mathbb{C}^{dn})$ act on the domain $\textrm{Dom}\,\mathcal{Y}_2=\textrm{Dom}\,\mathcal{X}$
by the rule $\mathcal{Y}_2\textbf{u}=\textrm{col}\,\lbrace a_1^*f\mathbf{u},\dots ,a_d^*f\mathbf{u}\rbrace $. Formally, 
$(\mathcal{Y}_2^*\mathcal{Y}+\mathcal{Y}^*\mathcal{Y}_2)\mathbf{u}=\sum_{j=1}^d\left( f^*a_jD_j(f\mathbf{u})+f^*D_j(a_j^*f\mathbf{u})\right)$.

Using the H\"older inequality, conditions  \eqref{f,h}, \eqref{a_j}, and the compactness of the embedding 
$H^1(\Omega )\hookrightarrow L_p(\Omega )$ for  $p=2\varrho 
/(\varrho -2)$, one can show that (cf. \cite[Subsection~5.2]{Su10}) for any $\nu >0$ there exist a constant $C(\nu)>0$ such that
\begin{equation}
\label{DOY_2u<=}
\Vert \mathcal{Y}_2\mathbf{u}\Vert ^2_{L_2(\mathbb{R}^d)}\leqslant
 \nu\Vert \mathcal{X}\mathbf{u}\Vert ^2_{L_2(\mathbb{R}^d)}+C(\nu)\Vert \mathbf{u}\Vert ^2_{L_2(\mathbb{R}^d)},\quad \mathbf{u}\in \textrm{Dom}\,\mathcal{X}.
\end{equation}
For $\nu$ fixed, the constant $C(\nu)$ depends on the norms $\Vert 
a_j\Vert _{L_\varrho (\Omega)},\;j=1,\dots ,d$, on $\Vert f\Vert 
_{L_\infty }$, $\Vert g^{-1}\Vert _{L_\infty }$, $\alpha _0$, $d$, 
$\varrho $, and on the parameters of the lattice $\Gamma $.

From \eqref{Yu<=} and \eqref{DOY_2u<=} one can derive the inequality
\begin{align}
\label{1089}
&2\varepsilon \vert \, \textrm{Re}(\mathcal{Y}\mathbf{u},\mathcal{Y}_2\mathbf{u})_{L_2(\mathbb{R}^d)}\vert \leqslant \frac{\kappa}{2}\Vert \mathcal{X}\mathbf{u}\Vert ^2_{L_2(\mathbb{R}^d)}+c_4\varepsilon ^2\Vert \mathbf{u}\Vert ^2_{L_2(\mathbb{R}^d)},
\\
&\mathbf{u}\in \textrm{Dom}\,\mathcal{X};
\quad
c_4=4\kappa ^{-1}c_1^2C(\nu)\enskip\text{for}\enskip\nu=\kappa ^2(16c_1^2)^{-1}.
\nonumber
\end{align}

\subsection{The operator $\mathcal{Q}_0$, the form $q[\mathbf{u},\mathbf{u}]$} 
\label{Subsection Q0 q}
Let $\mathcal{Q}_0$~be the operator acting in $L_2(\mathbb{R}^d;\mathbb{C}^n)$  
as multiplication by a $\Gamma$-periodic positive definite 
bounded matrix-valued function  
$\mathcal{Q}_0(\mathbf{x}):=f(\mathbf{x})^*f(\mathbf{x})$.

Suppose that $d\mu (\mathbf{x})=\lbrace d\mu 
_{jl}(\mathbf{x})\rbrace ,\,j,l=1,\dots ,n$ is a $\Gamma$-periodic Borel 
$\sigma $-finite measure in $\mathbb{R}^d$ with values in the class of
Hermitian $(n\times n)$-matrices. In other words, $d\mu _{jl}(\mathbf{x})$~is a
complex $\Gamma$-periodic measure in $\mathbb{R}^d$ and $d\mu 
_{jl}=d\mu _{lj}^*$. Assume that the measure $d\mu$ is such that the function 
$\vert v (\mathbf{x})\vert ^2$ is integrable with respect to each measure $d\mu _{jl}$ for any function $v\in H^1(\mathbb{R}^d)$.

In $L_2(\mathbb{R}^d;\mathbb{C}^n)$, consider the form
$$
q[\mathbf{u},\mathbf{u}]=\int\limits_{\mathbb{R}^d}\langle d\mu 
(\mathbf{x})f\mathbf{u},f\mathbf{u}\rangle ,\quad \mathbf{u}\in 
\textrm{Dom}\,\mathcal{X}.
$$ 
Assume that the measure $d\mu $ is subject to the following condition.

\begin{condition}
\label{Condition on mu}
\noindent
$1^\circ$. There exist constants $\widetilde{c}_2\geqslant 0$ and $\widehat{c}_3\geqslant 0$ such that for any $\mathbf{u},\mathbf{v}\in H^1(\Omega ;\mathbb{C}^n)$ we have
\begin{equation*}
\begin{split}
\Bigg\vert \int\limits_\Omega \langle d\mu (\mathbf{x})\mathbf{u},\mathbf{v}\rangle\Bigg\vert
&\leqslant
\left(\widetilde{c}_2\Vert \mathbf{D}\mathbf{u}\Vert ^2_{L_2(\Omega)}+\widehat{c}_3\Vert \mathbf{u}\Vert ^2 _{L_2(\Omega)}\right)^{1/2}
\\
&\times
\left(\widetilde{c}_2\Vert \mathbf{D}\mathbf{v}\Vert ^2_{L_2(\Omega)}+\widehat{c}_3\Vert \mathbf{v}\Vert ^2 _{L_2(\Omega)}\right)^{1/2}.
\end{split}
\end{equation*}

\noindent
$2^\circ$. We have
\begin{equation*}
\int\limits_\Omega \langle d\mu (\mathbf{x})\mathbf{u},\mathbf{u}\rangle
\geqslant
-\widetilde{c}\Vert \mathbf{D}\mathbf{u}\Vert ^2_{L_2(\Omega)}-\widehat{c}_0\Vert \mathbf{u}\Vert ^2_{L_2(\Omega)},\quad
\mathbf{u}\in H^1(\Omega ;\mathbb{C}^n),
\end{equation*}
with some constants $\widehat{c}_0\in\mathbb{R}$ and $\widetilde{c}$ such that $0\leqslant\widetilde{c}<\alpha _0\Vert g^{-1}\Vert ^{-1}_{L_\infty}$.
\end{condition}

Note that Condition~\ref{Condition on mu} implies estimates
{\allowdisplaybreaks
\begin{align}
\label{1034-1}
\begin{split}
\Bigg\vert \int\limits_\Omega \langle d\mu (\mathbf{x})f\mathbf{u},f\mathbf{v}\rangle\Bigg\vert
&\leqslant
\left(\widetilde{c}_2\Vert \mathbf{D}(f\mathbf{u})\Vert ^2_{L_2(\Omega)}+c_3\Vert \mathbf{u}\Vert ^2 _{L_2(\Omega)}\right)^{1/2}
\\
&\times
\left(\widetilde{c}_2\Vert \mathbf{D}(f\mathbf{v})\Vert ^2_{L_2(\Omega)}+c_3\Vert \mathbf{v}\Vert ^2 _{L_2(\Omega)}\right)^{1/2},
\end{split}
\\
\label{1034-2}
\int\limits_\Omega \langle d\mu (\mathbf{x})f\mathbf{u},f\mathbf{u}\rangle
&\geqslant
-\widetilde{c}\Vert \mathbf{D}(f\mathbf{u})\Vert ^2_{L_2(\Omega)}-c_0\Vert \mathbf{u}\Vert ^2_{L_2(\Omega)},
\end{align}
}
$\hspace{-4mm}$
$f\mathbf{u},f\mathbf{v}\in H^1(\Omega ;\mathbb{C}^n)$,  
with constant $c_0=\widehat{c}_0\Vert f\Vert _{L_\infty }^2$ if $\widehat{c}_0\geqslant 0$ and  $c_0=\widehat{c}_0\Vert f ^{-1}\Vert _{L_\infty}^{-2}$ if $\widehat{c}_0<0$; $c_3=\Vert f\Vert _{L_\infty }^2\widehat{c}_3$.

For $\mathbf{u},\mathbf{v}\in\textrm{Dom}\,\mathcal{X}$, we write inequalities 
\eqref{1034-1}, \eqref{1034-2} over shifted cells $\Omega +\mathbf{a},$ 
$\mathbf{a}\in\Gamma $, and sum up. We get
\begin{align*}
\begin{split}
\vert q[\mathbf{u},\mathbf{v}]\vert &\leqslant \left( \widetilde{c}_2\Vert \mathbf{D}(f\mathbf{u})\Vert _{L_2(\mathbb{R}^d)}^2 +c_3\Vert \mathbf{u} \Vert _{L_2(\mathbb{R}^d)}^2\right)^{1/2}
\\
&\times
\left( \widetilde{c}_2\Vert \mathbf{D}(f\mathbf{v})\Vert _{L_2(\mathbb{R}^d)}^2 +c_3\Vert \mathbf{v} \Vert _{L_2(\mathbb{R}^d)}^2\right)^{1/2},
\end{split}
\\
q[\mathbf{u},\mathbf{u}]&\geqslant -\widetilde{c}\Vert \mathbf{D}(f\mathbf{u})\Vert _{L_2(\mathbb{R}^d)}^2-c_0\Vert \mathbf{u}\Vert _{L_2(\mathbb{R}^d )}^2.
\end{align*}
By \eqref{<a[u,u]<}, this implies
\begin{align}
\label{DOcond.1.4-1}
\begin{split}
\vert q[\mathbf{u},\mathbf{v}]\vert
&\leqslant \left(c_2\Vert \mathcal{X}\mathbf{u}\Vert ^2_{L_2(\mathbb{R}^d)}+c_3\Vert \mathbf{u}\Vert _{L_2(\mathbb{R}^d)}^2\right)^{1/2}
\\
&\times
\left(c_2\Vert \mathcal{X}\mathbf{v}\Vert ^2_{L_2(\mathbb{R}^d)}+c_3\Vert \mathbf{v}\Vert _{L_2(\mathbb{R}^d)}^2\right)^{1/2},
\end{split}
\\
\label{DOcond.1.4-2}
q[\mathbf{u},\mathbf{u}]&\geqslant -(1-\kappa )\Vert \mathcal{X}\mathbf{u}\Vert _{L_2(\mathbb{R}^d)}^2-c_0\Vert \mathbf{u}\Vert ^2_{L_2(\mathbb{R}^d)},
\end{align}
$\mathbf{u},\mathbf{v}\in\textrm{Dom}\,\mathcal{X}$. Here 
$c_2=\widetilde{c}_2\alpha _0^{-1}\Vert g^{-1}\Vert _{L_\infty }$, 
$\kappa =1-\widetilde{c}\alpha _0^{-1}\Vert g^{-1}\Vert _{L_\infty }$, 
$0<\kappa \leqslant 1$.

Examples of measures $d\mu$ satisfying Condition~\ref{Condition on mu} were given in \cite[Subsection~5.5]{Su10}. Here we write only the main example.

\begin{example}
\label{Example Q in Ls}
Let the measure $d\mu$ be absolutely continuous with respect to Lebesgue measure: $d\mu (\mathbf{x})=\mathcal{Q}(\mathbf{x})\,d\mathbf{x}$, where $\mathcal{Q}(\mathbf{x})$ is a $\Gamma$-periodic Hermitian $(n\times n)$\nobreakdash-matrix-valued function in $\mathbb{R}^d$, and
\begin{equation}
\label{Q in Ls ex.}
\mathcal{Q}\in L_\sigma(\Omega),\quad \sigma =1\;\mbox{for}\;d=1;\quad \sigma >d/2\;\mbox{for}\;d\geqslant 2.
\end{equation}
Then
\begin{equation*}
q[\mathbf{u},\mathbf{u}]=\int\limits_{\mathbb{R}^d}\langle \mathcal{Q}(\mathbf{x})f(\mathbf{x})\mathbf{u}(\mathbf{x}),f(\mathbf{x})\mathbf{u}(\mathbf{x})\rangle\,d\mathbf{x},\quad
\mathbf{u}\in\mathrm{Dom}\,\mathcal{X}.
\end{equation*}
By the embedding theorem, under condition \eqref{Q in Ls ex.} for any $\nu >0$ there exists a positive constant $C_Q(\nu)$ such that
\begin{equation*}
\int\limits_\Omega \vert \mathcal{Q}(\mathbf{x})\vert \vert \mathbf{v}\vert ^2\,d\mathbf{x}
\leqslant
\nu \int\limits_\Omega \vert \mathbf{D}\mathbf{v}\vert ^2\,d\mathbf{x}
+C_Q(\nu)\int\limits_\Omega \vert \mathbf{v}\vert ^2\,d\mathbf{x},\quad
\mathbf{v}\in H^1(\Omega ;\mathbb{C}^n).
\end{equation*}
For $\nu$ fixed, the constant $C_Q(\nu)$ is controlled by $d$, $\sigma $, the norm $\Vert \mathcal{Q}\Vert _{L_\sigma(\Omega)}$ and the parameters of the lattice $\Gamma$. So, Condition~\textnormal{\ref{Condition on mu}} is satisfied and the constants can be chosen in the following way: $\widetilde{c}_2=1$, $\widehat{c}_3=C_Q(1)$, $\widetilde{c}=\nu$, and $\widehat{c}_0=C_Q(\nu)$ for $2\nu=\alpha _0\Vert g^{-1}\Vert ^{-1}_{L_\infty}$.
\end{example}

\subsection{The operator $\mathcal{B}(\varepsilon )$} 
\label{Subsection B(eps)}
In $L_2(\mathbb{R}^d;\mathbb{C}^n)$, we consider the quadratic form
\begin{equation}
\label{b(e)[u,u]}
\begin{split}
\mathfrak{b}(\varepsilon)[\mathbf{u},\mathbf{u}]
=\mathfrak{a}[\mathbf{u},\mathbf{u}]&+2\varepsilon\textrm{Re}\,(\mathcal{Y}\mathbf{u},\mathcal{Y}_2\mathbf{u})_{L_2(\mathbb{R}^d)}+\varepsilon ^2q[\mathbf{u},\mathbf{u}]
\\
&+\lambda \varepsilon ^2(\mathcal{Q}_0\mathbf{u},\mathbf{u})_{L_2(\mathbb{R}^d)},
\quad\mathbf{u}\in\textrm{Dom}\,\mathcal{X},
\end{split}
\end{equation}
where $0<\varepsilon \leqslant 1$, and the \textit{parameter $\lambda 
\in\mathbb{R}$ satisfies the following condition}:
\begin{equation}
\label{DOlambda}
\begin{split}
\lambda &>\Vert \mathcal{Q}_0 ^{-1}\Vert _{L_\infty }(c_0+c_4),\;\mbox{if}\;\lambda \geqslant 0,
\\
\lambda &>\Vert \mathcal{Q}_0 \Vert _{L_\infty}^{-1}(c_0+c_4),\;\mbox{if}\;\lambda <0\;(\mbox{and}\;c_0+c_4<0).
\end{split}
\end{equation}
Now we estimate the form \eqref{b(e)[u,u]} from below.
Let $\beta >0$ be defined by the formulas
\begin{equation}
\label{DObeta}
\begin{split}
\beta &=\lambda \Vert \mathcal{Q}_0^{-1}\Vert ^{-1}_{L_\infty}-c_0-c_4,\;\mbox{if}\;\lambda \geqslant 0,
\\
\beta &=\lambda \Vert \mathcal{Q}_0 \Vert _{L_\infty}-c_0-c_4,\;\mbox{if}\;\lambda <0\;(\mbox{and}\;c_0+c_4<0).
\end{split}
\end{equation}
Combining \eqref{1089}, \eqref{DOcond.1.4-2}, \eqref{DOlambda}, and 
\eqref{DObeta}, we arrive at
\begin{equation}
\label{b(e)[u,u]>=}
\mathfrak{b}(\varepsilon )[\mathbf{u},\mathbf{u}]\geqslant \frac{\kappa}{2}\mathfrak{a}[\mathbf{u},\mathbf{u}]+\varepsilon ^2\beta \Vert \mathbf{u}\Vert ^2_{L_2(\mathbb{R}^d)},
\quad \mathbf{u}\in \textrm{Dom}\, \mathcal{X}, \quad 0<\varepsilon \leqslant 1.
\end{equation}
Thus, the form $\mathfrak{b}(\varepsilon )$ is positive definite. 
Bringing together \eqref{Yu<=}, \eqref{DOY_2u<=} for $\nu =1$, and estimate  
\eqref{DOcond.1.4-1} for the quadratic form $q[\mathbf{u},\mathbf{u}]$, we obtain
\begin{equation}
\label{b(e)[u,u]<=}
\begin{split}
\mathfrak{b}(\varepsilon )[\mathbf{u},\mathbf{u}]&\leqslant (2+c_1^2+c_2)\mathfrak{a}[\mathbf{u},\mathbf{u}]\\
&\quad+\varepsilon ^2(C(1)+c_3+\vert \lambda \vert \Vert \mathcal{Q}_0 \Vert _{L_\infty })\Vert \mathbf{u}\Vert ^2_{L_2 (\mathbb{R}^d)},
\quad
\mathbf{u}\in\textrm{Dom}\,\mathcal{X}.
\end{split}
\end{equation}
Inequalities \eqref{b(e)[u,u]>=} and \eqref{b(e)[u,u]<=} imply that the form  
$\mathfrak{b}(\varepsilon )$ is closed.\textit{ The corresponding positive definite operator 
acting in the space $L_2(\mathbb{R}^d;\mathbb{C}^n)$ will be denoted by $\mathcal{B}(\varepsilon )$}. Formally,
\begin{equation}
\label{B(e)}
\begin{split}
\mathcal{B}(\varepsilon )&\!=\!\mathcal{A}\!+\!\varepsilon (\mathcal{Y}^*_2\mathcal{Y}\!+\!\mathcal{Y}^*\mathcal{Y}_2)\!+\!\varepsilon ^2f^*\mathcal{Q}f\!+\!\varepsilon ^2\lambda \mathcal{Q}_0
\\
&\!=\!f^*b(\mathbf{D})^*gb(\mathbf{D})f+\varepsilon 
\sum \limits _{j=1}^d f^*\left( a_jD_j+D_ja_j^* \right)f\!+\!\varepsilon ^2f^*\mathcal{Q}f\!+\!\varepsilon ^2\lambda \mathcal{Q}_0,
\end{split}
\end{equation}
where $\mathcal{Q}$ should be interpreted as a generalized matrix potential generated by the measure~$d\mu $.

\textit{For convenience of the further references, by the ,,initial data'' we call the following set of parameters:}
\begin{equation}
\label{dannie_zadachi}
\begin{split}
&d,\,m,\,n,\,\varrho ;\,\alpha _0,\,\alpha _1,\,\Vert g\Vert _{L_\infty},\,\Vert g^{-1}\Vert _{L_\infty} ,\,\Vert f\Vert _{L_\infty} ,\,\Vert f^{-1}\Vert _{L_\infty},\,\Vert a_j\Vert _{L_\varrho (\Omega )},\,
\\
&j=1,\dots ,d;\;
\widetilde{c},\,\widehat{c}_0,\,\widetilde{c}_2,\,\widehat{c}_3\,\mbox{ from Condition 4.1}; \;\lambda ; \mbox{ the parameters of the lattice }\Gamma .
\end{split}
\end{equation}
We will trace the dependence of constants in estimates on these data. 
The values of $c_1$, $C(1)$, $\kappa $, 
$c_2$, $c_3$, $c_4$, $c_0$, and $\beta$ are completely determined by initial data \eqref{dannie_zadachi}.

\section{Direct integral decomposition for the operator $\mathcal{B}(\varepsilon )$}

\subsection{The Gelfand transformation} 
The Gelfand transformation $\mathcal{U}$ is initially defined on functions from the Schwartz class $\mathbf{v}\in \mathcal{S}(\mathbb{R}^d;\mathbb{C}^n)$ 
by the rule
$$
\widetilde{\mathbf{v}}(\mathbf{k},\mathbf{x})=(\mathcal{U}\mathbf{v})(\mathbf{k},\mathbf{x})=\vert \widetilde{\Omega}\vert ^{-1/2}\sum\limits_{\mathbf{a}\in \Gamma}\textrm{exp}(-i\langle \mathbf{k},\mathbf{x}+\mathbf{a}\rangle )\mathbf{v}(\mathbf{x}+\mathbf{a}),\quad \mathbf{x}\in \Omega ,\; \mathbf{k}\in \widetilde{\Omega}.
$$
Moreover,
$\int _{\widetilde{\Omega}}\int _{\Omega}\vert 
\widetilde{\mathbf{v}}(\mathbf{k},\mathbf{x})\vert 
^2\,d\mathbf{x}\,d\mathbf{k}=\int _{\mathbb{R}^d}\vert 
\mathbf{v}(\mathbf{x})\vert ^2\,d\mathbf{x}$, so $\mathcal{U}$ can be extended by continuity to the unitary transformation 
\begin{equation*}
\mathcal{U}\; :\;L_2(\mathbb{R}^d;\mathbb{C}^n)
\rightarrow \int\limits_{\widetilde{\Omega}}\oplus L_2(\Omega;\mathbb{C}^n)\,d\mathbf{k}=:\mathcal{H}.
\end{equation*}
\textit{By $\widetilde{H}^1(\Omega;\mathbb{C}^n)$ we denote the subspace of functions from $H^1(\Omega;\mathbb{C}^n)$ whose 
$\Gamma$-periodic extension to  $\mathbb{R}^d$ belongs to the class 
$$
H_{\rm{loc}}^1(\mathbb{R}^d;\mathbb{C}^n).
$$
} The inclusion 
$\mathbf{v}\in H^1(\mathbb{R}^d;\mathbb{C}^n)$ is equivalent to
$$
\widetilde{\mathbf{v}}(\mathbf{k},\mathbf{\,\cdot\,})\in 
\widetilde{H}^1(\Omega;\mathbb{C}^n)
$$
 for a.~\!e. $\mathbf{k}\in 
\widetilde{\Omega}$ and
$$
\int\limits_{\widetilde{\Omega}}\int\limits_{\Omega}\left(\vert 
(\mathbf{D}+\mathbf{k)}\widetilde{\mathbf{v}}(\mathbf{k},\mathbf{x})\vert 
^2+\vert \widetilde{\mathbf{v}}(\mathbf{k},\mathbf{x})\vert 
^2\right)\,d\mathbf{x}\,d\mathbf{k}<\infty.
$$
Under the transformation $\mathcal{U}$, 
the operator of multiplication by a bounded periodic matrix-valued function in 
$L_2(\mathbb{R}^d;\mathbb{C}^n)$ turns into the operator of multiplication by the same function 
in the fibers of the direct integral   $\mathcal{H}$. 
The action of the operator $b(\mathbf{D})$ on the function $\mathbf{v}\in 
H^1(\mathbb{R}^d;\mathbb{C}^n)$ turns into the fiberwise action of the operator  
$b(\mathbf{D}+\mathbf{k})$ on 
$\widetilde{\mathbf{v}}(\mathbf{k},\mathbf{\,\cdot\,})\in 
\widetilde{H}^1(\Omega;\mathbb{C}^n)$.

\subsection{The operators $\mathcal{A}(\mathbf{k})$}
(See \cite[Subsection~2.2.1]{BSu}.) We set
\begin{equation}
\label{H, H_*}
\mathfrak{H}=L_2(\Omega;\mathbb{C}^n),\quad 
\mathfrak{H}_*=L_2(\Omega;\mathbb{C}^m),\quad 
\widetilde{\mathfrak{H}}=L_2(\Omega ;\mathbb{C}^{dn})
\end{equation}
and consider the closed operator
$\mathcal{X}(\mathbf{k})\, \colon \,\mathfrak{H}\rightarrow \mathfrak{H}_*,\; \mathbf{k}\in \mathbb{R}^d,$
defined by the relations
\begin{align}\label{X(k)}
&\mathcal{X}(\mathbf{k}) =hb(\mathbf{D}+\mathbf{k})f,\quad \mathbf{k}\in \mathbb{R}^d,
\\
\label{DomX(k)}
&\mathfrak{d} :=\textrm{Dom}\,\mathcal{X}(\mathbf{k})=\lbrace \mathbf{u}\in \mathfrak{H}\;:\;f\mathbf{u}\in \widetilde{H}^1(\Omega; \mathbb{C}^n)\rbrace.
\end{align}
The self-adjoint operator
$\mathcal{A}(\mathbf{k}):=\mathcal{X}(\mathbf{k})^*\mathcal{X}(\mathbf{k})\, 
\colon \,\mathfrak{H}\rightarrow \mathfrak{H},\;\mathbf{k}\in \mathbb{R}^d$, 
is generated by the closed quadratic form  
$\mathfrak{a}(\mathbf{k})[\mathbf{u},\mathbf{u}]:=\Vert 
\mathcal{X}(\mathbf{k})\mathbf{u}\Vert_{\mathfrak{H}_*}^2$, $\mathbf{u}\in 
\mathfrak{d}$, $\mathbf{k}\in \mathbb{R}^d$.
By \eqref{<b^*b<} and \eqref{f,h},
\begin{equation}
\label{<a(k)[u,u]<}
\begin{split}
\alpha _0 \Vert g^{-1}\Vert ^{-1}_{L_\infty } \Vert (\mathbf{D}+\mathbf{k})\mathbf{v}\Vert ^2_{L_2(\Omega)}\leqslant \mathfrak{a}(\mathbf{k})[\mathbf{u},\mathbf{u}]
\leqslant \alpha _1 \Vert g \Vert _{L_\infty} \Vert (\mathbf{D}+\mathbf{k})\mathbf{v}\Vert ^2_{L_2(\Omega)},
\\
\mathbf{v}=f\mathbf{u}\in \widetilde{H}^1(\Omega;\mathbb{C}^n).
\end{split}
\end{equation}
From \eqref{<a(k)[u,u]<} and the compactness of the embedding of $\widetilde{H}
^1(\Omega;\,\mathbb{C}^n)$ into $\mathfrak{H}$, it follows that the spectrum of the operator  
$\mathcal{A}(\mathbf{k})$ is discrete. Put 
$\mathfrak{N}:=\mathrm{Ker}\,\mathcal{A}(0)=\mathrm{Ker}\,\mathcal{X}(0)$. 
Inequalities \eqref{<a(k)[u,u]<} with $\mathbf{k}=0$ imply that
\begin{equation}
\label{KerA(k)=}
\mathfrak{N}=\mathrm{Ker}\,\mathcal{A}(0)=\lbrace \mathbf{u}\in L_2(\Omega;\mathbb{C}^n)\; :\; f\mathbf{u}=\mathbf{c}\in \mathbb{C}^n\rbrace ,\quad\mathrm{dim}\,\mathfrak{N}=n.
\end{equation}

As was shown in \cite[(2.2.11), (2.2.12)]{BSu}, 
\begin{equation}
\label{A(k)>}
\mathcal{A}(\mathbf{k})\geqslant c_*\vert \mathbf{k}\vert ^2I,\quad \mathbf{k}\in \mathrm{clos}\,\widetilde{\Omega};\quad c_*=\alpha _0\Vert f^{-1} \Vert_{L_\infty}^{-2} \Vert g^{-1}\Vert _{L_\infty}^{-1}.
\end{equation}
According to \cite[(2.2.14)]{BSu}, the distance $d^0$ from the point 
$\lambda _0 =0$ to the rest of the spectrum of the operator $\mathcal{A}(0)$ satisfies the estimate
\begin{equation}
\label{d^0>}
d^0\geqslant 4c_*r_0^2.
\end{equation}

\subsection{The operators $\mathcal{Y}(\mathbf{k})$ and $Y_2$} 
Consider the operator 
$\mathcal{Y}(\mathbf{k})\colon \mathfrak{H}\rightarrow 
\widetilde{\mathfrak{H}}$ acting on the domain  
$\textrm{Dom}\,\mathcal{Y}(\mathbf{k})=\mathfrak{d}$ by the rule 
\begin{equation}
\label{Y(k)u}
\mathcal{Y}(\mathbf{k})\mathbf{u}=(\mathbf{D}+\mathbf{k})f\mathbf{u}
=\mathrm{col}\,\lbrace (D_1+k_1)f\mathbf{u},\dots ,(D_d+k_d)f\mathbf{u}\rbrace ,\quad\mathbf{u}\in\mathfrak{d} .
\end{equation}
From the lower estimate \eqref{<a(k)[u,u]<}, it follows that
\begin{equation}
\label{Y(k)u<=}
\Vert \mathcal{Y}(\mathbf{k})\mathbf{u}\Vert _\mathfrak{H}\leqslant c_1\Vert \mathcal{X}(\mathbf{k})\mathbf{u}\Vert _{\mathfrak{H}_*},\quad\mathbf{u}\in\mathfrak{d},
\end{equation}
where the constant $c_1$ was defined in \eqref{Yu<=}.

Consider the operator $Y_2:\mathfrak{H}\rightarrow \widetilde{\mathfrak{H}}$ defined by the relation
\begin{equation}
\label{Y_2u=}
Y_2\mathbf{u}=\textrm{col}\,\lbrace a_1^*f\mathbf{u},\dots ,a_d^*f\mathbf{u}\rbrace ,\quad \mathrm{Dom}\,Y_2=\mathfrak{d}.
\end{equation}

As was shown in \cite[Subsection 5.7]{Su10}, for any $\nu >0$ there exist constants  
$C_j(\nu )>0$, $j=1,\dots ,d$, such that for $\mathbf{k}\in 
\mathbb{R}^d$ we have
\begin{equation*}
\Vert a_j^*\mathbf{v}\Vert ^2_{L_2(\Omega )}\leqslant \nu \Vert (\mathbf{D}+\mathbf{k})\mathbf{v}\Vert ^2_{L_2(\Omega )}+C_j(\nu )\Vert \mathbf{v}\Vert ^2_{L_2(\Omega )},\; \mathbf{v}\in \widetilde{H}^1(\Omega ;\mathbb{C}^n),\;j=1,\dots ,d.
\end{equation*}
Let $\mathbf{v}=f\mathbf{u}$, $\mathbf{u}\in \mathfrak{d}$. Then, summing up the mentioned inequalities  over
$j$ and taking \eqref{f,h} and  
\eqref{<a(k)[u,u]<} into account, we obtain that for any $\nu >0$ there exists a constant $C(\nu)>0$ (the same as in \eqref{DOY_2u<=}) such that
\begin{equation}
\label{DOcond.1.3}
\Vert Y_2\mathbf{u}\Vert ^2_{\widetilde{\mathfrak{H}}}\leqslant \nu \Vert \mathcal{X}(\mathbf{k})\mathbf{u}\Vert ^2_{\mathfrak{H}_*}+C(\nu)\Vert \mathbf{u}\Vert ^2_\mathfrak{H},\quad \mathbf{u}\in \mathfrak{d},\;\mathbf{k}\in \mathbb{R}^d.
\end{equation}

\subsection{The operator $\mathcal{Q}_0$ and the form $q_\Omega [\mathbf{u},\mathbf{u}]$} 
Let $\mathcal{Q}_0$~be the bounded operator in $\mathfrak{H}$, acting as multiplication by the matrix-valued function 
$\mathcal{Q}_0(\mathbf{x})=f(\mathbf{x})^*f(\mathbf{x})$.

In $L_2(\Omega ;\mathbb{C}^n)$, consider the form 
$$
q_\Omega 
[\mathbf{u},\mathbf{u}]=\int\limits_\Omega \langle d\mu 
(\mathbf{x})f\mathbf{u},f\mathbf{u}\rangle ,\quad \mathbf{u}\in 
\mathfrak{d}.$$ 
In \eqref{1034-1} and \eqref{1034-2}, we replace the function 
$f(\mathbf{x})\mathbf{u}(\mathbf{x})$ by
$f(\mathbf{x})\mathbf{u}(\mathbf{x})\exp (i\langle 
\mathbf{k},\mathbf{x}\rangle )$ (these functions simultaneously belong to the space $H^1(\Omega 
;\mathbb{C}^n)$) and $f(\mathbf{x})\mathbf{v}(\mathbf{x})$ by $f(\mathbf{x})\mathbf{v}(\mathbf{x})\exp (i\langle 
\mathbf{k},\mathbf{x}\rangle )$. By using \eqref{<a(k)[u,u]<}, we get
\begin{align}
\label{<q_Omega<-1}
\vert q_\Omega [\mathbf{u},\mathbf{v}]\vert &\leqslant
\left( c_2\Vert \mathcal{X}(\mathbf{k})\mathbf{u}\Vert 
^2_{\mathfrak{H}_*}+c_3\Vert \mathbf{u}\Vert ^2_\mathfrak{H}\right)^{1/2}
\left( c_2\Vert \mathcal{X}(\mathbf{k})\mathbf{v}\Vert 
^2_{\mathfrak{H}_*}+c_3\Vert \mathbf{v}\Vert ^2_\mathfrak{H}\right)^{1/2},
\\
\label{<q_Omega<-2}
q_\Omega [\mathbf{u},\mathbf{u}]&\geqslant -(1-\kappa )\Vert \mathcal{X}(\mathbf{k})\mathbf{u}\Vert ^2_{\mathfrak{H}_*}-c_0\Vert \mathbf{u}\Vert ^2_\mathfrak{H};\quad \mathbf{u}, \mathbf{v}\in \mathfrak{d},\quad\mathbf{k}\in \mathbb{R}^d .
\end{align}
Here the constants $\kappa$, $c_0$, $c_2$, and $c_3$ are the same as in \eqref{DOcond.1.4-1} and \eqref{DOcond.1.4-2}.

\subsection{The operator pencil $\mathcal{B}(\mathbf{k},\varepsilon )$} 
In $\mathfrak{H}$, consider the quadratic form 
\begin{align*}
\mathfrak{b}(\mathbf{k},\varepsilon )[\mathbf{u},\mathbf{u}]
=\mathfrak{a}(\mathbf{k})[\mathbf{u},\mathbf{u}]&+2\varepsilon\textrm{Re}\,(\mathcal{Y}(\mathbf{k})\mathbf{u},Y_2\mathbf{u})_{\widetilde{\mathfrak{H}}}
\\
&+\varepsilon ^2 q_\Omega [\mathbf{u},\mathbf{u}]+\lambda \varepsilon ^2(\mathcal{Q}_0\mathbf{u},\mathbf{u})_\mathfrak{H},
\quad \mathbf{u}\in\mathfrak{d}.
\end{align*}
From \eqref{DOlambda}, \eqref{DObeta},
\eqref{Y(k)u<=}, \eqref{DOcond.1.3}, and \eqref{<q_Omega<-2} it follows that
\begin{equation}
\label{b(k,e)[u,u]>=}
\mathfrak{b}(\mathbf{k},\varepsilon )[\mathbf{u},\mathbf{u}]\geqslant \frac{\kappa}{2}\mathfrak{a}(\mathbf{k})[\mathbf{u},\mathbf{u}]+\beta \varepsilon ^2\Vert \mathbf{u}\Vert ^2_\mathfrak{H},\quad \mathbf{u}\in\mathfrak{d}.
\end{equation}
Next, from \eqref{Y(k)u<=}, \eqref{DOcond.1.3} for $\nu =1$, and 
 \eqref{<q_Omega<-1} we derive that
\begin{equation}
\label{b(k,e)[u,u]<=}
\begin{split}
\mathfrak{b}(\mathbf{k},\varepsilon )[\mathbf{u},\mathbf{u}]&\leqslant (2+c_1^2+c_2)\mathfrak{a}(\mathbf{k})[\mathbf{u},\mathbf{u}]
\\
&\quad+(C(1)+c_3+\vert \lambda \vert \Vert \mathcal{Q}_0\Vert _{L_\infty })\varepsilon ^2\Vert \mathbf{u}\Vert ^2_\mathfrak{H},\quad \mathbf{u}\in \mathfrak{d}.
\end{split}
\end{equation}
Inequalities \eqref{b(k,e)[u,u]>=} and \eqref{b(k,e)[u,u]<=} show that the form $\mathfrak{b}(\mathbf{k},\varepsilon )$ is positive definite and closed on the domain 
\eqref{DomX(k)}. The self-adjoint positive-definite operator in  $\mathfrak{H}$ corresponding to this form is denoted by  
$\mathcal{B}(\mathbf{k},\varepsilon )$. Formally we can write
\begin{equation}\label{B(k,e)}
\begin{split}
\mathcal{B}(\mathbf{k},\varepsilon )&\!=\!\mathcal{A}(\mathbf{k})\!+\!\varepsilon (Y_2^*\mathcal{Y}(\mathbf{k})\!+\!\mathcal{Y}(\mathbf{k})^*Y_2)\!+\!\varepsilon ^2 f^*\mathcal{Q}f\!+\!\lambda \varepsilon ^2 \mathcal{Q}_0
\\
&\!=\!f^*b(\mathbf{D}\!+\!\mathbf{k})^*gb(\mathbf{D}\!+\!\mathbf{k})f
\!+\!\varepsilon \sum \limits _{j=1}^d f^*\big(a_j(D_j\!+\!k_j)\!+\!(D_j\!+\!k_j)a_j^*\big)f\\
&\qquad\qquad\qquad\qquad\qquad\quad\!+\!\varepsilon ^2f^*\mathcal{Q}f\!+\!\lambda \varepsilon ^2 f^*f.
\end{split}
\end{equation}

\subsection{The direct integral decomposition for the operator $\mathcal{B}(\varepsilon)$} 
\label{Subsection direct integral decomposition}
Under the Gelfand transformation $\mathcal{U}$, the operator 
\eqref{B(e)} acting in the space $L_2(\mathbb{R}^d;\mathbb{C}^n)$ is decomposed into the direct integral of operators~\eqref{B(k,e)} acting in 
$L_2(\Omega ;\mathbb{C}^n)$:
$$
\mathcal{U}\mathcal{B}(\varepsilon )\mathcal{U}^{-1}=\int\limits
_{\widetilde{\Omega }}\oplus \mathcal{B}(\mathbf{k},\varepsilon 
)\,d\mathbf{k}.
$$ 
This means the following. Let  
$\widetilde{\mathbf{u}}=\mathcal{U}\mathbf{u},$ where $\mathbf{u}\in 
\mathrm{Dom}\,\mathfrak{b}(\varepsilon)=\mathrm{Dom}\,\mathcal{X}$. Then
\begin{align}
\label{tilde u}
&\widetilde{\mathbf{u}}(\mathbf{k},\,\cdot\,)\in \mathfrak{d}\; \mbox{for a.~e. } \mathbf{k}\in \widetilde{\Omega},
\\
\label{b[u,u]=int}
&\mathfrak{b}(\varepsilon )[\mathbf{u},\mathbf{u}]=\int\limits_{\widetilde{\Omega}}\mathfrak{b}(\mathbf{k},\varepsilon 
)[\widetilde{\mathbf{u}}(\mathbf{k},\,\cdot\, 
),\,\widetilde{\mathbf{u}}(\mathbf{k},\,\cdot\, )]\,d\mathbf{k}.
\end{align}
Conversely, if for $\widetilde{\mathbf{u}}\in \mathcal{H}$ one has
\eqref{tilde u} and the integral in \eqref{b[u,u]=int} is finite, then
$\mathbf{u}\in \mathrm{Dom}\,\mathfrak{b}(\varepsilon )=\mathrm{Dom}\,\mathcal{X}$ and  
\eqref{b[u,u]=int} is fulfilled.

\section{Incorporation of the operators $\mathcal{B}(\mathbf{k},\varepsilon )$ \\
into the abstract scheme}

The content of Subsections~\ref{Subsection vkluchenie-1}--\ref{Subsection N(k,eps)} is borrowed from \cite{Su14}.

\subsection{} 
\label{Subsection vkluchenie-1}
For $d>1$, the operator $\mathcal{B}(\mathbf{k},\varepsilon )$ depends on the multi-dimensional parameter~$\mathbf{k}$. According to \cite[Chapter~2]{BSu}, we distinct a one-dimensional parameter $t$, setting $\mathbf{k}=t\boldsymbol{\theta}$, 
$t=\vert \mathbf{k}\vert $, $\boldsymbol{\theta}\in\mathbb{S}^{d-1}.$ 
We apply the scheme of Chapter~1.  In this case, the resulting objects depend on an additional parameter $\boldsymbol{\theta}$ and we must make our considerations and estimates uniform in $\boldsymbol{\theta}$. 
The spaces $\mathfrak{H}$, $\mathfrak{H}_*$, and $\widetilde{\mathfrak{H}}$ 
are defined in \eqref{H, H_*}. Put $X(t)=X(t;\boldsymbol{\theta} 
):=\mathcal{X}(t\boldsymbol{\theta}).$ By~\eqref{X(k)}, we have
$X(t;\boldsymbol{\theta} )=X_0+tX_1(\boldsymbol{\theta} ),$ where
\begin{equation}
\label{6.0}
X_0=\mathcal{X}(0)=h b(\mathbf{D})f,\quad 
\mathrm{Dom}\,X_0=\mathfrak{d};\quad X_1(\boldsymbol{\theta })=h b(\boldsymbol{\theta})f .
\end{equation}
Next, we set 
$A(t)=A(t;\boldsymbol{\theta}):=\mathcal{A}(t\boldsymbol{\theta} )$. In accordance with~\eqref{KerA(k)=}, the kernel 
$\mathfrak{N}=\mathrm{Ker}\,X_0=\mathrm{Ker}\,\mathcal{A}(0)$ in $n$-dimensional. 
Condition~\ref{Condition lambda0} is satisfied. The quantity $d^0$ admits estimate \eqref{d^0>}. As was shown in \cite[Chapter~2, \S 3]{BSu}, the condition $n\leqslant 
n_*=\mathrm{dim\,Ker}\,X_0^*$ is also satisfied. 
Estimate \eqref{A(k)>} corresponds to inequality \eqref{A(t)>= abstract}.

Next, the role of $Y(t)$ is played by the operator $Y(t;\boldsymbol{\theta 
}):=\mathcal{Y}(t\boldsymbol{\theta})$. By \eqref{Y(k)u}, we have the identity
$Y(t;\boldsymbol{\theta })=Y_0+tY_1(\boldsymbol{\theta })$, where
\begin{equation}
\label{Y_0, Y_1(theta)}
\begin{split}
&Y_0\mathbf{u}=\mathbf{D}(f\mathbf{u})=\textrm{col}\,\lbrace D_1f\mathbf{u},\dots ,D_df\mathbf{u}\rbrace ,\quad \mathrm{Dom}\,Y_0=\mathfrak{d};
\\
&Y_1(\boldsymbol{\theta })\mathbf{u}=\mathrm{col}\,\lbrace \theta 
_1f\mathbf{u},\dots ,\theta _df\mathbf{u}\rbrace .
\end{split}
\end{equation}
Condition~\ref{Condition Y(t)} holds true because of estimate \eqref{Y(k)u<=}. The operator $Y_2$ is defined in~\eqref{Y_2u=}. Condition~\ref{Condition Y2} is satisfied due to \eqref{DOcond.1.3}. The role of the form $\mathfrak{q}$ from Subsection~\ref{Subsection form g abstract scheme} is played by the form $q_\Omega$. Condition~\ref{Condition q}  
is valid due to \eqref{<q_Omega<-1}, \eqref{<q_Omega<-2}. The role of the operator $Q_0$ from Subsection~\ref{Subsection B(t,eps) abstract} is played by the operator of multiplication by the matrix-valued function  $\mathcal{Q}_0(\mathbf{x})$. 
The restriction \eqref{lambda condition} on the parameter $\lambda $ holds true by 
\eqref{DOlambda}. 

Finally, the role of the operator pencil $B(t,\varepsilon )$ (see
Subsection~\ref{Subsection B(t,eps) abstract}) is played by the operator family \eqref{B(k,e)}: 
$B(t,\varepsilon ;\boldsymbol{\theta }):=\mathcal{B}(t\boldsymbol{\theta 
},\varepsilon )$.

Thus, all the assumptions of the abstract scheme are fulfilled.

\subsection{} 
According to Subsection~\ref{Subsection to parameters tau and vartheta}, we should fix a positive number~$\delta $ such that 
$\delta <\kappa d^0/13$. Taking \eqref{A(k)>} and \eqref{d^0>} into account,  
we put 
\begin{equation*}
\delta = \frac{1}{4}\kappa c_*r_0^2=\frac{1}{4}\kappa \alpha _0\Vert f^{-1}\Vert ^{-2}_{L_\infty }\Vert g ^{-1}\Vert ^{-1}_{L_\infty }r_0^2.
\end{equation*}
Note that \eqref{<b^*b<}, \eqref{f,h}, \eqref{6.0}, and \eqref{Y_0, Y_1(theta)} imply the relations
\begin{equation}
\label{X_1(theta)<=, Y_1(theta)=}
\Vert X_1(\boldsymbol{\theta } )\Vert \leqslant \alpha _1^{1/2}\Vert g\Vert ^{1/2}_{L_\infty }\Vert f \Vert _{L_\infty },\quad \Vert Y_1(\boldsymbol{\theta })\Vert =\Vert f\Vert _{L_\infty },\quad \boldsymbol{\theta }\in \mathbb{S}^{d-1}.
\end{equation}

Now the right-hand side in estimate \eqref{tau 0 in abstract scheme} takes the form 
$$\delta ^{1/2}((2+c_1^2+c_2)\Vert X_1(\boldsymbol{\theta })\Vert ^2 +C(1)+c_3+\vert \lambda \vert \Vert f\Vert ^2_{L_\infty })^{-1/2}$$ 
and depends on $\boldsymbol{\theta }$. We choose the following value for the constant $\tau _0$ suitable for all $\boldsymbol{\theta }\in \mathbb{S}^{d-1}$:
\begin{equation}
\label{DOtau_0}
\tau _0=\delta ^{1/2}((2+c_1^2+c_2)\alpha _1\Vert g \Vert _{L_\infty }\Vert f \Vert ^2_{L_\infty }+C(1)+c_3+\vert \lambda \vert \Vert f \Vert ^2_{L_\infty })^{-1/2}.
\end{equation}

By \eqref{A(k)>} and 
\eqref{b(k,e)[u,u]>=}, the operator $B(t,\varepsilon ;\boldsymbol{\theta 
})$ satisfies a condition of the form~\eqref{B(tau,theta)>=}:
\begin{align}
\label{B(t,eps,theta)>=}
\mathcal{B}(\mathbf{k},\varepsilon )&=B(t,\varepsilon ;\boldsymbol{\theta }) \geqslant 
\check{c}_*(t^2+\varepsilon ^2)I,\quad \mathbf{k}=t\boldsymbol{\theta 
}\in\widetilde{\Omega},\quad 0<\varepsilon \leqslant 1;
\\
\label{DOcheck c_*}
\check{c}_* &=\frac{1}{2}\min \lbrace \kappa c_*;2\beta \rbrace.
\end{align}

\subsection{The case when $f=\mathbf{1}_n$}
\label{Subsection f=1n}

\textit{In what follows, all the objects corresponding to the case of $f=\mathbf{1}_n$, 
are marked by the upper sing} ``$\widehat{\quad}$''. We have
$\widehat{\mathfrak{H}}=\mathfrak{H}=L_2(\Omega;\mathbb{C}^n)$. 
In accordance with Subsection~\ref{Subsection vkluchenie-1}, 
$\widehat{X}(t;\boldsymbol{\theta} 
)=\widehat{X}_0+t\widehat{X}_1(\boldsymbol{\theta} )$, where
\begin{equation}
\label{hat X0 for DO's}
\widehat{X}_0=h b(\mathbf{D}),\quad 
\mathrm{Dom}\,\widehat{X}_0=\widetilde{H}^1(\Omega;\mathbb{C}^n),
\end{equation}
 and $\widehat{X}_1(\boldsymbol{\theta })$~is a bounded operator of multiplication by the matrix 
$h (\mathbf{x}) b(\boldsymbol{\theta})$:
\begin{equation}
\label{X_1 hat for DO's}
\widehat{X}_1(\boldsymbol{\theta })=h b(\boldsymbol{\theta}).
\end{equation}
Formally, 
$\widehat{A}(t;\boldsymbol{\theta})=\widehat{X}(t;\boldsymbol{\theta})^*\widehat{X}(t;\boldsymbol{\theta})$. 
In the case when $f=\mathbf{1}_n$, the kernel \eqref{KerA(k)=} coincides with the subspace of constants $\widehat{\mathfrak{N}}=\lbrace \mathbf{u}\in 
\mathfrak{H}\; :\;\mathbf{u}=\mathbf{c}\in \mathbb{C}^n\rbrace $. 
The orthogonal projection $\widehat{P}$ of the space 
$\mathfrak{H}=L_2(\Omega;\mathbb{C}^n)$ onto the subspace  
$\widehat{\mathfrak{N}}=\mathbb{C}^n$~is the operator of averaging over the cell 
$\Omega$: 
\begin{equation}
\label{P= for DO's}
\widehat{P}\mathbf{u}=\vert \Omega \vert ^{-1}\int\limits_{\Omega}\mathbf{u}(\mathbf{x})\,d\mathbf{x}.
\end{equation}

Next, 
$\widehat{Y}(t;\boldsymbol{\theta})=\widehat{Y}_0+t\widehat{Y}_1(\boldsymbol{\theta}): 
{\mathfrak{H}}\rightarrow \widetilde{\mathfrak{H}}$, where 
\begin{align}
\label{Y0 for DO's}
\widehat{Y}_0\mathbf{u}&=\mathbf{Du}=\mathrm{col}\, \lbrace D_1\mathbf{u},\dots ,D_d\mathbf{u} \rbrace,
\quad
\mathrm{Dom}\,\widehat{Y}_0=\widetilde{H}^1(\Omega ;\mathbb{C}^n),
\\
\label{hat Y1 for DO's}
\widehat{Y}_1(\boldsymbol{\theta })\mathbf{u}&=\mathrm{col}\,\lbrace \theta _1\mathbf{u},\dots ,\theta _d\mathbf{u}\rbrace.
\end{align}
The operator  $\widehat{Y}_2: {\mathfrak{H}}\rightarrow 
\widetilde{\mathfrak{H}}$ acts by the rule 
\begin{equation}
\label{hat Y2 for DO's}
\widehat{Y}_2\mathbf{u}=\mathrm{col}\,\lbrace a_1^*\mathbf{u},\dots 
,a_d^*\mathbf{u}\rbrace,\quad  \mathrm{Dom}\,\widehat{Y}_2=\widetilde{H}^1(\Omega ;\mathbb{C}^n). 
\end{equation}
The role of the form 
$\widehat{\mathfrak{q}}[\mathbf{u},\mathbf{u}]$ is played by the form $\int _\Omega 
\langle d\mu (\mathbf{x})\mathbf{u},\mathbf{u}\rangle$, $\mathbf{u}\in\widetilde{H}^1(\Omega ;\mathbb{C}^n)$. The role of the operator  
$\widehat{Q}_0$ is played by the identity operator~$I$.

The operator pencil $\widehat{B}(t,\varepsilon ;\boldsymbol{\theta })$ 
is formally given by the expression 
$$
\widehat{B}(t,\varepsilon ;\boldsymbol{\theta })=\widehat{A}(t;\boldsymbol{\theta })+\varepsilon (\widehat{Y}_2^*\widehat{Y}(t;\boldsymbol{\theta })+\widehat{Y}(t;\boldsymbol{\theta })^*\widehat{Y}_2)+\varepsilon ^2 \mathcal{Q}+\lambda\varepsilon ^2I.
$$
(We emphasise that the formal object $\mathcal{Q}$ is not changed.)

\subsection{The case when $f\neq \mathbf{1}_n$. The implementation of the assumptions of \S\ref{Section bordered exp abstract} }

Now we return to consideration of the operators $\mathcal{B}(\varepsilon)$ of the general form \eqref{B(e)} and the corresponding families $B(t,\varepsilon 
;\boldsymbol{\theta})$ defined in Subsection~\ref{Subsection vkluchenie-1}. We hold the upper mark ``$\widehat{\quad}$'' to indicate the objects corresponding to the case where 
$f=\mathbf{1}_n$ and 
$b(\mathbf{D}),\,g,\,a_j,\,j=1,\dots,d,\,\lambda$, and $d\mu$ are the same as before.

Let us show that the operator families $B(t,\varepsilon;\boldsymbol{\theta})$ and $\widehat{B}(t,\varepsilon ;\boldsymbol{\theta })$ satisfy the assumptions from \S\ref{Section bordered exp abstract} of the abstract scheme. Indeed, identity \eqref{B=M*hat-B M abstract} corresponds to the obvious equality $B(t,\varepsilon;\boldsymbol{\theta})=f^*\widehat{B}(t,\varepsilon;\boldsymbol{\theta})f$. The role of the isomorphism $M$ is played now by the operator of multiplication by the matrix-valued function $f(\mathbf{x})$. Then the role of the operator $G$ (see \eqref{G= abstract scheme}) is played by the operator of multiplication by the matrix-valued function $(f(\mathbf{x})f(\mathbf{x})^*)^{-1}$. 
The block of $G$ in the kernel $\widehat{\mathfrak{N}}=\mathbb{C}^n$~is the operator of multiplication by the constant matrix 
$\overline{G}=\vert \Omega\vert ^{-1}\int 
_{\Omega}(f(\mathbf{x})f(\mathbf{x})^*)^{-1}\,d\mathbf{x}$. The role of the operator  
$M_0$ (see \eqref{M0}) is played by the operator of multiplication by the matrix  
\begin{equation}
\label{f0}
f_0:=(\overline{G})^{-1/2}=(\underline{ff^*})^{1/2}.
\end{equation}
Note that 
\begin{equation}
\label{f_0<f}
\vert f_0\vert  \leqslant \Vert f\Vert _{L_\infty},\quad \vert f_0^{-1}\vert \leqslant \Vert f^{-1}\Vert _{L_\infty}.
\end{equation}

\subsection{On the constants in estimates}

Our aim is to apply the results of Chapter~1 to the operator $B(t,\varepsilon;\boldsymbol{\theta})$ depending on the additional parameter $\boldsymbol{\theta}$. To realise the abstract scheme, we need to make all the estimates from Chapter~1 for this operator uniform in~$\boldsymbol{\theta}$, i. e., for each constant to choose a value suitable for all $\boldsymbol{\theta}\in\mathbb{S}^{d-1}$. In Chapter~1, it was traced that (see Remarks \ref{Remark constants Sec.1} and \ref{Remark CG} and Theorem~\ref{Theorem bordered exp final in abstrach scheme}) the values of the constants $C_N$, $C_5$, $C_7$, $C_8$, and $C_9$ are controlled in terms of polynomials with positive numerical coefficients and  the variables \eqref{data for abstract scheme} and $\check{c}_*^{-1}$, and the constant $C_G$ is controlled in terms of the same variables and on $\Vert M ^{-1}\Vert$. In the case under consideration, $\Vert M ^{-1}\Vert =\Vert \mathcal{Q}_0\Vert ^{1/2} _{L_\infty}$, the data \eqref{data for abstract scheme} reduce to the set of data $\delta ,\delta ^{-1/2}$, $\tau _0$, $\kappa ^{1/2}$, $\kappa ^{-1/2}$, $c_1$, $c_2^{1/2}$, $c_3^{1/2}$, $C(1)^{1/2}$, $\vert \lambda\vert$, $\Vert X_1(\boldsymbol{\theta})\Vert $, $\Vert Y_1(\boldsymbol{\theta})\Vert $, and $ \Vert \mathcal{Q}_0\Vert _{L_\infty}$. The quantities $\delta ,\delta ^{-1/2}$, $\kappa ^{1/2}$, $\kappa ^{-1/2}$, $c_1$, $c_2^{1/2}$, $c_3^{1/2}$, and $C(1)^{1/2}$ do not depend on $\boldsymbol{\theta}$ and can be controlled in terms of the initial data \eqref{dannie_zadachi}; the number $\tau _0$ is already chosen independent on $\boldsymbol{\theta}$ (see \eqref{DOtau_0}). By using \eqref{X_1(theta)<=, Y_1(theta)=}, instead of $\Vert X_1(\boldsymbol{\theta})\Vert $ we can take $\alpha _1^{1/2}\Vert g\Vert ^{1/2}_{L_\infty}\Vert f\Vert _{L_\infty}$, and instead of $\Vert Y_1(\boldsymbol{\theta})\Vert $ we can take the norm $\Vert f\Vert _{L_\infty}$. Thus, the following observation holds true. 

\begin{remark}
\label{Remark Constants Chapter 2}
The values of constants $C_G$, $C_5$, $C_7$, $C_8$, and $C_9$ for the operator $B(t,\varepsilon;\boldsymbol{\theta})$ can be chosen independent on $\boldsymbol{\theta}\in\mathbb{S}^{d-1}$.
\end{remark}

\section{The effective characteristics} 

\subsection{The operators $\widehat{Z}(\boldsymbol{\theta})$, $\widehat{\widetilde{Z}}$, and $\widehat{R}(\boldsymbol{\theta})$}

The operator $\widehat{Z}$ defined in Subsection~\ref{Subsection Z and tilde-Z-abstract} now depends on $\boldsymbol{\theta}$. We define a $\Gamma$-periodic $(n\times m)$-matrix-valued function~$\Lambda (\mathbf{x})$ as the weak solution of the equation 
\begin{equation}
\label{Lambda}
b(\mathbf{D})^*g(\mathbf{x})(b(\mathbf{D})\Lambda (\mathbf{x})+\mathbf{1}_m)=0,\quad \int\limits_\Omega \Lambda (\mathbf{x})\,d\mathbf{x}=0.
\end{equation}
As was shown in \cite[Subsection 6.3]{Su10}, the operator $\widehat{Z}(\boldsymbol{\theta}):\mathfrak{H}\rightarrow\mathfrak{H}$ is defined by the equality
\begin{equation}
\label{Z= for DO's}
\widehat{Z}(\boldsymbol{\theta})=\Lambda b(\boldsymbol{\theta})\widehat{P},
\end{equation}
where $\widehat{P}$ is the projector \eqref{P= for DO's}.

According to \cite[Subsection 6.3]{Su10} 
\begin{equation}
\label{tilde Z = for DO's}
\widehat{\widetilde{Z}}=\widetilde{\Lambda}\widehat{P},
\end{equation} 
where $\widetilde{\Lambda }(\mathbf{x})$~is a $\Gamma $-periodic
$(n\times n)$-matrix-valued solution of the problem
\begin{equation}
\label{tildeLambda}
b(\mathbf{D})^*g(\mathbf{x})b(\mathbf{D})\widetilde{\Lambda }(\mathbf{x})+\sum \limits _{j=1}^dD_ja_j(\mathbf{x})^*=0,\quad \int\limits_{\Omega }\widetilde{\Lambda }(\mathbf{x})\,d\mathbf{x}=0.
\end{equation}

The operator $\widehat{R}(\boldsymbol{\theta}):\widehat{\mathfrak{N}}\rightarrow\mathfrak{N}_*$ acts as multiplication by the matrix-valued function 
\begin{equation}
\label{hat R for DO's}
\widehat{R}(\boldsymbol{\theta})=h(b(\mathbf{D})\Lambda +\mathbf{1}_m)b(\boldsymbol{\theta}).
\end{equation}

\subsection{The operator $\widehat{S}(\boldsymbol{\theta})$. The effective matrix}

The spectral germ~$\widehat{S}$ introduced in Subsection~\ref{Subsection R and S in abstract scheme} now depends on  
$\boldsymbol{\theta}$. According to \cite[Chapter~3,~\S1]{BSu}, the operator 
$\widehat{S}(\boldsymbol{\theta}) :\widehat{\mathfrak{N}}\rightarrow 
\widehat{\mathfrak{N}}$ acts as the operator of multiplication by the matrix  
$b(\boldsymbol{\theta})^*g^0b(\boldsymbol{\theta})$, 
$\boldsymbol{\theta}\in \mathbb{S}^{d-1}.$ Here $g^0$~is a constant  
$(m\times m)$-matrix called the \textit{effective matrix} and defined by the expression 
\begin{equation}
\label{g^0}
g^0=\vert \Omega \vert ^{-1}\int\limits_{\Omega}\widetilde{g}(\mathbf{x})\,d\mathbf{x} ,\quad  \widetilde{g}(\mathbf{x}):= g(\mathbf{x})(b(\mathbf{D})\Lambda(\mathbf{x})+\mathbf{1}_m).
\end{equation}

We need the following properties of the effective matrix, see \cite[Chapter~3, \S1]{BSu}.

\begin{proposition}
The effective matrix $g^0$ is subject to estimates 
$$\underline{g}\leqslant g^0\leqslant\overline{g}.$$ 
\textnormal{(}Here $\underline{g}$ and $\overline{g}$ are defined according to \eqref{underline and overline definition}.\textnormal{)} For $m=n$ one has $g^0=\underline{g}$.
\end{proposition}

\begin{proposition}
The equality $g^0=\overline{g}$ is equivalent to the relations 
\begin{equation}
\label{g0=overline g}
b(\mathbf{D})^*\mathbf{g}_k(\mathbf{x})=0,\quad k=1,\dots,m,
\end{equation}
for the columns $\mathbf{g}_k(\mathbf{x})$ of the matrix $g(\mathbf{x})$.
\end{proposition}

\begin{proposition}
The identity $g^0=\underline{g}$ is equivalent to the representations
\begin{equation*}
{\mathbf l}_k(\mathbf{x}) = {\mathbf l}_k^0 + b(\mathbf{D}) {\mathbf w}_k,\ \ {\mathbf l}_k^0\in \mathbb{C}^m,\ \
{\mathbf w}_k \in \widetilde{H}^1(\Omega;\mathbb{C}^m),\ \ k=1,\dots,m,
\end{equation*}
where ${\mathbf l}_k(\mathbf{x})$, $k=1,\dots,m,$ are the columns of the matrix $g(\mathbf{x})^{-1}$.
\end{proposition}

\subsection{The operator $\mathcal{B}^0(\mathbf{k},\varepsilon)$}
\label{Subsection B0(k,eps)}

The operator $\widehat{L}(t,\varepsilon)$ defined according to \eqref{L(t,eps)} and acting in the space $\widehat{\mathfrak{N}}$ now depends on $\boldsymbol{\theta}$. It turns out (see \cite[(7.2), (7.3), (7.8)]{Su10}) that
\begin{equation*}
\begin{split}
\widehat{L}(\mathbf{k},\varepsilon)
&=b(\mathbf{k})^*g^0b(\mathbf{k})\\
&-\varepsilon ( b(\mathbf{k})^*V+V^*b(\mathbf{k}))+\varepsilon\sum  _{j=1}^d(\overline{a_j+a_j^*})k_j 
+\varepsilon ^2(-W+\overline{\mathcal{Q}}+\lambda I),
\end{split}
\end{equation*}
where the constant matrices $(\overline{a_j+a_j^*})$ are defined in accordance with \eqref{underline and overline definition},
{\allowdisplaybreaks
\begin{align}
\label{V}
V &:=\vert \Omega \vert ^{-1}\int\limits_{\Omega}(b(\mathbf{D})\Lambda(\mathbf{x}))^*g(\mathbf{x})b(\mathbf{D})\widetilde{\Lambda}(\mathbf{x})\,d\mathbf{x},
\\
\label{W}
W &:=\vert \Omega \vert ^{-1}\int 
_{\Omega}(b(\mathbf{D})\widetilde{\Lambda}(\mathbf{x}))^*g(\mathbf{x})b(\mathbf{D})\widetilde{\Lambda}(\mathbf{x})\,d\mathbf{x},
\\
\label{overlineQ}
\overline{\mathcal{Q}}&:=\vert \Omega \vert ^{-1}\int\limits_\Omega d\mu(\mathbf{x}).
\end{align}
}

In consistence with \eqref{A(k)>}, 
$\widehat{\mathcal{A}}(\mathbf{k})\geqslant \widehat{c}_*\vert 
\mathbf{k}\vert ^2I$, $\mathbf{k}\in \widetilde{\Omega}$,
where $\widehat{c}_* = \alpha_0 \|g^{-1}\|^{-1}_{L_\infty}$.
Note that the constants $c_*$ and $\widehat{c}_*$
are related by the identity $c_*=\Vert f^{-1}\Vert^{-2}_{L_\infty} \widehat{c}_*$.
In accordance with \eqref{beta<=||M^-1||^-2hat beta}, $\beta \leqslant \Vert 
f^{-1}\Vert ^{-2}_{L_\infty}\widehat{\beta}$.
By \eqref{DOcheck c_*}, 
$\check{c}_*=\frac{1}{2}\min \lbrace \kappa c_*;2\beta \rbrace$,
 $\widehat{\check{c}}_* =\frac{1}{2}\min \lbrace \kappa \widehat{c}_*;2 \widehat{\beta} \rbrace$. So,
$\check{c}_*\leqslant \Vert f^{-1}\Vert 
^{-2}_{L_\infty}\widehat{\check{c}}_*$. According to \eqref{L(t,eps)>=}, 
$$\widehat{L}(\mathbf{k},\varepsilon ) \geqslant \break
\widehat{\check{c}}_*(\vert \mathbf{k}\vert ^2+\varepsilon 
^2)\mathbf{1}_n.$$
Together with \eqref{f_0<f} this implies the estimate
\begin{equation}
\label{f_0Lf_0>=}
f_0 \widehat{L}(\mathbf{k},\varepsilon )f_0\geqslant \check{c}_*(\vert \mathbf{k}\vert ^2+\varepsilon ^2)\mathbf{1}_n,\quad \mathbf{k}\in \mathbb{R}^d.
\end{equation}

Set 
\begin{align}
\widehat{\mathcal{A}}^0(\mathbf{k}) &=b(\mathbf{D}+\mathbf{k})^*g^0b(\mathbf{D}+\mathbf{k}),\quad \widehat{\mathcal{Y}}^0(\mathbf{k})=-b(\mathbf{D}+\mathbf{k})^*V+\sum\limits_{j=1}^d\overline{a_j}(D_j+k_j),
\nonumber
\\
\label{hat B0 (k,eps)}
\widehat{\mathcal{B}}^0(\mathbf{k},\varepsilon ) &=\widehat{\mathcal{A}}^0(\mathbf{k})+\varepsilon(\widehat{\mathcal{Y}}^0(\mathbf{k})+\widehat{\mathcal{Y}}^0(\textbf{k})^*)+\varepsilon ^2(\overline{\mathcal{Q}}-W+\lambda I).
\end{align}
Then 
\begin{equation}
\label{LP=B0P}
\widehat{L}(\mathbf{k},\varepsilon)\widehat{P}=\widehat{\mathcal{B}}^0(\mathbf{k},\varepsilon)\widehat{P}.
\end{equation}

Denote 
\begin{equation}
\label{B0 (k,eps)}
\mathcal{B}^0(\mathbf{k},\varepsilon):=f_0\widehat{\mathcal{B}}^0(\mathbf{k},\varepsilon)f_0,
\end{equation}
where $\widehat{\mathcal{B}}^0(\mathbf{k},\varepsilon)$ is the operator \eqref{hat B0 (k,eps)}. Since the symbol of the operator $\mathcal{B}^0(\mathbf{k},\varepsilon)$ is subject to estimate \eqref{f_0Lf_0>=}, using the Fourier series expansion, one can show that 
\begin{equation}
\label{B0(k,eps)>=}
\mathcal{B}^0(\mathbf{k},\varepsilon)
\geqslant \check{c}_*(\vert \mathbf{k}\vert ^2+\varepsilon ^2)I,\quad \mathbf{k}\in \widetilde{\Omega}.
\end{equation}

\subsection{The case when $f\neq \mathbf{1}_n$. The operators $\widehat{Z}_G(\boldsymbol{\theta})$ and $\widehat{\widetilde{Z}}_G$} 

Now we return to the analysis of the general case where $f\neq\mathbf{1}_n$. We give a realisation for the operators discussed in Subsection~\ref{Subsection ZG and tilde ZG abstract scheme}. Define a $\Gamma$-periodic $(n\times 
m)$-matrix-valued function $\Lambda _{G} (\mathbf{x})$ as a (weak) solution of the problem
\begin{equation}
\label{Lambda-G problem}
b(\mathbf{D})^*g(\mathbf{x})\left( b(\mathbf{D})\Lambda _{G} (\mathbf{x})+\mathbf{1}_m\right) =0,\quad \int  _\Omega G(\mathbf{x})\Lambda _{G} (\mathbf{x})\,d\mathbf{x}=0.
\end{equation}
Cf.~\cite[\S5]{BSu05}.
Obviously, $\Lambda _{G} (\mathbf{x})$ differs from the solution $\Lambda 
(\mathbf{x})$ of problem \eqref{Lambda} by the constant summand:
\begin{equation}
\label{Lambda_rho=Lambda+Lambda_rho^0}
\Lambda _{G} (\mathbf{x})=\Lambda (\mathbf{x})+\Lambda _{G} ^0,\quad \Lambda _G ^0=-\left(\overline{G}\right)^{-1}\left( \overline{G\Lambda}\right).
\end{equation}
In \cite[Subsection 7.3]{BSu05} it was obtained that
\begin{equation}
\label{Lambda_rho^0<=}
\vert \Lambda _{G} ^0\vert \leqslant \mathfrak{C}_{G} =m^{1/2}(2r_0)^{-1}\alpha _0^{-1/2}\Vert g\Vert ^{1/2}_{L_\infty} \Vert g^{-1}\Vert ^{1/2}_{L_\infty}\Vert f\Vert ^2_{L_\infty}\Vert f^{-1}\Vert ^2_{L_\infty} .
\end{equation}

According to \cite[\S5]{BSu05}, the role of the operator $\widehat{Z}_G$ from Subsection~\ref{Subsection ZG and tilde ZG abstract scheme} 
is played by the operator 
\begin{equation}
\label{hat ZG for DO's}
\widehat{Z}_{G}(\boldsymbol{\theta})=\Lambda _{G} 
b(\boldsymbol{\theta})\widehat{P}.
\end{equation} 
By using the relation $b(\mathbf{D})\widehat{P}=0$,  
we get $t\widehat{Z}_{G}(\boldsymbol{\theta})=\Lambda _{G} 
b(\mathbf{D}+\mathbf{k})\widehat{P}$, $\mathbf{k}\in \mathbb{R}^d$.

Next, similarly to \eqref{tildeLambda}, we consider a $\Gamma$-periodic $(n\times n)$-matrix-valued solution $\widetilde{\Lambda}_G(\mathbf{x})$ of the problem
\begin{equation}
\label{tilde-Lambda-G problem}
b(\mathbf{D})^*g(\mathbf{x})b(\mathbf{D})\widetilde{\Lambda }_G(\mathbf{x}) +\sum\limits _{j=1}^d D_j a_j(\mathbf{x})^*=0,\quad \int  _\Omega G (\mathbf{x})\widetilde{\Lambda}_G (\mathbf{x})\,d\mathbf{x}=0.
\end{equation}
Note that 
\begin{equation}
\label{tildeLambda_rho=tildeLambda+}
\widetilde{\Lambda}_G(\mathbf{x})=\widetilde{\Lambda}(\mathbf{x})+\widetilde{\Lambda}^0_G ,\quad \widetilde{\Lambda }^0_G =-\left(\overline{G}\right)^{-1}\left(\overline{G \widetilde{\Lambda}}\right).
\end{equation}
In \cite[(7.12)]{M_AA}, it was shown that
\begin{equation}
\label{tildeLambda_rho^0<=}
\vert \widetilde{\Lambda }_G^0 \vert \leqslant \widetilde{\mathfrak{C}}_G =(2r_0)^{-1}C_an^{1/2}\alpha _0^{-1}\Vert g^{-1}\Vert _{L_\infty}\Vert f\Vert ^2_{L_\infty} \Vert f ^{-1}\Vert ^2_{L_\infty}\vert \Omega \vert ^{-1/2}.
\end{equation}
Here $C_a^2=\sum _{j=1}^d\int _\Omega \vert a_j(\mathbf{x})\vert ^2\,d\mathbf{x}$. 
By the definition of the operator $\widehat{\widetilde{Z}}_G$ and the matrix-valued function $\widetilde{\Lambda }_G$,  it follows that
\begin{equation}
\label{hat tilde ZG for DO's}
\widehat{\widetilde{Z}}_G=\widetilde{\Lambda}_G \widehat{P}.
\end{equation}

\subsection{The operator $\widehat{N}_G(\mathbf{k},\varepsilon)$}
\label{Subsection N(k,eps)}

The operator $\widehat{N}_G(t,\varepsilon)$ defined according \eqref{NG(t,eps):=}--\eqref{N_G22 abstract} now depends on $\boldsymbol{\theta}$ and looks like this:
\begin{equation}
\label{N=... chapter 2 start}
\widehat{N}_G(t,\varepsilon;\boldsymbol{\theta})=t^3 \widehat{N}_{G,11}(\boldsymbol{\theta})+t^2\varepsilon \widehat{N}_{G,12}(\boldsymbol{\theta})+t\varepsilon ^2 \widehat{N}_{G,21}(\boldsymbol{\theta})+\varepsilon ^3 \widehat{N}_{G,22},
\end{equation}
where the operators $\widehat{N}_{G,jl}(\boldsymbol{\theta})$ are defined by  \eqref{NG11 abstract}--\eqref{N_G22 abstract} with the operators $\widehat{X}_1$, $\widehat{Z}_G$, $\widehat{R}$, and $\widehat{Y}_1$ depending on $\boldsymbol{\theta}$. (It is clear that $\widehat{N}_{G,22}$ does not depend on $\boldsymbol{\theta}$.)

By using equalities \eqref{X_1 hat for DO's}, \eqref{P= for DO's}, \eqref{hat R for DO's}, and \eqref{hat ZG for DO's}, one can show (cf. \cite[Subsection~5.3]{BSu05}) that the first summand in the right-hand side of \eqref{N=... chapter 2 start} has the form
\begin{equation*}
\widehat{N}_{G,11}(\mathbf{k}):=t^3 \widehat{N}_{G,11}(\boldsymbol{\theta})
= b(\mathbf{k})^* M_G(\mathbf{k})b(\mathbf{k})\widehat{P},
\end{equation*}
where 
$
M_G(\mathbf{k})=\overline{\Lambda _G ^* b(\mathbf{k})^*\widetilde{g}}
+\overline{\widetilde{g}^*b(\mathbf{k})\Lambda _G} $. 
Note that $M_G(\mathbf{k})$ is a Hermitian\break$(m\times m)$-matrix-valued function of $\mathbf{k}$ homogeneous of the degree one. Thus, $M_G(\mathbf{k})$ is a symbol of a self-adjoint first order DO $M_G(\mathbf{D})$ with constant coefficients.

Put $\widehat{N}_{G,12}(\mathbf{k}):=t^2 \widehat{N}_{G,12}(\boldsymbol{\theta})$. Using identities \eqref{hat X0 for DO's}--\eqref{hat Y2 for DO's}, \eqref{hat R for DO's}, \eqref{hat ZG for DO's}, and \eqref{hat tilde ZG for DO's}, by analogy with \cite[Subsection 5.2]{Su14}, one can show that
\begin{equation*}
\varepsilon \widehat{N}_{G,12}(\mathbf{k})=\varepsilon\bigl(
b(\mathbf{k})^*T_{G,0}b(\mathbf{k})+M_{G,1}(\mathbf{k})b(\mathbf{k})+b(\mathbf{k})^*M_{G,1}(\mathbf{k})^*
\bigr)\widehat{P},
\end{equation*}
where
\begin{align*}
M_{G,1}(\mathbf{k})&=\overline{\widetilde{\Lambda}_G^*b(\mathbf{k})^*\widetilde{g}}
+\overline{ (b(\mathbf{D})\widetilde{\Lambda}_G)^*g b(\mathbf{k})\Lambda_G}
+\sum _{j=1}^d \overline{(a_j+a_j^*)\Lambda _G}k_j,
\\
T_{G,0}&=2\sum _{j=1}^d \mathrm{Re}\,\overline{\Lambda _G ^* a_j D_j\Lambda _G}.
\end{align*}

Let $\widehat{N}_{G,21}(\mathbf{k}):=t\widehat{N}_{G,21}(\boldsymbol{\theta})$.  
By analogy with \cite[Subsection 5.2]{Su14}, one can show that 
\begin{equation*}
\begin{split}
\varepsilon ^2 \widehat{N}_{G,21}(\mathbf{k})
&=\varepsilon ^2 
\bigl(
M_{G,2}(\mathbf{k})+M_{G,2}(\mathbf{k})^*+T^*_Gb(\mathbf{k})+b(\mathbf{k})^*T_G\bigr)\widehat{P}
\\
&+2\varepsilon ^2 \sum _{j=1}^d \mathrm{Re}\,\overline{(a_j+a_j^*)\widetilde{\Lambda}_G}k_j\widehat{P},
\end{split}
\end{equation*}
where
\begin{align*}
M_{G,2}(\mathbf{k})&=\overline{(b(\mathbf{D})\widetilde{\Lambda}_G)^* g b(\mathbf{k})\widetilde{\Lambda}_G},
\\
T_G&=\sum _{j=1}^d \left(
\overline{\Lambda ^* _G a_j(D_j\widetilde{\Lambda}_G)}
+\overline{(D_j\Lambda _G)^* a_j^*\widetilde{\Lambda}_G}\right)
+\overline{\Lambda ^*_G \mathcal{Q}}+\lambda\overline{\Lambda ^*_G}.
\end{align*}
Here $\overline{{\Lambda}_G^* \mathcal{Q}}=\vert \Omega\vert ^{-1}\int _\Omega {\Lambda}_G(\mathbf{x})^*\,d\mu(\mathbf{x})$.

Finally, similarly to \cite[(5.30), (5.31)]{Su14}, 
$
\widehat{N}_{G,22}=(\widetilde{T}_G+\widetilde{T}_G^*)\widehat{P}$. 
Here
\begin{equation*}
\widetilde{T}_G=\sum _{j=1}^d \overline{\widetilde{\Lambda}_G^*a_j(D_j\widetilde{\Lambda}_G)}
+\overline{\widetilde{\Lambda}_G^*\mathcal{Q}}+\lambda \overline{\widetilde{\Lambda}^*_G},
\end{equation*}
where $\overline{\widetilde{\Lambda}_G^* \mathcal{Q}}=\vert \Omega\vert ^{-1}\int _\Omega\widetilde{\Lambda}_G(\mathbf{x})^*\,d\mu(\mathbf{x})$.

As a result, the operator $\widehat{N}_G(\mathbf{k},\varepsilon)=\widehat{N}_G(t,\varepsilon;\boldsymbol{\theta})$ defined in \eqref{N=... chapter 2 start} is represented as 
\begin{equation*}
\widehat{N}_G(\mathbf{k},\varepsilon)=\widehat{N}_{G,11}(\mathbf{k})+\varepsilon \widehat{N}_{G,12}(\mathbf{k})+\varepsilon ^2 \widehat{N}_{G,21}(\mathbf{k})+\varepsilon ^3 \widehat{N}_{G,22}.
\end{equation*}
According to Remark~\ref{Remark Constants Chapter 2}, for the operator $\widehat{N}_G(\mathbf{k},\varepsilon)$, estimate of the form \eqref{NG<= abstract scheme} holds true and the constant $C_G$ can be chosen independent on~$\boldsymbol{\theta}$. So, for $\mathbf{k}\in\mathbb{R}^d$ and $0<\varepsilon \leqslant 1$ we have
\begin{equation}
\label{N(k,eps)<=}
\Vert \widehat{N}_G(\mathbf{k},\varepsilon)\Vert _{\mathfrak{H}\rightarrow\mathfrak{H}} 
=\vert \widehat{N}_G(\mathbf{k},\varepsilon)\vert
\leqslant 
C_G(\vert \mathbf{k}\vert ^2 +\varepsilon ^2)^{3/2}.
\end{equation}
(We taken into account that, for the operator of multiplication by a constant matrix, the operator norm in $L_2 (\Omega;\mathbb{C}^n)$ matches with the matrix norm.)

Since the projection $\widehat{P}$ stands in the expressions defining the operators $\widehat{N}_{G,11}(\mathbf{k})$, $\widehat{N}_{G,12}(\mathbf{k})$, and $\widehat{N}_{G,21}(\mathbf{k})$, we can replace $\mathbf{k}$ by $\mathbf{D}+\mathbf{k}$ in these expressions: 
$\widehat{N}_G(\mathbf{k},\varepsilon)=\mathcal{N}(\mathbf{k},\varepsilon)\widehat{P}$, 
where $\mathcal{N}(\mathbf{k},\varepsilon)$ is a self-adjoint third order DO:
\begin{equation*}
\mathcal{N}(\mathbf{k},\varepsilon)=\mathcal{N}_{11}(\mathbf{D}+\mathbf{k})+\varepsilon\mathcal{N}_{12}(\mathbf{D}+\mathbf{k})+\varepsilon ^2 \mathcal{N}_{21}(\mathbf{D}+\mathbf{k})+\varepsilon ^3 \mathcal{N}_{22}.
\end{equation*}
Here the summands in the right-hand side are DO's of the third, second, first, and zero order, respectively, given by the relations
\begin{align}
\mathcal{N}_{11}(\mathbf{D}+\mathbf{k})&=b(\mathbf{D}+\mathbf{k})^*M_G(\mathbf{D}+\mathbf{k})b(\mathbf{D}+\mathbf{k}),
\nonumber
\\
\begin{split}
\mathcal{N}_{12}(\mathbf{D}+\mathbf{k})&=b(\mathbf{D}+\mathbf{k})^*T_{G,0}b(\mathbf{D}+\mathbf{k})+M_{G,1}(\mathbf{D}+\mathbf{k})b(\mathbf{D}+\mathbf{k})
\\
&+b(\mathbf{D}+\mathbf{k})^*M_{G,1}(\mathbf{D}+\mathbf{k})^*,
\end{split}
\nonumber
\\
\begin{split}
\mathcal{N}_{21}(\mathbf{D}+\mathbf{k})&=
M_{G,2}(\mathbf{D}+\mathbf{k})+M_{G,2}(\mathbf{D}+\mathbf{k})^*+{T}^*_Gb(\mathbf{D}+\mathbf{k})
+b(\mathbf{D}+\mathbf{k})^* {T}_G
\\
&+2\sum _{j=1}^d \mathrm{Re}\,\overline{(a_j+a_j^*)\widetilde{\Lambda}_G}(D_j+k_j),
\end{split}
\nonumber
\\
\label{N22 for DO's}
\mathcal{N}_{22}&=\widetilde{T}_G+\widetilde{T}^*_G.
\end{align}

\begin{remark}
\label{Remark N(k,eps)=0}
If the $\Gamma$-periodic solutions of problems \eqref{Lambda-G problem} and \eqref{tilde-Lambda-G problem} are equal to zero\textnormal{:} $\Lambda _G =0$, $\widetilde{\Lambda}_G=0$, then, by construction, $\mathcal{N}(\mathbf{k},\varepsilon)=0$.
\end{remark}

\section{Approximation for the operator $fe^{-\mathcal{B}(\mathbf{k},\varepsilon)s}f^*$}

\subsection{}

Now we apply Theorem~\ref{Theorem bordered exp final in abstrach scheme} to the operator $\mathcal{B}(\mathbf{k},\varepsilon)$. 
Taking equalities \eqref{Z= for DO's}, \eqref{tilde Z = for DO's}, \eqref{LP=B0P}, and \eqref{B0 (k,eps)} into account, from Theorem~\ref{Theorem bordered exp final in abstrach scheme} with $\vert \tau\vert =(\vert \mathbf{k}\vert ^2+\varepsilon ^2)^{1/2}\leqslant\tau _0$ 
we deduce that
\begin{equation}
\label{Ch 2 Th for small k}
\begin{split}
fe^{-\mathcal{B}(\mathbf{k},\varepsilon)s}f^*=f_0e^{-\mathcal{B}^0(\mathbf{k},\varepsilon)s}f_0\widehat{P}+\mathcal{K}(\mathbf{k},\varepsilon,s)+f\mathcal{R}(\mathbf{k},\varepsilon,s)f^*.
\end{split}
\end{equation}
Here  
the operator $\mathcal{B}^0(\mathbf{k},\varepsilon)$ was defined in \eqref{B0 (k,eps)},
\begin{equation}
\label{K(k,eps)}
\begin{split}
\mathcal{K}(\mathbf{k},\varepsilon,s)
&:=
\left( \Lambda _G b(\mathbf{D}+\mathbf{k})+\varepsilon\widetilde{\Lambda}_G\right)f_0e^{-\mathcal{B}^0(\mathbf{k},\varepsilon)s}f_0\widehat{P}
\\
&+f_0e^{-\mathcal{B}^0(\mathbf{k},\varepsilon)s}f_0\widehat{P} \left( b(\mathbf{D}+\mathbf{k})^*\Lambda _G ^*+\varepsilon\widetilde{\Lambda}_G^*\right)
\\
&-\int\limits_0^s f_0e^{-\mathcal{B}^0(\mathbf{k},\varepsilon)(s-\widetilde{s})}f_0\mathcal{N}(\mathbf{k},\varepsilon)f_0e^{-\mathcal{B}^0(\mathbf{k},\varepsilon)\widetilde{s}}f_0\widehat{P}\,d\widetilde{s},
\end{split}
\end{equation}
and the remainder term $f\mathcal{R} (\mathbf{k},\varepsilon,s)f^*$ for $(\vert \mathbf{k}\vert ^2+\varepsilon ^2)^{1/2}\leqslant\tau _0$ satisfies
\begin{align}
\label{fR2f small k s>=0}
&\Vert f\mathcal{R} (\mathbf{k},\varepsilon,s)f^*\Vert _{L_2(\Omega)\rightarrow L_2(\Omega)} \!\leqslant\! C_7 \Vert f\Vert ^2_{L_\infty}\!(s\!+\!1)^{\!-1}e^{-\check{c}_*(\vert \mathbf{k}\vert ^2\! +\varepsilon ^2)s/2}\!,\quad s\!\geqslant 0;
\\
\label{fR2f small k s>0}
&\Vert f\mathcal{R} (\mathbf{k},\varepsilon,s)f^*\Vert _{L_2(\Omega)\rightarrow L_2(\Omega)} \!\leqslant\! C_5\Vert f\Vert ^2_{L_\infty}s^{-1}e^{-\check{c}_*(\vert \mathbf{k}\vert ^2 +\varepsilon ^2)s/2},\quad s> 0.
\end{align}

\subsection{Estimates for $\vert \mathbf{k}\vert ^2 +\varepsilon ^2>\tau _0^2$}

If $\mathbf{k}\in\widetilde{\Omega}$, $0<\varepsilon\leqslant 1$, and $\vert \mathbf{k}\vert ^2 +\varepsilon ^2 >\tau _0^2$, it suffices to estimate each operator in  \eqref{Ch 2 Th for small k} separately. By \eqref{B(t,eps,theta)>=} and the elementary inequalities $e^{-\alpha}\leqslant \alpha ^{-1}e^{-\alpha /2}$ for $\alpha >0$ and $e^{-\alpha}\leqslant 2 (1+\alpha)^{-1}e^{-\alpha /2} $ for $\alpha\geqslant 0$, we see that for $\mathbf{k}\in\widetilde{\Omega}$,  $\vert \mathbf{k}\vert ^2 +\varepsilon ^2>\tau _0^2$ one has 
\begin{align}
\label{exp B(k,eps) big k, s>0}
\begin{split}
\Vert fe^{-\mathcal{B}(\mathbf{k},\varepsilon)s}f^*\Vert  _{L_2(\Omega)\rightarrow L_2(\Omega)} &\leqslant \Vert f\Vert ^2 _{L_\infty}e^{-\check{c}_*(\vert \mathbf{k}\vert ^2+\varepsilon ^2) s}\\
&\leqslant 
\tau _0^{-2}\check{c}_* ^{-1} s^{-1}e^{-\check{c}_*(\vert \mathbf{k}\vert ^2+\varepsilon ^2) s/2}\Vert f\Vert ^2 _{L_\infty},\quad s>0;
\end{split}
\\
\label{exp B(k,eps) big k, s>=0}
\begin{split}
\Vert fe^{-\mathcal{B}(\mathbf{k},\varepsilon)s}f^*\Vert  _{L_2(\Omega)\rightarrow L_2(\Omega)} &\leqslant 
2\max\lbrace 1;\check{c}_*^{-1}\tau _0^{-2}\rbrace 
\\
&\times (1+s)^{-1}e^{-\check{c}_*(\vert \mathbf{k}\vert ^2+\varepsilon ^2) s/2}
\Vert f\Vert ^2 _{L_\infty},\quad s\geqslant 0.
\end{split}
\end{align}
Similarly, by \eqref{f_0<f} and \eqref{B0(k,eps)>=} for $\mathbf{k}\in\widetilde{\Omega}$,  $\vert \mathbf{k}\vert ^2 +\varepsilon ^2>\tau _0^2$ we have
\begin{align}
\label{exp B0(k,eps) big k, s>0}
\Vert f_0e^{-\mathcal{B}^0(\mathbf{k},\varepsilon)s}f_0\widehat{P}\Vert  _{L_2(\Omega)\rightarrow L_2(\Omega)} & \!\leqslant \!\tau _0^{\!-2}\check{c}_*^{-1}s^{\!-1}e^{\!-\check{c}_*(\vert \mathbf{k}\vert ^2+\varepsilon ^2) s/2}\Vert f\Vert ^2 _{L_\infty},\!\quad s\!>\!0;
\\
\label{exp B0(k,eps) big k, s>=0}
\begin{split}
\Vert f_0e^{-\mathcal{B}^0(\mathbf{k},\varepsilon)s}f_0\widehat{P}\Vert  _{L_2(\Omega)\rightarrow L_2(\Omega)} &\!\leqslant\! 2\max\lbrace 1;\check{c}_*^{-1}\tau _0^{\!-2}\rbrace 
\\
&\!\!\times (1\!+\!s)^{\!-1}e^{\!-\check{c}_*(\vert \mathbf{k}\vert ^2\!+\varepsilon ^2) s/2}
\Vert f\Vert ^2 _{L_\infty},\!\quad s\!\geqslant\! 0.
\end{split}
\end{align}
By \eqref{KG (t,eps)<= for s>0}, \eqref{KG (t,eps)<= for s>=0} and Remark~\ref{Remark Constants Chapter 2}, for the corrector we have estimates with constants independent of $\boldsymbol{\theta}$:
{\allowdisplaybreaks\begin{align}
\label{K (k,eps)<= for s>0}
\begin{split}
\Vert & \mathcal{K}(\mathbf{k},\varepsilon,s)\Vert   _{L_2(\Omega)\rightarrow L_2(\Omega)}
\leqslant C_8\Vert f\Vert _{L_\infty} ^2 s^{-1}(\vert \mathbf{k}\vert ^2+\varepsilon ^2)^{-1/2}e^{-\check{c}_*(\vert \mathbf{k}\vert ^2+\varepsilon ^2)s/2}
\\
&\leqslant  C_8\tau _0^{-1}\Vert f\Vert  _{L_\infty} ^2 s^{-1}e^{-\check{c}_*( \vert \mathbf{k}\vert ^2+\varepsilon ^2)s/2}
,\\ 
&s>0,\quad 0<\varepsilon\leqslant 1,\quad \mathbf{k}\in\widetilde{\Omega},\quad \vert \mathbf{k}\vert ^2 +\varepsilon ^2 >\tau _0^2;
\end{split}
\\
\label{K (k,eps)<= for s>=0}
\begin{split}
\Vert  &\mathcal{K}(\mathbf{k},\varepsilon,s)\Vert   _{L_2(\Omega)\rightarrow L_2(\Omega)}
\\
&\leqslant C_9\Vert f\Vert  _{L_\infty} ^2 (\vert \mathbf{k}\vert ^2+\varepsilon ^2)^{1/2}\left(1+\check{c}_*( \vert \mathbf{k}\vert ^2+\varepsilon ^2)s\right)^{-1}e^{-\check{c}_*(\vert \mathbf{k}\vert ^2+\varepsilon ^2)s/2}
\\
&\leqslant
C_9\Vert f\Vert  _{L_\infty} ^2 (r_1^2+1)^{1/2}\max\lbrace 1;\check{c}_*^{-1}\tau _0^{-2}\rbrace
(1+s)^{-1}e^{-\check{c}_*(\vert \mathbf{k}\vert ^2+\varepsilon ^2)s/2},
\\
&s\geqslant 0,\quad 0<\varepsilon\leqslant 1,\quad \mathbf{k}\in\widetilde{\Omega},\quad \vert \mathbf{k}\vert ^2 +\varepsilon ^2 >\tau _0^2.
\end{split}
\end{align}
}

Combining \eqref{fR2f small k s>0}, \eqref{exp B(k,eps) big k, s>0}, \eqref{exp B0(k,eps) big k, s>0}, and \eqref{K (k,eps)<= for s>0}, we conclude that for $s>0$ and any $\mathbf{k}\in\widetilde{\Omega}$ the representation \eqref{Ch 2 Th for small k} holds true with the error estimate 
\begin{align}
\label{fR2f all k s>0}
\Vert f\mathcal{R} (\mathbf{k},\varepsilon,s)f^*\Vert  _{L_2(\Omega)\rightarrow L_2(\Omega)}\!\leqslant\! C_{10} s^{-1}e^{-\check{c}_*(\vert \mathbf{k}\vert ^2 +\varepsilon ^2)s/2},\quad s\!>\! 0, \quad \mathbf{k}\!\in\!\widetilde{\Omega}.
\end{align}
Here $C_{10} =\max \lbrace  C_5 ; 2\tau _0^{-2}\check{c}_*^{-1}+C_8\tau _0^{-1}\rbrace \Vert f\Vert ^2_{L_\infty}$.

By \eqref{fR2f small k s>=0}, \eqref{exp B(k,eps) big k, s>=0}, \eqref{exp B0(k,eps) big k, s>=0}, and \eqref{K (k,eps)<= for s>=0}, for the remainder term in \eqref{Ch 2 Th for small k} for $s\geqslant 0$ and $\mathbf{k}\in\widetilde{\Omega}$ we have the estimate
\begin{equation}
\label{fR2f all k s>=0}
\Vert f\mathcal{R} (\mathbf{k},\varepsilon,s)f^*\Vert _{L_2(\Omega)\rightarrow L_2(\Omega)} \!\!\leqslant\! C_{11} (1\!+\!s)^{\!-1}e^{\!-\check{c}_*(\vert \mathbf{k}\vert ^2 \!+\varepsilon ^2)s/2},\!\quad s\!\geqslant\! 0, \!\quad \mathbf{k}\!\in\!\widetilde{\Omega},
\end{equation}
with the constant $C_{11}=\max\left\lbrace C_7 ;(4+C_9(1+r_1^2)^{1/2})\max\lbrace 1;\check{c}_*^{-1}\tau _0^{-2}\rbrace \right\rbrace \Vert f\Vert ^2_{L_\infty}$.

We have got the following result.

\begin{theorem}
\label{Theorem exp of B(k,eps)}
Let $\mathcal{B}(\mathbf{k},\varepsilon)$ and $\mathcal{B}^0(\mathbf{k},\varepsilon)$ be the operators \eqref{B(k,e)} and \eqref{B0 (k,eps)}, respectively, and let $\mathcal{K}(\mathbf{k},\varepsilon,s)$ be the corrector \eqref{K(k,eps)}. Then for $s\geqslant 0$, $0<\varepsilon\leqslant 1$, and $\mathbf{k}\in\widetilde{\Omega}$ we have the representation \eqref{Ch 2 Th for small k},  and the remainder term satisfies estimates \eqref{fR2f all k s>0} and \eqref{fR2f all k s>=0}, where the constants  $C_{10}$ and $C_{11}$ are controlled in terms of the initial data \eqref{dannie_zadachi}.
\end{theorem}

\section{Approximation for the operator $fe^{-\mathcal{B}(\varepsilon)s}f^*$}

\label{Section f exp B(eps) f*}

\subsection{}

\label{Subsection exp into direct integral}

Now we return to the analysis of the operator $\mathcal{B}(\varepsilon)$ acting in $L_2(\mathbb{R}^d;\mathbb{C}^n)$ and defined in Subsection~\ref{Subsection B(eps)}. Approximation for the bordered operator exponential $fe^{-\mathcal{B}(\varepsilon)s}f^*$ can be derived from Theorem~\ref{Theorem exp of B(k,eps)} by using the direct integral decomposition.

We introduce the effective operator with constant coefficients:
\begin{equation}
\label{B0(eps)}
\begin{split}
\mathcal{B}^0(\varepsilon)&:=
f_0b(\mathbf{D})^*g^0b(\mathbf{D})f_0
\\
&+\varepsilon f_0\bigg(-b(\mathbf{D})^*V-V^*b(\mathbf{D})+\sum_{j=1}^d \overline{(a_j+a_j^*)}D_j\bigg)f_0
\\
&
+\varepsilon ^2 f_0(-W+\overline{\mathcal{Q}}+\lambda I)f_0.
\end{split}
\end{equation}
The symbol of the operator \eqref{B0(eps)} is the matrix $f_0\widehat{L}(\mathbf{k},\varepsilon)f_0$ (see Subsection~\ref{Subsection B0(k,eps)}).

Define the third order DO with constant coefficients:
\begin{align}
\label{9.1a}
\mathcal{N}(\varepsilon)&=\mathcal{N}_{11}(\mathbf{D})+\varepsilon\mathcal{N}_{12}(\mathbf{D})+\varepsilon ^2 \mathcal{N}_{21}(\mathbf{D})+\varepsilon ^3\mathcal{N}_{22},
\\
\mathcal{N}_{11}(\mathbf{D})&=b(\mathbf{D})^*M_G(\mathbf{D})b(\mathbf{D}),
\nonumber
\\
\mathcal{N}_{12}(\mathbf{D})&=b(\mathbf{D})^*T_{G,0}b(\mathbf{D})+M_{G,1}(\mathbf{D})b(\mathbf{D})+b(\mathbf{D})^*M_{G,1}(\mathbf{D})^*,
\nonumber
\\
\begin{split}
\mathcal{N}_{21}(\mathbf{D})&=M_{G,2}(\mathbf{D})+M_{G,2}(\mathbf{D})^*+T_G^*b(\mathbf{D})+b(\mathbf{D})^*T_G
\\
&+2\sum _{j=1}^d \mathrm{Re}\,\overline{(a_j+a_j^*)\widetilde{\Lambda}_G}D_j.
\end{split}
\nonumber
\end{align}
(Recall that the matrix $\mathcal{N}_{22}$ was defined in \eqref{N22 for DO's}.)

We use the direct integral decomposition for the operator $\mathcal{B}(\varepsilon)$, see Subsection~\ref{Subsection direct integral decomposition}. Then for the operator exponential one has the representation
\begin{equation*}
e^{-\mathcal{B}(\varepsilon)s}=\mathcal{U}^{-1}\bigg(\int\limits_{\widetilde{\Omega}}\oplus e^{-\mathcal{B}(\mathbf{k},\varepsilon)s}\,d\mathbf{k}\bigg)\mathcal{U}.
\end{equation*}
A similar identity holds true for $e^{-\mathcal{B}^0(\varepsilon)s}$. The operator 
$$
\Lambda _G b(\mathbf{D})+\varepsilon\widetilde{\Lambda}_G
$$
decomposes into the direct integral with fibers 
$\Lambda _G b(\mathbf{D}+\mathbf{k})+\varepsilon\widetilde{\Lambda}_G$. For the operator $\mathcal{N}(\varepsilon)$, we have decomposition into the direct integral of the operators  $\mathcal{N}(\mathbf{k},\varepsilon)$. Thus, Remark~\ref{Remark N(k,eps)=0} implies the following observation.

\begin{remark}
\label{Remark N(eps)=0}
If $\Lambda _G =0$ and $\widetilde{\Lambda}_G=0$, then  $\mathcal{N}(\varepsilon)=0$.
\end{remark}

Define a bounded operator $\Pi$ in $L_2(\mathbb{R}^d;\mathbb{C}^n)$ by the relation  
$$
\Pi 
=\mathcal{U}^{-1}[\widehat{P}]\mathcal{U},
$$
where $[\widehat{P}]$ is the operator acting in the space $\mathcal{H}=\int _{\widetilde{\Omega}}\oplus L_2 (\Omega ;\mathbb{C}^n)\,d\mathbf{k}$ fiber-wise as the operator $\widehat{P}$ of averaging over the cell. As was obtained in \cite[Subsection~6.1]{BSu05}, the operator $\Pi$ can be written as 
\begin{equation}
\label{Pi}
(\Pi \mathbf{u})(\mathbf{x})=(2\pi )^{-d/2}\int\limits_{\widetilde{\Omega}} e^{i\langle \mathbf{x},\boldsymbol{\xi}\rangle}(\mathcal{F}\mathbf{u})(\boldsymbol{\xi})\,d\boldsymbol{\xi},
\end{equation}
where $(\mathcal{F}\mathbf{u})(\,\cdot\,)$ is the Fourier image of the function $\mathbf{u}$. Thus, $\Pi$ is PDO with the symbol $\chi _{\widetilde{\Omega}}(\boldsymbol{\xi})$, where $\chi _{\widetilde{\Omega}}$ is the characteristic function of the set $\widetilde{\Omega}$. This operator is a smoothing one. Note that $\Pi$ commutes with differential operators with constant coefficients.

Theorem~\ref{Theorem exp of B(k,eps)} implies the identity
\begin{equation}
\label{exp B(e) approx with Pi}
\begin{split}
fe^{-\mathcal{B}(\varepsilon)s}f^*
=
f_0e^{-\mathcal{B}^0(\varepsilon)s}f_0\Pi +\mathcal{K}(\varepsilon,s)+\mathcal{R}(\varepsilon,s).
\end{split}
\end{equation}
Here the corrector $\mathcal{K}(\varepsilon,s)=\mathcal{U}^{-1}\left(\int _{\widetilde{\Omega}}\oplus \mathcal{K}(\mathbf{k},\varepsilon,s) \,d\mathbf{k}\right)\mathcal{U}$ has the form
\begin{equation}
\label{K(eps,s) with Pi}
\begin{split}
\mathcal{K}(\varepsilon,s)
&=
\left( \Lambda _G b(\mathbf{D})+\varepsilon\widetilde{\Lambda}_G\right)f_0e^{-\mathcal{B}^0(\varepsilon)s}f_0\Pi
\\
&+f_0e^{-\mathcal{B}^0(\varepsilon)s}f_0\Pi \left( b(\mathbf{D})^*\Lambda _G ^*+\varepsilon\widetilde{\Lambda}_G^*\right)
\\
&-\int\limits_0^s f_0e^{-\mathcal{B}^0(\varepsilon)(s-\widetilde{s})}f_0\mathcal{N}(\varepsilon)f_0e^{-\mathcal{B}^0(\varepsilon)\widetilde{s}}f_0\Pi\,d\widetilde{s} ,
\end{split}
\end{equation}
and the operator  $\mathcal{R} (\varepsilon,s)=  \mathcal{U}^{-1}\left(\int _{\widetilde{\Omega}}\oplus f \mathcal{R} (\mathbf{k},\varepsilon,s)f^*\,d\mathbf{k}\right)\mathcal{U}$ is subject to the estimates 
\begin{align*}
\begin{split}
\Vert \mathcal{R} (\varepsilon,s)\Vert _{L_2(\mathbb{R}^d)\rightarrow L_2(\mathbb{R}^d)}
&=
\esssup_{\mathbf{k}\in\widetilde{\Omega}}\Vert f \mathcal{R} (\mathbf{k},\varepsilon,s)f^*\Vert _{L_2(\Omega)\rightarrow L_2(\Omega)}
\\
&\leqslant C_{10} s^{-1}e^{-\check{c}_*\varepsilon ^2 s/2},\quad s> 0; 
\end{split}
\\
\Vert \mathcal{R} (\varepsilon,s)\Vert _{L_2(\mathbb{R}^d)\rightarrow L_2(\mathbb{R}^d)}&\leqslant C_{11} (1+s)^{-1}e^{-\check{c}_*\varepsilon ^2 s/2},\quad s\geqslant 0.
\end{align*}

\subsection{Elimination of the operator $\Pi$} Now we discuss the possibility to remove the smoothing operator from the corrector in the approximation \eqref{exp B(e) approx with Pi}. Put 
\begin{equation}
\label{Xi (eps,s)}
\Xi (\varepsilon ,s):=f_0e^{-\mathcal{B}^0(\varepsilon)s}f_0(I-\Pi).
\end{equation}
Since the matrix $f_0\widehat{L}(\boldsymbol{\xi},\varepsilon)f_0$ is the symbol of the operator $\mathcal{B}^0(\varepsilon)$ and $\Pi$ is a PDO with the symbol $\chi _{\widetilde{\Omega}}(\boldsymbol{\xi})$, by \eqref{f_0<f}, \eqref{f_0Lf_0>=}, and the elementary inequality $e^{-\alpha}\leqslant 2(1+\alpha )^{-1}e^{-\alpha /2}$, $\alpha \geqslant 0$, 
for $s\geqslant 0$ we have 
\begin{equation*}
\begin{split}
\Vert \Xi (\varepsilon ,s)\Vert _{L_2(\mathbb{R}^d)\rightarrow L_2(\mathbb{R}^d)}
&\leqslant \Vert f\Vert ^2 _{L_\infty}\sup _{\boldsymbol{\xi}\in\mathbb{R}^d} \vert e^{-f_0\widehat{L}(\boldsymbol{\xi},\varepsilon)f_0s}\vert (1-\chi _{\widetilde{\Omega}}(\boldsymbol{\xi}))
\\
&\leqslant \Vert f\Vert ^2 _{L_\infty}\sup _{\boldsymbol{\xi}\in\mathbb{R}^d,\,\vert \boldsymbol{\xi}\vert \geqslant r_0} e^{-\check{c}_*(\vert \boldsymbol{\xi}\vert ^2 +\varepsilon ^2)s}
\\
&\leqslant
2\Vert f\Vert ^2 _{L_\infty} \max\lbrace 1;\check{c}_*^{-1}r_0^{-2}\rbrace(1+s)^{-1}e^{-\check{c}_*\varepsilon ^2 s/2}.
\end{split}
\end{equation*}
Thus, the smoothing operator $\Pi$ can always be removed from the principal term of approximation \eqref{exp B(e) approx with Pi}.

For $s>0$, we also can replace $\Pi$ by $I$ in the third term of the corrector~\eqref{K(eps,s) with Pi}. Note that $\mathcal{N}(\varepsilon)$ is DO with the symbol $\widehat{N}_G(\boldsymbol{\xi},\varepsilon)$ (see Subsection~\ref{Subsection N(k,eps)}) satisfying  estimate \eqref{N(k,eps)<=}. From this, \eqref{f_0<f}, and estimate~\eqref{f_0Lf_0>=} for the symbol of the operator $\mathcal{B}^0(\varepsilon)$ we derive that
\begin{equation*}
\begin{split}
\Bigl\Vert & \int\limits_0^s  f_0e^{-\mathcal{B}^0(\varepsilon)(s-\widetilde{s})}f_0\mathcal{N}(\varepsilon)f_0e^{-\mathcal{B}^0(\varepsilon)\widetilde{s}}f_0(I-\Pi)\,d\widetilde{s}
\Bigr\Vert _{L_2(\mathbb{R}^d)\rightarrow L_2(\mathbb{R}^d)}
\\
&\leqslant
\Vert f\Vert ^4 _{L_\infty}\sup _{\boldsymbol{\xi}\in\mathbb{R}^d}\int\limits_0^s\vert e^{-f_0\widehat{L}(\boldsymbol{\xi},\varepsilon)f_0(s-\widetilde{s})}\vert \vert \widehat{N}_G(\boldsymbol{\xi},\varepsilon)\vert \vert e^{-f_0\widehat{L}(\boldsymbol{\xi},\varepsilon)f_0\widetilde{s}}\vert\,d\widetilde{s}(1-\chi _{\widetilde{\Omega}}(\boldsymbol{\xi}))
\\
&\leqslant
\Vert f\Vert ^4_{L_\infty}\sup _{\boldsymbol{\xi}\in\mathbb{R}^d,\vert \boldsymbol{\xi}\vert \geqslant r_0}s e^{-\check{c}_*(\vert \boldsymbol{\xi}\vert ^2+\varepsilon ^2)s}C_G(\vert \boldsymbol{\xi}\vert ^2+\varepsilon ^2)^{3/2}
\\
&\leqslant 3 \check{c}_*^{-2}r_0^{-1}C_G\Vert f\Vert ^4 _{L_\infty}s^{-1}e^{-\check{c}_*\varepsilon ^2 s/2},\quad s>0.
\end{split}
\end{equation*}
(We have used the elementary inequality $e^{-\alpha}\leqslant 3\alpha ^{-2}e^{-\alpha /2}$, $\alpha >0$.)

Now, we discuss the possibility of elimination of the operator $\Pi$ from the other terms of the corrector. 
Since the matrix-valued functions $\Lambda _G$ and $\Lambda$, $\widetilde{\Lambda}_G$ and $\widetilde{\Lambda}$ differ by a constant summand  (see \eqref{Lambda_rho=Lambda+Lambda_rho^0} and \eqref{tildeLambda_rho=tildeLambda+}), we use the multiplier properties of the matrix-valued functions  $\Lambda$ and $\widetilde{\Lambda}$ to eliminate the smoothing operator from the terms of the corrector containing  $\Lambda _G$ and $\widetilde{\Lambda} _G$. 
The following result was obtained in \cite[Proposition~6.8]{BSu05} for $d\leqslant 4$ and \cite[(7.19) and (7.20)]{V} for $d>4$.

\begin{lemma}
\label{Lemma Lambda multiplier}
Let $\Lambda$ be a $\Gamma$-periodic $(n\times m)$-matrix-valued solution of the problem \eqref{Lambda}. Let $l=1$ for $d\leqslant 4$ and $l=d/2-1$ for $d>4$. Then the operator  $[\Lambda]$ of multiplication by the matrix-valued function $\Lambda$  maps $H^l(\mathbb{R}^d;\mathbb{C}^m)$ into $L_2(\mathbb{R}^d;\mathbb{C}^n)$ continuously, and
\begin{equation*}
\Vert [\Lambda ]\Vert _{H^l(\mathbb{R}^d)\rightarrow L_2(\mathbb{R}^d)}\leqslant \mathfrak{C}_\Lambda .
\end{equation*}
The constant $ \mathfrak{C}_\Lambda$ is controlled in terms of $m$, $n$, $d$, $\alpha _0$, $\Vert g\Vert _{L_\infty}$, $\Vert g^{-1}\Vert _{L_\infty}$, and the parameters of the lattice $\Gamma$.
\end{lemma}

For $d\leqslant 6$, we apply Proposition~6.9 from \cite{Su14}, and for 
$d> 6$ we rely on Lemma~6.5($1^\circ$) from \cite{MSu}.

\begin{lemma}
\label{Lemma tilde Lambda multiplier}
Let $\widetilde{\Lambda}$ be a $\Gamma$-periodic $(n\times n)$-matrix-valued solution of the problem \eqref{tildeLambda}. Put $\sigma =2$ for $d\leqslant 6$ and $\sigma =d/2-1$ for $d>6$. Then the operator $[\widetilde{\Lambda}]$ of multiplication by the matrix-valued function $\widetilde{\Lambda}$ maps $H^\sigma (\mathbb{R}^d;\mathbb{C}^n)$ into $L_2(\mathbb{R}^d;\mathbb{C}^n)$ continuously and
\begin{equation*}
\Vert [\widetilde{\Lambda} ]\Vert _{H^\sigma(\mathbb{R}^d)\rightarrow L_2(\mathbb{R}^d)}\leqslant \mathfrak{C}_{\widetilde{\Lambda}} .
\end{equation*}
The constant $\mathfrak{C}_{\widetilde{\Lambda}}$ is controlled in terms of the initial data \eqref{dannie_zadachi}.
\end{lemma} 

The following result was obtained in \cite[Proposition 8.3]{M_AA}.

\begin{proposition}
Let $\Xi (\varepsilon ,s)$ be the operator \eqref{Xi (eps,s)}. Then for $s>0$ and all $l>0$ the operators $b(\mathbf{D})\Xi (\varepsilon ,s)$ and $\varepsilon \Xi (\varepsilon ,s)$ are bounded from $L_2(\mathbb{R}^d;\mathbb{C}^n)$ into $H^l (\mathbb{R}^d;\mathbb{C}^n)$ and
\begin{align}
\label{b(D)Xi(eps,s)<=}
\Vert b(\mathbf{D})\Xi (\varepsilon ,s)\Vert _{L_2(\mathbb{R}^d)\rightarrow H^l(\mathbb{R}^d)}
&\leqslant \alpha _1^{1/2}\mathscr{C}(l)s^{-(l+1)/2}e^{-\check{c}_*\varepsilon ^2 s/2},
\\
\label{eps Xi(eps,s)<=}
\varepsilon\Vert \Xi (\varepsilon ,s)\Vert _{L_2(\mathbb{R}^d)\rightarrow H^l(\mathbb{R}^d)}
&\leqslant \mathscr{C}(l)s^{-(l+1)/2}e^{-\check{c}_*\varepsilon ^2 s/2}.
\end{align}
The constant $\mathscr{C}(l)$ depends only on $l$ and on initial data \eqref{dannie_zadachi}.
\end{proposition}

Now we can consider the first summand in \eqref{K(eps,s) with Pi}. (The second summand is adjoint to the first one and does not require additional considerations.) By \eqref{Lambda_rho=Lambda+Lambda_rho^0},
\begin{equation}
\label{part corr with Lambda G=}
\begin{split}
\Lambda _G b(\mathbf{D}) \Xi (\varepsilon ,s)
=\Lambda b(\mathbf{D}) \Xi (\varepsilon ,s)
+\Lambda _G^0 b(\mathbf{D})\Xi (\varepsilon ,s).
\end{split}
\end{equation}
Combining \eqref{<b^*b<}, \eqref{f_0<f}, \eqref{f_0Lf_0>=}, \eqref{Lambda_rho^0<=}, and \eqref{Xi (eps,s)}, we get
\begin{equation}
\label{9.8a}
\begin{split}
\Vert \Lambda _G^0 b(\mathbf{D})&\Xi (\varepsilon ,s)\Vert _{L_2(\mathbb{R}^d)\rightarrow L_2(\mathbb{R}^d)}
\\
&\leqslant
\mathfrak{C}_G\Vert f\Vert ^2_{L_\infty}\sup _{\boldsymbol{\xi}\in\mathbb{R}^d}\vert b(\boldsymbol{\xi})\vert e^{-\check{c}_*(\vert \boldsymbol{\xi}\vert ^2 +\varepsilon ^2)s}(1-\chi _{\widetilde{\Omega}}(\boldsymbol{\xi}))
\\
&\leqslant \alpha _1^{1/2}\check{c}_*^{-1}r_0^{-1}\mathfrak{C}_G\Vert f\Vert ^2  _{L_\infty} s^{-1}e^{-\check{c}_*\varepsilon ^2 s/2},\quad s>0.
\end{split}
\end{equation}
In what follows, we are only interested in large values of $s$, so now it suffices to consider  $s\geqslant 1$. 
Lemma~\ref{Lemma Lambda multiplier} and estimate \eqref{b(D)Xi(eps,s)<=} imply that
\begin{equation}
\label{part with Lambda<=}
\Vert \Lambda b(\mathbf{D})\Xi (\varepsilon ,s)\Vert _{L_2(\mathbb{R}^d)\rightarrow L_2(\mathbb{R}^d)}
\leqslant C_{12}s^{-1}e^{-\check{c}_*\varepsilon ^2 s/2},\quad s\geqslant 1.
\end{equation}
Here $C_{12}=\mathfrak{C}_\Lambda \alpha _1^{1/2}\mathscr{C}(l)$, where $l=1$ for $d\leqslant 4$ and $l=d/2-1$ for $d> 4$. 
Combining \eqref{part corr with Lambda G=}--\eqref{part with Lambda<=}, we get 
\begin{equation*}
\Vert \Lambda _G b(\mathbf{D}) \Xi (\varepsilon ,s)\Vert _{L_2(\mathbb{R}^d)\rightarrow L_2(\mathbb{R}^d)}
\leqslant
C_{13}s^{-1}e^{-\check{c}_*\varepsilon ^2 s/2},\quad s\geqslant 1,
\end{equation*}
where $C_{13}=\alpha _1^{1/2}\check{c}_*^{-1}r_0^{-1}\mathfrak{C}_G\Vert f\Vert ^2  _{L_\infty}+C_{12}$.

Similarly, using Lemma~\ref{Lemma tilde Lambda multiplier} and 
\eqref{tildeLambda_rho=tildeLambda+}, \eqref{tildeLambda_rho^0<=}, and \eqref{eps Xi(eps,s)<=}, we arrive at the estimate 
\begin{equation*}
\Vert \varepsilon \widetilde{\Lambda}_G\Xi (\varepsilon ,s)\Vert _{L_2(\mathbb{R}^d)\rightarrow L_2(\mathbb{R}^d)}
\leqslant C_{14}s^{-1}e^{-\check{c}_*\varepsilon ^2 s/2},\quad s\geqslant 1,
\end{equation*}
with the constant $C_{14}=\check{c}_*^{-1}r_0^{-1}\widetilde{\mathfrak{C}}_G\Vert f\Vert ^2_{L_\infty}+\mathfrak{C}_{\widetilde{\Lambda}}\mathscr{C}(\sigma )$, where $\sigma=2$ for $d\leqslant 6$ and $\sigma =d/2-1$ for $d> 6$.  Thus, we have shown that for $s\geqslant 1$ the smoothing operator $\Pi$ can be removed from all the terms of the corrector \eqref{K(eps,s) with Pi}.

\begin{theorem}
Let $\mathcal{B}(\varepsilon)$ and $\mathcal{B}^0(\varepsilon)$ be the operators \eqref{B(e)} and \eqref{B0(eps)}, respectively. 

\noindent
$1^\circ$. 
Let $\mathcal{K}(\varepsilon ,s)$ be the operator \eqref{K(eps,s) with Pi}. Then for $0<\varepsilon \leqslant 1$ and $s\geqslant 0$ we have the estimate 
\begin{equation}
\label{epx B(eps) appr no Pi in principal part Th}
\Vert fe^{-\mathcal{B}(\varepsilon)s}f^*-f_0e^{-\mathcal{B}^0(\varepsilon )s}f_0-\mathcal{K}(\varepsilon,s)\Vert _{L_2(\mathbb{R}^d)\rightarrow L_2(\mathbb{R}^d)}
\leqslant C_{15}(1+s)^{-1}e^{-\check{c}_*\varepsilon ^2 s/2}
\end{equation}
with the constant $C_{15}=C_{11}+2\Vert f\Vert ^2_{L_\infty}\max\lbrace 1;\check{c}_*^{-1}r_0^{-2}\rbrace$.

\noindent
$2^\circ$. 
Denote
\begin{equation}
\label{K(eps,s) without Pi}
\begin{split}
\mathcal{K}^0(\varepsilon,s)
:&=
\left( \Lambda _G b(\mathbf{D})+\varepsilon\widetilde{\Lambda}_G\right)f_0e^{-\mathcal{B}^0(\varepsilon)s}f_0
\\
&+f_0e^{-\mathcal{B}^0(\varepsilon)s}f_0\left( b(\mathbf{D})^*\Lambda _G ^*+\varepsilon\widetilde{\Lambda}_G ^*\right)
\\
&-\int\limits_0^s f_0e^{-\mathcal{B}^0(\varepsilon)(s-\widetilde{s})}f_0\mathcal{N}(\varepsilon)f_0e^{-\mathcal{B}^0(\varepsilon)\widetilde{s}}f_0 \,d\widetilde{s}.
\end{split}
\end{equation}
Then for $0<\varepsilon \leqslant 1$ and $s\geqslant 1$ we have
\begin{equation}
\label{appr exp B(eps) no Pi}
\Vert fe^{-\mathcal{B}(\varepsilon)s}f^*-f_0e^{-\mathcal{B}^0(\varepsilon )s}f_0-\mathcal{K}^0(\varepsilon,s)\Vert _{L_2(\mathbb{R}^d)\rightarrow L_2(\mathbb{R}^d)}
\leqslant C_{16}s^{-1}e^{-\check{c}_*\varepsilon ^2 s/2},
\end{equation}
where $C_{16}=C_{10}+2\Vert f\Vert ^2_{L_\infty}\max\lbrace 1;\check{c}_*^{-1}r_0^{-2}\rbrace +3\check{c}_*^{-2}r_0^{-1}C_G\Vert f\Vert ^4_{L_\infty}+2C_{13}+2C_{14}$.
\end{theorem}

\section*{Chapter 3. Homogenization problem for parabolic systems}

\section{Homogenization of the operator $f^\varepsilon e^{-\mathcal{B}_\varepsilon s} (f^\varepsilon)^*$}

\label{Section main results in operator terms}

\subsection{The problem setting} For any $\Gamma$-periodic function $\phi (\mathbf{x})$, $\mathbf{x}\in\mathbb{R}^d$, we use the notation $\phi ^\varepsilon (\mathbf{x}):=\phi (\mathbf{x}/\varepsilon )$, $\varepsilon >0$. In $L_2(\mathbb{R}^d;\mathbb{C}^n)$, we consider the operator 
\begin{equation*}
\mathcal{A}_\varepsilon =f^\varepsilon (\mathbf{x})^* b(\mathbf{D})^*g^\varepsilon (\mathbf{x})b(\mathbf{D})f^\varepsilon (\mathbf{x}),
\end{equation*}
corresponding to the quadratic form
\begin{equation}
\label{a_eps [u,u]}
\mathfrak{a}_\varepsilon [\mathbf{u},\mathbf{u}]=\int\limits_{\mathbb{R}^d}\langle g^\varepsilon b(\mathbf{D})f^\varepsilon \mathbf{u},b(\mathbf{D})f^\varepsilon \mathbf{u}\rangle\,d\mathbf{x}
\end{equation}
on the domain $\mathfrak{d}_\varepsilon =\lbrace \mathbf{u}\in L_2(\mathbb{R}^d;\mathbb{C}^n) : f^\varepsilon\mathbf{u}\in H^1(\mathbb{R}^d;\mathbb{C}^n)\rbrace$. 
The form \eqref{a_eps [u,u]} is subject to estimates that are similar to inequalities  \eqref{<a[u,u]<}:
\begin{equation}
\label{<a_eps[u,u]<}
\alpha _0\Vert g^{-1}\Vert _{L_\infty}^{-1}\Vert \textbf{D}(f^\varepsilon \textbf{u})\Vert ^2_{L_2(\mathbb{R}^d)}\leqslant \mathfrak{a}_\varepsilon [\textbf{u},\textbf{u}]\leqslant \alpha _1\Vert g \Vert _{L_\infty}\Vert \textbf{D}(f^\varepsilon \textbf{u})\Vert ^2_{L_2(\mathbb{R}^d)}.
\end{equation}

Next, let $\mathcal{Y}_\varepsilon : L_2(\mathbb{R}^d;\mathbb{C}^n)\rightarrow L_2(\mathbb{R}^d;\mathbb{C}^{dn})$ be the~operator, acting by the rule
\begin{equation*}
\mathcal{Y}_\varepsilon\mathbf{u}=\mathbf{D}(f^\varepsilon\mathbf{u})=\mathrm{col}\,\lbrace D_1(f^\varepsilon \mathbf{u}),\dots , D_d (f^\varepsilon\mathbf{u})\rbrace,\quad\mathrm{Dom}\,\mathcal{Y}_\varepsilon =\mathfrak{d}_\varepsilon ,
\end{equation*}
and let $\mathcal{Y}_{2,\varepsilon} : L_2(\mathbb{R}^d;\mathbb{C}^n)\rightarrow L_2(\mathbb{R}^d;\mathbb{C}^{dn})$ be the operator of multiplication by the $(dn\times d)$-matrix-valued function consisting of the blocks $a_j^\varepsilon (\mathbf{x})^*f^\varepsilon (\mathbf{x}) $, $j=1,\dots ,d$:
\begin{equation*}
\mathcal{Y}_{2,\varepsilon}\mathbf{u}=\mathrm{col}\,\lbrace (a_1^\varepsilon )^*f^\varepsilon\mathbf{u},\dots , (a_d^\varepsilon )^*f^\varepsilon\mathbf{u}\rbrace,\quad \mathbf{u}\in \mathfrak{d}_\varepsilon.
\end{equation*}

Let $d\mu$ a matrix-valued measure on $\mathbb{R}^d$ defined in Subsection~\ref{Subsection Q0 q}. By using it, we build the measure  $d\mu ^\varepsilon$ as follows. For any Borel set $\Delta \subset\mathbb{R}^d$, consider the set $\varepsilon ^{-1}\Delta=\lbrace\mathbf{y}=\varepsilon ^{-1}\mathbf{x} : \mathbf{x}\in\Delta\rbrace$ and put $\mu ^\varepsilon (\Delta)=\varepsilon ^d\mu (\varepsilon ^{-1}\Delta )$. Define the quadratic form $q_\varepsilon$ by the rule
\begin{equation*}
q_\varepsilon [\mathbf{u},\mathbf{u}]=\int\limits_{\mathbb{R}^d}\langle d\mu ^\varepsilon (\mathbf{x})(f^\varepsilon\mathbf{u})(\mathbf{x}),(f^\varepsilon\mathbf{u})(\mathbf{x})\rangle,\quad\mathbf{u}\in \mathfrak{d}_\varepsilon.
\end{equation*}

All the assumptions of Subsections~\ref{Subsection Factorized DOs}--\ref{Subsection Q0 q} are assumed to be satisfied. In the space $L_2(\mathbb{R}^d;\mathbb{C}^n)$, consider the quadratic form 
\begin{equation}
\label{b_eps [u,u]}
\mathfrak{b}_\varepsilon [\mathbf{u},\mathbf{u}]=\mathfrak{a}_\varepsilon [\mathbf{u},\mathbf{u}]
+2\mathrm{Re}\,(\mathcal{Y}_\varepsilon \mathbf{u},\mathcal{Y}_{2,\varepsilon}\mathbf{u})_{L_2(\mathbb{R}^d)}+q_\varepsilon [\mathbf{u},\mathbf{u}],\quad\mathbf{u}\in \mathfrak{d}_\varepsilon.
\end{equation}

Let $T_\varepsilon$ be the scaling transformation unitary in $L_2(\mathbb{R}^d;\mathbb{C}^n)$: 
\begin{equation}
\label{T_eps}
(T_\varepsilon \mathbf{u})(\mathbf{x})=\varepsilon ^{d/2}\mathbf{u}(\varepsilon\mathbf{x}),\quad\varepsilon >0.
\end{equation}
The forms \eqref{frak a[u,u]} and \eqref{a_eps [u,u]}, \eqref{b(e)[u,u]} and \eqref{b_eps [u,u]} satisfy the obvious identities
\begin{equation*}
\mathfrak{a}_\varepsilon [\mathbf{u},\mathbf{u}]
=
\varepsilon ^{-2}\mathfrak{a}[T_\varepsilon\mathbf{u},T_\varepsilon \mathbf{u}],\quad 
\mathfrak{b}_\varepsilon [\mathbf{u},\mathbf{u}]
=
\varepsilon ^{-2}\mathfrak{b}(\varepsilon)[T_\varepsilon\mathbf{u},T_\varepsilon \mathbf{u}],\quad \mathbf{u}\in \mathfrak{d}_\varepsilon .
\end{equation*}
Together with estimates \eqref{b(e)[u,u]>=} and \eqref{b(e)[u,u]<=}, this implies that for  $\mathbf{u}\in \mathfrak{d}_\varepsilon$ we have
\begin{align}
\label{b_eps [u,u]>=}
&\mathfrak{b}_\varepsilon [\mathbf{u},\mathbf{u}]\geqslant \frac{\kappa}{2}\mathfrak{a}_\varepsilon [\mathbf{u},\mathbf{u}]+\beta \Vert \mathbf{u}\Vert ^2_{L_2(\mathbb{R}^d)},
 \\
\label{b_eps [u,u]<=}
&\mathfrak{b}_\varepsilon [\mathbf{u},\mathbf{u}]\leqslant (2+c_1^2+c_2)\mathfrak{a}_\varepsilon [\mathbf{u},\mathbf{u}]+(C(1)+c_3+\vert \lambda \vert \Vert \mathcal{Q}_0 \Vert _{L_\infty })\Vert \mathbf{u}\Vert ^2_{L_2(\mathbb{R}^d)}.
\end{align}
By \eqref{<a_eps[u,u]<}, \eqref{b_eps [u,u]>=}, and \eqref{b_eps [u,u]<=}, the form $\mathfrak{b}_\varepsilon$ is closed and positive definite. By $\mathcal{B}_\varepsilon$ we denote the corresponding self-adjoint operator in the space $L_2(\mathbb{R}^d;\mathbb{C}^n)$. Formally,
\begin{equation}
\label{B_eps}
\begin{split}
\mathcal{B}_\varepsilon &=\mathcal{A}_\varepsilon + (\mathcal{Y}^*_{2,\varepsilon}\mathcal{Y}_\varepsilon +\mathcal{Y}^*_\varepsilon \mathcal{Y}_{2,\varepsilon})+(f^\varepsilon)^*\mathcal{Q}^\varepsilon
f^\varepsilon +\lambda \mathcal{Q}_0^\varepsilon 
\\
&=(f^\varepsilon)^*b(\mathbf{D})^*g^\varepsilon b(\mathbf{D})f^\varepsilon + 
\sum \limits _{j=1}^d (f^\varepsilon)^*\left( a_j^\varepsilon D_j+D_j(a_j^\varepsilon)^* \right)f^\varepsilon 
\\
&+ (f^\varepsilon)^*\mathcal{Q}^\varepsilon
f^\varepsilon +\lambda \mathcal{Q}_0^\varepsilon ,
\end{split}
\end{equation}
where $\mathcal{Q}^\varepsilon$ should be interpreted as a generalised matrix-valued potential generated by the rapidly oscillating measure~$d\mu ^\varepsilon$.

\subsection{The principal term of approximation for $f^\varepsilon e^{-\mathcal{B}_\varepsilon s} (f^\varepsilon)^*$} 
The principal term of approximation was found at \cite[Subsection~9.3]{M_AA}. To formulate the result, we need to write out the effective operator:
\begin{equation}
\label{B0}
\begin{split}
\mathcal{B}^0&=
f_0b(\mathbf{D})^*g^0b(\mathbf{D})f_0
\\&+ f_0\bigg(-b(\mathbf{D})^*V-V^*b(\mathbf{D})+\sum _{j=1}^d \overline{(a_j+a_j^*)}D_j\bigg)f_0
\\
&+ f_0(-W+\overline{\mathcal{Q}}+\lambda I)f_0.
\end{split}
\end{equation}
Here the constant matrices $f_0$, $g^0$, $V$, $W$, $(\overline{a_j+a_j^*})$, and $\overline{\mathcal{Q}}$ are defined according to \eqref{f0}, \eqref{g^0}, \eqref{V}, \eqref{W}, \eqref{underline and overline definition}, and \eqref{overlineQ}. The symbol of the operator \eqref{B0} is the matrix $f_0\widehat{L}(\boldsymbol{\xi} ,1)f_0$ (see Subsection~\ref{Subsection B0(k,eps)}). Thus, estimate \eqref{f_0Lf_0>=} implies that the operator $\mathcal{B}^0$ is positive definite: $\mathcal{B}^0\geqslant \check{c}_*I$.

In \cite[Theorem~9.1]{M_AA}, the following result was obtained.

\begin{theorem}
\label{Theorem principal term of approx exp B_eps}
Let the assumptions of Subsections~\textnormal{\ref{Subsection Factorized DOs}--\ref{Subsection Q0 q}} be satisfied. Let $\mathcal{B}_\varepsilon$ be the operator \eqref{B_eps} and let $\mathcal{B}^0$ be the effective operator \eqref{B0}. Then for $s\geqslant 0$ and $0<\varepsilon\leqslant 1$ we have
\begin{equation*}
\Vert f^\varepsilon e^{-\mathcal{B}_\varepsilon s}(f^\varepsilon )^*-f_0e^{-\mathcal{B}^0s}f_0\Vert _{L_2(\mathbb{R}^d)\rightarrow L_2(\mathbb{R}^d)}
\leqslant
C_{17}\varepsilon (s+\varepsilon ^2)^{-1/2}e^{-\check{c}_*s/2}.
\end{equation*}
The constant $C_{17}$ depends only on the initial data \eqref{dannie_zadachi}.
\end{theorem}

\subsection{Approximation for the operator $f^\varepsilon e^{-\mathcal{B}_\varepsilon s} (f^\varepsilon)^*$ with the corrector taken into account}

Now we are going to derive a sharper approximation for the operator exponential from the results of 
\S\ref{Section f exp B(eps) f*} with the help of the scaling transformation \eqref{T_eps}. The operators \eqref{B_eps} and \eqref{B(e)}, \eqref{B0} and \eqref{B0(eps)} 
are related by the identities
$\mathcal{B}_\varepsilon =\varepsilon ^{-2}T_\varepsilon ^* \mathcal{B}(\varepsilon)T_\varepsilon$ and 
$\mathcal{B}^0=\varepsilon ^{-2}T_\varepsilon ^*\mathcal{B}^0(\varepsilon)T_\varepsilon $. 
Thus,
\begin{align}
\label{exp scaling}
&f^\varepsilon e^{-\mathcal{B}_\varepsilon s} (f^\varepsilon)^*=T_\varepsilon ^* f e^{-\mathcal{B}(\varepsilon)\varepsilon ^{-2}s}f^*T_\varepsilon ,
\\
\label{eff exp scaling}
& f_0e^{-\mathcal{B}^0s}f_0=T_\varepsilon ^* f_0e^{-\mathcal{B}^0(\varepsilon)\varepsilon ^{-2}s}f_0 T_\varepsilon .
\end{align}

Next, by $\Pi _\varepsilon$ we define a PDO acting in $L_2(\mathbb{R}^d;\mathbb{C}^n)$ whose symbol is the characteristic function $\chi _{\widetilde{\Omega}/\varepsilon}(\boldsymbol{\xi})$ of the set $\widetilde{\Omega}/\varepsilon$:
\begin{equation}
\label{Pi_eps}
(\Pi _\varepsilon \mathbf{f})(\mathbf{x})=(2\pi )^{-d/2}\int\limits_{\widetilde{\Omega}/\varepsilon}e^{i\langle \mathbf{x},\boldsymbol{\xi}\rangle}(\mathcal{F}\mathbf{f})(\boldsymbol{\xi})\,d\boldsymbol{\xi}.
\end{equation}
The operators \eqref{Pi} and \eqref{Pi_eps} are related via the formula 
\begin{equation}
\label{Pi scaling tozd}
\Pi _\varepsilon =T_\varepsilon ^*\Pi T_\varepsilon.
\end{equation}

Denote 
\begin{equation*}
\mathcal{N}:=\mathcal{N}_{11}(\mathbf{D})+\mathcal{N}_{12}(\mathbf{D})+\mathcal{N}_{21}(\mathbf{D})+\mathcal{N}_{22}.
\end{equation*}
Then
\begin{equation}
\label{N tozd scaling}
\mathcal{N}=\varepsilon ^{-3}T_\varepsilon ^*\mathcal{N}(\varepsilon)T_\varepsilon ,
\end{equation}
where the operator $\mathcal{N}(\varepsilon )$ was defined in \eqref{9.1a}.

We point out the identities
\begin{equation}
\label{Lambda i tilde-Lambda scaling tozd}
[\Lambda _G ^\varepsilon]b(\mathbf{D})=\varepsilon ^{-1}T_\varepsilon ^*[\Lambda _G]b(\mathbf{D})T_\varepsilon,
\quad
[\widetilde{\Lambda} _G^\varepsilon]=T_\varepsilon ^*[\widetilde{\Lambda}_G]T_\varepsilon .
\end{equation}

We introduce the correctors:
{\allowdisplaybreaks
\begin{align}
\label{K_eps (s) with Pi}
\begin{split}
\mathcal{K}_\varepsilon(s)
:&=
\left( \Lambda _G ^\varepsilon b(\mathbf{D})+\widetilde{\Lambda}_G^\varepsilon\right)f_0e^{-\mathcal{B}^0 s}f_0\Pi _\varepsilon
\\
&+f_0e^{-\mathcal{B}^0s}f_0\Pi _\varepsilon \left( b(\mathbf{D})^*(\Lambda _G^\varepsilon) ^*+(\widetilde{\Lambda}_G^\varepsilon)^*\right)
\\
&-\int\limits_0^s f_0e^{-\mathcal{B}^0(s-\widetilde{s})}f_0\mathcal{N}f_0e^{-\mathcal{B}^0\widetilde{s}}f_0\Pi _\varepsilon \,d\widetilde{s} ,
\end{split}
\\
\label{K_eps 0(s) no Pi}
\begin{split}
\mathcal{K}_\varepsilon ^0(s)
:&=
\left( \Lambda _G ^\varepsilon b(\mathbf{D})+\widetilde{\Lambda}_G^\varepsilon\right)f_0e^{-\mathcal{B}^0 s}f_0
\\
&+f_0e^{-\mathcal{B}^0s}f_0 \left( b(\mathbf{D})^*(\Lambda _G^\varepsilon) ^*+(\widetilde{\Lambda}_G^\varepsilon)^*\right)
\\
&-\int\limits_0^s f_0e^{-\mathcal{B}^0(s-\widetilde{s})}f_0\mathcal{N}f_0e^{-\mathcal{B}^0\widetilde{s}}f_0 \,d\widetilde{s}.
\end{split}
\end{align}
}

Combining \eqref{exp scaling}, \eqref{eff exp scaling}, and \eqref{Pi scaling tozd}--\eqref{Lambda i tilde-Lambda scaling tozd}, we conclude that the operators~\eqref{K(eps,s) with Pi} and \eqref{K_eps (s) with Pi} are related by the identity
\begin{equation}
\label{K scaling}
\varepsilon \mathcal{K}_\varepsilon (s)=T_\varepsilon ^* \mathcal{K}(\varepsilon ,\varepsilon ^{-2}s)T_\varepsilon .
\end{equation}
For the operators \eqref{K(eps,s) without Pi} and \eqref{K_eps 0(s) no Pi}, a similar identity holds true. 
Now \eqref{epx B(eps) appr no Pi in principal part Th}, \eqref{appr exp B(eps) no Pi}, \eqref{exp scaling}, \eqref{eff exp scaling}, and \eqref{K scaling} imply the following result.

\begin{theorem}
\label{Theorem exp B_eps final}
Let the assumptions of Theorem~\textnormal{\ref{Theorem principal term of approx exp B_eps}} be satisfied. Let $\mathcal{K}_\varepsilon (s)$ and $\mathcal{K}_\varepsilon ^0(s)$ be the operators \eqref{K_eps (s) with Pi} and \eqref{K_eps 0(s) no Pi}, respectively. Then for $0<\varepsilon \leqslant 1$ and $s\geqslant 0$ we have the approximation
\begin{equation}
\label{Th Main result s>=0}
\begin{split}
\Vert f^\varepsilon e^{-\mathcal{B}_\varepsilon s}(f^\varepsilon )^*&-f_0e^{-\mathcal{B}^0s}f_0 -\varepsilon \mathcal{K}_\varepsilon (s)\Vert _{L_2(\mathbb{R}^d)\rightarrow L_2(\mathbb{R}^d)}
\\
&\leqslant C_{15}\varepsilon ^2 (s+\varepsilon ^2)^{-1}e^{-\check{c}_* s/2}.
\end{split}
\end{equation}
For $0<\varepsilon \leqslant 1$ and $s\geqslant \varepsilon ^2$ we have
\begin{equation*}
\Vert f^\varepsilon e^{-\mathcal{B}_\varepsilon s}(f^\varepsilon )^*-f_0e^{-\mathcal{B}^0s}f_0 -\varepsilon \mathcal{K}_\varepsilon ^0 (s)\Vert _{L_2(\mathbb{R}^d)\rightarrow L_2(\mathbb{R}^d)}
\leqslant C_{16}\varepsilon ^2 s^{-1}e^{-\check{c}_* s/2}.
\end{equation*}
The constants $C_{15}$, $C_{16}$, and $\check{c}_*$ depend only on the initial data \eqref{dannie_zadachi}.
\end{theorem}

\subsection{The case when the corrector is equal to zero}

Assume that $g^0=\overline{g}$, i.~e., relations \eqref{g0=overline g} hold true. Then the $\Gamma$-periodic solutions of problems \eqref{Lambda}, \eqref{Lambda-G problem} are equal to zero: $\Lambda (\mathbf{x})=\Lambda _G(\mathbf{x})=0$. In addition, we assume that $\sum _{j=1}^d D_j a_j(\mathbf{x})^*=0$. Then the $\Gamma$-periodic solutions of problems \eqref{tildeLambda}, \eqref{tilde-Lambda-G problem} are also equal to zero: $\widetilde{\Lambda} (\mathbf{x})=\widetilde{\Lambda} _G(\mathbf{x})=0$. Thus, \eqref{V} and \eqref{W} imply that $V=0$ and $W=0$. The effective operator \eqref{B0} takes the form
\begin{equation}
\label{B0 with V=W=0}
\begin{split}
\mathcal{B}^0=
f_0b(\mathbf{D})^*g^0b(\mathbf{D})f_0+f_0\sum _{j=1}^d \overline{(a_j+a_j^*)}D_jf_0
+ f_0(\overline{\mathcal{Q}}+\lambda I)f_0.
\end{split}
\end{equation}
Using \eqref{N tozd scaling} and Remark~\ref{Remark N(eps)=0}, we conclude that $\mathcal{N}=0$ in the case under consideration. Thus, all the corrector terms in \eqref{K_eps (s) with Pi} are equal to zero and \eqref{Th Main result s>=0} implies the following result.

\begin{proposition}
Under the assumptions of Theorem~\textnormal{\ref{Theorem exp B_eps final}}, let relations  \eqref{g0=overline g} and \break$\sum _{j=1}^d D_j a_j(\mathbf{x})^*=0$ hold true. Then for $0<\varepsilon \leqslant 1$ and $s\geqslant 0$ we have
\begin{equation*}
\Vert f^\varepsilon e^{-\mathcal{B}_\varepsilon s}(f^\varepsilon )^*-f_0e^{-\mathcal{B}^0s}f_0 \Vert _{L_2(\mathbb{R}^d)\rightarrow L_2(\mathbb{R}^d)}
\leqslant C_{15}\varepsilon ^2 (s+\varepsilon ^2)^{-1}e^{-\check{c}_* s/2}.
\end{equation*}
Here $\mathcal{B}^0$ is the operator \eqref{B0 with V=W=0}.
\end{proposition}

\section{Homogenization for solutions of the Cauchy problem for parabolic systems}

\label{Section Homogenization of parabolic systems}

\subsection{Application to the Cauchy problem for a~homogeneous equation} Let $\widehat{\mathcal{B}}_\varepsilon $ be the operator of the form \eqref{B_eps} with $f=\mathbf{1}_n$:
\begin{equation*}
\begin{split}
\widehat{\mathcal{B}}_\varepsilon =
b(\mathbf{D})^*g^\varepsilon b(\mathbf{D})+ 
\sum \limits _{j=1}^d \left( a_j^\varepsilon D_j+D_j(a_j^\varepsilon)^* \right)
+ \mathcal{Q}^\varepsilon
 +\lambda I.
\end{split}
\end{equation*}

Consider the Cauchy problem
\begin{equation}
\label{Cauchy problem}
G^\varepsilon (\mathbf{x})\frac{\partial \mathbf{u}_\varepsilon (\mathbf{x},s)}{\partial s}=-\widehat{\mathcal{B}}_\varepsilon \mathbf{u}_\varepsilon (\mathbf{x},s),\quad s>0;
\quad
G^\varepsilon (\mathbf{x})\mathbf{u}_\varepsilon (\mathbf{x},0)=\boldsymbol{\phi}(\mathbf{x}),
\end{equation}
where  $\boldsymbol{\phi}\in L_2(\mathbb{R}^d;\mathbb{C}^n)$, and the $\Gamma$-periodic $(n\times n)$-matrix-valued function $G$ is assumed to be bounded and positive definite. We rewrite $G(\mathbf{x})$ in a factorised form: $G(\mathbf{x})^{-1}=f(\mathbf{x})f(\mathbf{x})^*$, where $f(\mathbf{x})$ is a $\Gamma$-periodic matrix-valued function. 
Put $\mathbf{v}_\varepsilon =(f^\varepsilon )^{-1}\mathbf{u}_\varepsilon$. Then $\mathbf{v}_\varepsilon (\mathbf{x},s)$ is the solution of the problem 
\begin{equation*}
\frac{\partial \mathbf{v}_\varepsilon (\mathbf{x},s)}{ \partial s}=-(f^\varepsilon )^*\widehat{\mathcal{B}}_\varepsilon f^\varepsilon\mathbf{v}_\varepsilon (\mathbf{x},s),\quad s>0;\quad 
\mathbf{v}_\varepsilon (\mathbf{x},0)=(f^\varepsilon )^*\boldsymbol{\phi}(\mathbf{x}).
\end{equation*}
Let $\mathcal{B}_\varepsilon =(f^\varepsilon)^*\widehat{\mathcal{B}}_\varepsilon f^\varepsilon$. Then $\mathbf{v}_\varepsilon (\,\cdot\, ,s)=e^{-\mathcal{B}_\varepsilon s}(f^\varepsilon )^*\boldsymbol{\phi}$. 
This implies that the solution $\mathbf{u}_\varepsilon =f^\varepsilon \mathbf{v}_\varepsilon $ of the problem \eqref{Cauchy problem} admits the representation
\begin{equation}
\label{u-eps= F=0}
\mathbf{u}_\varepsilon (\,\cdot\, ,s)=f^\varepsilon e^{-\mathcal{B}_\varepsilon s}(f^\varepsilon )^*\boldsymbol{\phi}.
\end{equation}

The corresponding effective problem has the form
\begin{equation}
\label{eff Cauchy problem}
\overline{G}\frac{\partial \mathbf{u}_0(\mathbf{x},s)}{\partial s}=-\widehat{\mathcal{B}}^0\mathbf{u}_0(\mathbf{x},s),\quad s>0;\quad \overline{G}\mathbf{u}_0(\mathbf{x},0)=\boldsymbol{\phi}(\mathbf{x}). 
\end{equation}
Here
\begin{equation*}
\begin{split}
\widehat{\mathcal{B}}^0=
b(\mathbf{D})^*g^0b(\mathbf{D}) -b(\mathbf{D})^*V-V^*b(\mathbf{D})+\sum _{j=1}^d \overline{(a_j+a_j^*)}D_j 
 -W+\overline{\mathcal{Q}}+\lambda I.
\end{split}
\end{equation*}
Let $f_0=(\overline{G})^{-1/2}$ and $ \mathcal{B}^0 =f_0\widehat{\mathcal{B}}^0 f_0$. 
Then
\begin{equation}
\label{u-0= F=0}
\mathbf{u}_0(\,\cdot\, ,s)=f_0e^{-\mathcal{B}^0 s}f_0\boldsymbol{\phi}.
\end{equation}
From \eqref{u-eps= F=0}, \eqref{u-0= F=0}, and Theorem~\ref{Theorem exp B_eps final} we derive the following result.

\begin{theorem}
Let $\mathbf{u}_\varepsilon $ and $\mathbf{u}_0$ be the solutions of problems \eqref{Cauchy problem} and \eqref{eff Cauchy problem}, respectively. Then for $0<\varepsilon\leqslant 1$ and $s\geqslant \varepsilon ^2$ we have the approximation
\begin{equation}
\label{Th solutions homogenization}
\begin{split}
\Bigl\Vert  \mathbf{u}_\varepsilon (\,\cdot\, ,s)&-\mathbf{u}_0(\,\cdot\, ,s)
-\varepsilon 
\Bigl(\Lambda _G^\varepsilon b(\mathbf{D})\mathbf{u}_0(\,\cdot\, ,s)+\widetilde{\Lambda}_G^\varepsilon \mathbf{u}_0(\,\cdot\, ,s)
\\
&+f_0e^{-\mathcal{B}^0s}f_0\bigl((\Lambda _G^\varepsilon b(\mathbf{D}))^*+(\widetilde{\Lambda}_G^\varepsilon)^*\bigr)\boldsymbol{\phi}
\\
&-\int\limits_0^s f_0e^{-\mathcal{B}^0(s-\widetilde{s})}f_0\mathcal{N}\mathbf{u}_0(\,\cdot\, ,\widetilde{s})\,d\widetilde{s}
\Bigr)
\Bigr\Vert _{L_2(\mathbb{R}^d)}
\\
&\leqslant C_{16}\varepsilon ^2 s^{-1}e^{-\check{c}_* s/2}\Vert \boldsymbol{\phi}\Vert _{L_2(\mathbb{R}^d)}.
\end{split}
\end{equation}
\end{theorem}

Note that the third term of the corrector in \eqref{Th solutions homogenization}, i. e., the function
$$
\mathbf{v}^{(3)}(\,\cdot\, ,s):=\int\limits_0^s f_0e^{-\mathcal{B}^0(s-\widetilde{s})}f_0\mathcal{N}\mathbf{u}_0(\,\cdot\, ,\widetilde{s})\,d\widetilde{s},
$$
is the solution of the problem
\begin{equation*}
\overline{G}\frac{\partial \mathbf{v}^{(3)}(\mathbf{x},s)}{\partial s}=-\widehat{\mathcal{B}}^0\mathbf{v}^{(3)}(\mathbf{x},s)+\mathcal{N}\mathbf{u}_0(\mathbf{x},s)
,\quad s>0;\quad \overline{G}\mathbf{v}^{(3)}(\mathbf{x},0)=0. 
\end{equation*}

\subsection{The Cauchy problem for an inhomogeneous equation}

Now, we consider a more general Cauchy problem
\begin{equation}
\label{Problem with F}
\begin{split}
G^\varepsilon (\mathbf{x})\frac{\partial \mathbf{u}_\varepsilon (\mathbf{x},s)}{\partial s}&=-\widehat{\mathcal{B}}_\varepsilon\mathbf{u}_\varepsilon (\mathbf{x},s)+\mathbf{F}(\mathbf{x},s),\quad s\in (0,T);
\\
G^\varepsilon (\mathbf{x})\mathbf{u}_\varepsilon (\mathbf{x},0)&=\boldsymbol{\phi}(\mathbf{x}),
\end{split}
\end{equation}
where $\mathbf{x}\in\mathbb{R}^d$, $s\in (0,T)$, $T>0$. Let
$$
\boldsymbol{\phi}\in L_2(\mathbb{R}^d;\mathbb{C}^n),\quad \mathbf{F} \in L_p((0,T);L_2(\mathbb{R}^d;\mathbb{C}^n))=:\mathcal{H}_p
$$
for some $1<p\leqslant\infty$. Let $\mathbf{u}_0$ be the solution of the ``homogenized'' problem
\begin{equation}
\label{effective Problem with F}
\overline{G}\frac{\partial \mathbf{u}_0 (\mathbf{x},s)}{\partial s}=-\widehat{\mathcal{B}}^0\mathbf{u}_0 (\mathbf{x},s)+\mathbf{F}(\mathbf{x},s),\quad s\in (0,T);
\quad
\overline{G}\mathbf{u}_0(\mathbf{x},0)=\boldsymbol{\phi}(\mathbf{x}).
\end{equation}

By analogy with the proof of representation \eqref{u-eps= F=0}, one can shown that
\begin{align}
\label{u-eps= with F}
&\mathbf{u}_\varepsilon (\,\cdot\, ,s)= f^\varepsilon e^{-\mathcal{B}_\varepsilon s}(f^\varepsilon )^*\boldsymbol{\phi}+\int\limits_0^s f^\varepsilon e^{-\mathcal{B}_\varepsilon (s-\widetilde{s})}(f^\varepsilon )^*\mathbf{F}(\,\cdot\, ,\widetilde{s})\,d\widetilde{s},
\\
\label{u-0 = with F}
&\mathbf{u}_0(\,\cdot\, ,s) =f_0 e^{-\mathcal{B}^0 s}f_0\boldsymbol{\phi}
+\int\limits_0^s f_0e^{-\mathcal{B}^0 (s-\widetilde{s})}f_0\mathbf{F}(\,\cdot\, ,\widetilde{s})\,d\widetilde{s}.
\end{align}

The principal term of approximation for the solution $\mathbf{u}_\varepsilon$ was obtained in \cite[Theorem~10.1]{M_AA}:
\begin{equation}
\label{Th 10.1 from M}
\begin{split}
\Vert \mathbf{u}_\varepsilon (\,\cdot\, ,s)-\mathbf{u}_0(\,\cdot\, ,s)\Vert _{L_2(\mathbb{R}^d)}
&\leqslant
C_{17}\varepsilon (s+\varepsilon ^2 )^{-1/2}e^{-\check{c}_* s/2}\Vert \boldsymbol{\phi}\Vert _{L_2(\mathbb{R}^d)}
\\
&+C_{18}\theta _1(\varepsilon ,p)\Vert \mathbf{F}\Vert _{\mathcal{H}_p},
\end{split}
\end{equation}
where the constant $C_{18}$ depends only on $p$ and the problem data \eqref{dannie_zadachi},
\begin{equation}
\label{theta-1 (eps,p)}
\theta _1(\varepsilon ,p)=
\begin{cases}
\varepsilon ^{2-2/p}, &1<p<2,
\\
\varepsilon(1+\vert \ln\varepsilon \vert )^{1/2} , &p=2,
\\
\varepsilon , &2<p\leqslant\infty .
\end{cases}
\end{equation}

Repeating considerations of \cite[Subsection~10.1]{V}, from Theorem~\ref{Theorem exp B_eps final} and representations \eqref{u-eps= with F}, \eqref{u-0 = with F} we derive the following result.

\begin{theorem}
\label{Theorem solutions F neq 0}
Let $\mathbf{u}_\varepsilon$ and $\mathbf{u}_0$ be the solutions of problems \eqref{Problem with F} and \eqref{u-0 = with F}, respectively, where $\boldsymbol{\phi}\in L_2(\mathbb{R}^d;\mathbb{C}^n)$ and $\mathbf{F} \in \mathcal{H}_p$. Then for $0<s<T$, $0<\varepsilon\leqslant 1$ and $2<p\leqslant\infty$ we have
\begin{equation*}
\begin{split}
\mathbf{u}_\varepsilon (\,\cdot\, ,s)
&=\mathbf{u}_0(\,\cdot\, ,s)
+\varepsilon \Bigl(\Lambda _G^\varepsilon b(\mathbf{D})\Pi _\varepsilon \mathbf{u}_0(\,\cdot\, ,s) +\widetilde{\Lambda}_G^\varepsilon \Pi _\varepsilon \mathbf{u}_0(\,\cdot\, ,s)
\\
&+f_0e^{-\mathcal{B}^0s}f_0\bigl( (\Lambda _G^\varepsilon b(\mathbf{D})\Pi _\varepsilon )^* +(\widetilde{\Lambda}_G^\varepsilon \Pi _\varepsilon )^*\bigr)\boldsymbol{\phi}\Bigr)
\\
&-\varepsilon \int\limits_0^s f_0e^{-\mathcal{B}^0(s-\widetilde{s})}f_0\mathcal{N}\Pi _\varepsilon
f_0e^{-\mathcal{B}^0\widetilde{s}}f_0 \boldsymbol{\phi}\,d\widetilde{s}
\\
&+\varepsilon \int\limits_0^s f_0e^{-\mathcal{B}^0(s-\widetilde{s})}f_0 \Bigl(
(\Lambda ^\varepsilon _G b(\mathbf{D})\Pi _\varepsilon )^* +(\widetilde{\Lambda}^\varepsilon _G\Pi _\varepsilon )^*\Bigr)\mathbf{F}(\,\cdot\, ,\widetilde{s})\,d\widetilde{s}
\\
&-\varepsilon \int\limits_0^s \,d\widetilde{s}\int\limits_0^{s-\widetilde{s}}\,ds' f_0e^{-\mathcal{B}^0(s-\widetilde{s}-s')}f_0
\mathcal{N}\Pi _\varepsilon f_0e^{-\mathcal{B}^0s'}f_0 \mathbf{F}(\,\cdot\, ,\widetilde{s})+\widetilde{\mathbf{u}}_\varepsilon(\,\cdot\, ,s),
\end{split}
\end{equation*}
where the remainder term $\widetilde{\mathbf{u}}_\varepsilon(\,\cdot\, ,s)$ admits the estimate
\begin{equation*}
\Vert \widetilde{\mathbf{u}}_\varepsilon(\,\cdot\, ,s)\Vert _{L_2(\mathbb{R}^d)}
\leqslant
C_{15}\varepsilon ^2(s+\varepsilon ^2)^{-1}e^{-\check{c}_*s/2}\Vert \boldsymbol{\phi}\Vert _{L_2(\mathbb{R}^d)}
+C_{19}\varepsilon ^{2/p'}\theta _2(\varepsilon ,p)\Vert\mathbf{F}\Vert _{\mathcal{H}_p},
\end{equation*}
$p^{-1}+(p')^{-1}=1$. 
Here the constant $C_{19}$ depends only on $p$ and the initial data \eqref{dannie_zadachi}, $\theta _2(\varepsilon ,p)=1$ for $p<\infty$ and $\theta _2(\varepsilon ,p)=1+\vert\ln\varepsilon\vert$ for $p=\infty$.

For $0<s<T$, $0<\varepsilon\leqslant 1$, and $1<p\leqslant 2$ we have the representation
\begin{equation*}
\begin{split}
\mathbf{u}_\varepsilon (\,\cdot\, ,s)&=\mathbf{u}_0(\,\cdot\, ,s)
+\varepsilon \Bigl(\Lambda _G^\varepsilon b(\mathbf{D})\Pi _\varepsilon f_0e^{-\mathcal{B}^0s}f_0 \boldsymbol{\phi} +\widetilde{\Lambda}_G^\varepsilon \Pi _\varepsilon f_0e^{-\mathcal{B}^0s}f_0 \boldsymbol{\phi}
\\
&+f_0e^{-\mathcal{B}^0s}f_0\bigl( (\Lambda _G^\varepsilon b(\mathbf{D})\Pi _\varepsilon )^* +(\widetilde{\Lambda}_G^\varepsilon \Pi _\varepsilon )^*\bigr)\boldsymbol{\phi}\Bigr)
\\
&-\varepsilon \int\limits_0^s f_0e^{-\mathcal{B}^0(s-\widetilde{s})}f_0\mathcal{N}\Pi _\varepsilon
f_0e^{-\mathcal{B}^0\widetilde{s}}f_0 \boldsymbol{\phi}\,d\widetilde{s}+\widetilde{\mathbf{u}}_\varepsilon(\,\cdot\, ,s),
\end{split}
\end{equation*}
and
\begin{equation*}
\Vert \widetilde{\mathbf{u}}_\varepsilon(\,\cdot\, ,s)\Vert _{L_2(\mathbb{R}^d)}
\leqslant
C_{15}\varepsilon ^2(s+\varepsilon ^2)^{-1}e^{-\check{c}_*s/2}\Vert \boldsymbol{\phi}\Vert _{L_2(\mathbb{R}^d)}
+C_{18}\theta _1(\varepsilon ,p)\Vert \mathbf{F}\Vert _{\mathcal{H}_p}.
\end{equation*}
Here $\theta _1(\varepsilon ,p)$ is quantity \eqref{theta-1 (eps,p)}.
\end{theorem}

\begin{remark}
The different approximations of the solution $\mathbf{u}_\varepsilon$ for different conditions on $p$ in Theorem~\textnormal{\ref{Theorem solutions F neq 0}} are caused by the fact that the corrector term  taken into account in the first summand on the right-hand side of \eqref{u-eps= with F} makes the approximation more precise compared to \eqref{Th 10.1 from M} for any $1<p\leqslant \infty$, while the corrector taken into account in the second summand  does the same only for $2<p\leqslant\infty$.
\end{remark}

\section{Example of application of the general results: \\the scalar elliptic operator}

\label{Section Example}

In the present section, we consider an example of application of the general method. For elliptic problems, this example was earlier studied in \cite{Su10,Su14}. Other examples also can  be found there.

\subsection{The scalar elliptic operator}
\label{Subsection scalar elliptic operator}

Consider the case, when $n=1$, $m=d$, $b(\mathbf{D})=\mathbf{D}$, and $g(\mathbf{x})$ is a  $\Gamma$-periodic symmetric $(d\times d)$-matrix with \textit{real entries}, that is bounded and positive definite. Then the operator $\mathcal{A}_\varepsilon$ has the form 
$$
\mathcal{A}_\varepsilon =\mathbf{D}^*g^\varepsilon (\mathbf{x})\mathbf{D}=-\mathrm{div}\,g^\varepsilon (\mathbf{x})\nabla.
$$
Obviously, in the case under consideration, $\alpha _0=\alpha _1=1$; see \eqref{<b^*b<}.

Next, let $\mathbf{A}(\mathbf{x})=\mathrm{col}\,\lbrace A_1(\mathbf{x}),\dots ,A_d(\mathbf{x})\rbrace$, where $A_j(\mathbf{x})$ are $\Gamma$-periodic real-valued functions, and
\begin{equation*}
A_j\in L_\varrho (\Omega),\quad \varrho =2\;\mbox{for}\;d=1,\quad \varrho >d\;\mbox{for}\;d\geqslant 2;\quad j=1,\dots,d.
\end{equation*}
Let $v(\mathbf{x})$ and $\mathcal{V}(\mathbf{x})$ be real-valued $\Gamma$-periodic functions such that
\begin{equation*}
v,\mathcal{V}\in L_\sigma(\Omega),
\quad
\sigma =1\;\mbox{for}\;d=1,\quad \sigma >\frac{d}{2}\;\mbox{for}\;d\geqslant 2;
\quad \int\limits_\Omega v(\mathbf{x})\,d\mathbf{x}=0.
\end{equation*}

In $L_2(\mathbb{R}^d)$, we consider the operator $\mathfrak{B}_\varepsilon$, formally given by the differential expression
\begin{equation}
\label{example of B-eps}
\mathfrak{B}_\varepsilon =
(\mathbf{D}-\mathbf{A}^\varepsilon (\mathbf{x}))^*g^\varepsilon (\mathbf{x})(\mathbf{D}-\mathbf{A}^\varepsilon (\mathbf{x}))
+\varepsilon ^{-1}v^\varepsilon (\mathbf{x})+\mathcal{V}^\varepsilon (\mathbf{x}).
\end{equation}
The precise definition of the operator $\mathfrak{B}_\varepsilon$ is given via the quadratic form
$$
\mathfrak{b}_\varepsilon [u,u]\!=\!\!\int\limits_{\mathbb{R}^d}
\Bigl(\langle g^\varepsilon (\mathbf{D}\!-\!\mathbf{A}^\varepsilon )u,(\mathbf{D}\!-\!\mathbf{A}^\varepsilon)u\rangle
\!+\!(\varepsilon ^{-1}v^\varepsilon \!+\!\mathcal{V}^\varepsilon )\vert u\vert ^2\Bigr)\,d\mathbf{x},
\quad
u\in H^1(\mathbb{R}^d).
$$
The operator \eqref{example of B-eps} can be treated as a periodic Schr\"odinger operator with the metric $g^\varepsilon$, the magnetic potential $\mathbf{A}^\varepsilon$, and the electric potential $\varepsilon ^{-1}v^\varepsilon +\mathcal{V}^\varepsilon$ containing a ``singular'' first summand.

In \cite[Subsection~13.1]{Su10}, it was obtained that the operator  \eqref{example of B-eps} can be written in a required form
\begin{equation*}
\mathfrak{B}_\varepsilon =\mathbf{D}^*g^\varepsilon (\mathbf{x})\mathbf{D}
+\sum _{j=1}^d \bigl(a_j^\varepsilon (\mathbf{x})D_j+D_j(a_j^\varepsilon (\mathbf{x}))^*\bigr)+\mathcal{Q}^\varepsilon (\mathbf{x}).
\end{equation*}
The real-valued function $\mathcal{Q}(\mathbf{x})$ is defined as
\begin{equation}
\label{Q scalar example}
\mathcal{Q}(\mathbf{x})=\mathcal{V}(\mathbf{x})+\langle g(\mathbf{x})\mathbf{A}(\mathbf{x}),\mathbf{A}(\mathbf{x})\rangle.
\end{equation}
The complex-valued functions $a_j(\mathbf{x})$ are given by
\begin{equation}
\label{a_j example}
a_j(\mathbf{x})=-\eta _j(\mathbf{x})+i\zeta _j(\mathbf{x}),
\end{equation}
where the functions $\eta _j(\mathbf{x})$ are the components of the vector-valued function $\boldsymbol{\eta}(\mathbf{x})=g(\mathbf{x})\mathbf{A}(\mathbf{x})$ and the functions $\zeta _j(\mathbf{x})$ are defined in terms of the $\Gamma$-periodic solution of the equation $\Delta\Phi (\mathbf{x})=v(\mathbf{x})$ by the relation $\zeta _j(\mathbf{x})=-\partial _j\Phi (\mathbf{x})$. Moreover,
\begin{equation}
\label{v= example}
v(\mathbf{x})=-\sum _{j=1}^d \partial _j\zeta _j(\mathbf{x}).
\end{equation}
It is easily seen that the functions \eqref{a_j example} satisfy the assumption \eqref{a_j} for a~suitable exponent $\widetilde{\varrho}$ (depending on $\varrho$ and $\sigma$); herewith, the norms $\Vert a_j\Vert _{L_{\widetilde{\varrho}}(\Omega)}$ are controlled in terms of $\Vert g\Vert _{L_\infty}$, $\Vert \mathbf{A}\Vert _{L_\varrho (\Omega)}$, $\Vert v\Vert _{L_\sigma(\Omega)}$, and the parameters of the lattice $\Gamma$. The function \eqref{Q scalar example} is subject to condition \eqref{Q in Ls ex.} with a suitable factor $\widetilde{\sigma }=\min\lbrace\sigma ;\varrho /2\rbrace$. Thus, now Example~\ref{Example Q in Ls} is implemented.

Let the parameter $\lambda$ be chosen in accordance with condition  \eqref{DOlambda} with $c_0$ and $c_4$ corresponding to the operator \eqref{example of B-eps}. Denote $\mathcal{B}_\varepsilon :=\mathfrak{B}_\varepsilon +\lambda I$. We are interested in approximation of the exponential for this operator. In the case under consideration, initial data \eqref{dannie_zadachi} match with the data:
\begin{equation}
\label{problem data for scalar elliptic operator}
\begin{split}
d,\varrho,\sigma ; \Vert g\Vert _{L_\infty}, \Vert g^{-1}\Vert _{L_\infty}, \Vert \mathbf{A}\Vert _{L_\varrho (\Omega)}, \Vert v\Vert _{L_\sigma(\Omega)}, \Vert \mathcal{V}\Vert _{L_\sigma(\Omega)},\lambda ;
\\
\mbox{the parameters of the lattice }\Gamma.
\end{split}
\end{equation}

\subsection{The effective operator}

\label{Subsection eff op for scalar case}

We write out the effective operator. In our case, the $\Gamma$-periodic solution of problem \eqref{Lambda} is a row-matrix: $\Lambda (\mathbf{x})=i\Psi (\mathbf{x})$, $\Psi (\mathbf{x})=(\psi _1(\mathbf{x}),\dots,\psi _d(\mathbf{x}))$, where $\psi _j\in \widetilde{H}^1(\Omega)$ is the solution of the problem
\begin{equation*}
\mathrm{div}\,g(\mathbf{x})(\nabla \psi _j(\mathbf{x})+\mathbf{e}_j)=0,\quad
\int\limits_\Omega \psi _j(\mathbf{x})\,d\mathbf{x}=0.
\end{equation*}
Here $\mathbf{e}_j$, $j=1,\dots,d$, form the standard basis in $\mathbb{R}^d$. It is clear that the functions $\psi _j(\mathbf{x})$ are real-valued and the elements of the row-matrix $\Lambda (\mathbf{x})$ are purely imaginary. By \eqref{g^0}, the columns of the $(d\times d)$-matrix-valued function $\widetilde{g}(\mathbf{x})$ are the vector-valued functions $g(\mathbf{x})(\nabla \psi _j(\mathbf{x})+\mathbf{e}_j)$, $j=1,\dots,d$. The effective matrix is defined in accordance with \eqref{g^0}: $g^0=\vert \Omega\vert ^{-1}\int _\Omega \widetilde{g}(\mathbf{x})\,d\mathbf{x}$. It is clear that $\widetilde{g}(\mathbf{x})$ and $g^0$ have real entries.

By \eqref{a_j example} and \eqref{v= example}, the periodic solution of problem \eqref{tildeLambda} can be represented as $\widetilde{\Lambda}(\mathbf{x})=\widetilde{\Lambda}_1(\mathbf{x})+i\widetilde{\Lambda}_2(\mathbf{x})$, where the real-valued $\Gamma$-periodic functions $\widetilde{\Lambda}_1(\mathbf{x})$ and $\widetilde{\Lambda}_2(\mathbf{x})$ are the solutions of the problems
\begin{align*}
&-\mathrm{div}\,g(\mathbf{x})\nabla \widetilde{\Lambda}_1(\mathbf{x})+v(\mathbf{x})=0,\quad\int\limits_\Omega \widetilde{\Lambda}_1(\mathbf{x})\,d\mathbf{x}=0;
\\
&-\mathrm{div}\,g(\mathbf{x})\nabla \widetilde{\Lambda}_2(\mathbf{x})+\mathrm{div}\,g(\mathbf{x})\mathbf{A}(\mathbf{x})=0,\quad \int\limits_\Omega \widetilde{\Lambda}_2(\mathbf{x})\,d\mathbf{x}=0.
\end{align*}
The column $V$ (see \eqref{V}) has the form $V=V_1+iV_2$, where $V_1$, $V_2$ are the columns with real entries defined by
\begin{align*}
&V_1=\vert \Omega\vert ^{-1}\int\limits_\Omega (\nabla \Psi (\mathbf{x}))^tg(\mathbf{x})\nabla\widetilde{\Lambda}_2(\mathbf{x})\,d\mathbf{x},
\\
&V_2=-\vert \Omega\vert ^{-1}\int\limits_\Omega (\nabla\Psi (\mathbf{x}))^tg(\mathbf{x})\nabla \widetilde{\Lambda}_1(\mathbf{x})\,d\mathbf{x}.
\end{align*}
By \eqref{W}, the constant $W$ can be written as
\begin{equation*}
W=\vert \Omega\vert ^{-1}\int\limits_\Omega \Bigl(
\langle g(\mathbf{x})\nabla \widetilde{\Lambda}_1(\mathbf{x}),\nabla \widetilde{\Lambda}_1(\mathbf{x})\rangle
+\langle g(\mathbf{x})\nabla \widetilde{\Lambda}_2(\mathbf{x}),\nabla\widetilde{\Lambda}_2(\mathbf{x})\rangle\Bigr)\,d\mathbf{x}.
\end{equation*}
The effective operator for $\mathcal{B}_\varepsilon$ acts by the rule
\begin{equation*}
\mathcal{B}^0u=-\mathrm{div}\,g^0\nabla u +2i\langle \nabla u, V_1+\overline{\boldsymbol{\eta}}\rangle
+(-W+\overline{\mathcal{Q}}+\lambda)u,\quad u\in H^2(\mathbb{R}^d).
\end{equation*}
The corresponding differential expression can be written as follows
\begin{equation*}
\mathcal{B}^0=(\mathbf{D}-\mathbf{A}^0)^*g^0(\mathbf{D}-\mathbf{A}^0)+\mathcal{V}^0+\lambda,
\end{equation*}
where
\begin{equation*}
\mathbf{A}^0=(g^0)^{-1}(V_1+\overline{g\mathbf{A}}),\quad\mathcal{V}^0=\overline{\mathcal{V}}+\overline{\langle g\mathbf{A},\mathbf{A}\rangle }-\langle g^0\mathbf{A}^0,\mathbf{A}^0\rangle -W.
\end{equation*}

\subsection{The operator $\mathcal{N}$} In the case under consideration, the structure of the operator $\mathcal{N}$ was found in \cite[(9.12)--(9.15)]{Su14}. We write out the result:
\begin{equation*}
\mathcal{N}=\sum _{k,l=1}^d \mathcal{N}_{12,kl}D_kD_l+\sum _{k=1}^d\mathcal{N}_{21,k}D_k+\mathcal{N}_{22},
\end{equation*}
where
\begin{align*}
\mathcal{N}_{12,kl}&=2\overline{\widetilde{\Lambda}_1\widetilde{g}_{kl}}
+\overline{v\psi _k\psi _l}-\sum _{j=1}^d\overline{(g_{jl}\psi _k+g_{jk}\psi _l)\partial _j\widetilde{\Lambda}_1},
\\
\mathcal{N}_{21,k}&=2\sum _{j=1}^d\overline{g_{jk}\Bigl(\widetilde{\Lambda}_1\partial _j\widetilde{\Lambda}_2-\widetilde{\Lambda}_2 \partial _j \widetilde{\Lambda}_1\Bigr)}
\\
&+2\overline{\psi _k\langle\boldsymbol{\eta},\nabla \widetilde{\Lambda}_1\rangle}
-2\overline{\widetilde{\Lambda}_1\langle\boldsymbol{\eta},\nabla \psi _k\rangle}
+2\overline{v\widetilde{\Lambda}_2\psi _k}
-4\overline{\eta _k\widetilde{\Lambda}_1},
\\
\mathcal{N}_{22}&=2\overline{\widetilde{\Lambda}_2\langle \boldsymbol{\eta},\nabla\widetilde{\Lambda}_1\rangle}
-2\overline{\widetilde{\Lambda}_1\langle\boldsymbol{\eta},\nabla\widetilde{\Lambda}_2\rangle}
+\overline{v\Bigl(\widetilde{\Lambda}_1^2+\widetilde{\Lambda}_2^2\Bigr)}
+2\overline{\widetilde{\Lambda}_1(\mathcal{Q}+\lambda)}.
\end{align*}

\subsection{Approximation for the exponential}

According to the equality $\Lambda (\mathbf{x})=i\Psi (\mathbf{x})$, the operator \eqref{K_eps 0(s) no Pi} can be written as
\begin{equation}
\label{K_eps 0 for scalar ell op}
\begin{split}
\mathcal{K}_\varepsilon ^0 (s)&
\!=\!(\Psi ^\varepsilon \nabla \!+\!\widetilde{\Lambda}^\varepsilon )e^{-\mathcal{B}^0s}
\!+\!e^{-\mathcal{B}^0s}(\Psi ^\varepsilon \nabla \!+\!\widetilde{\Lambda}^\varepsilon )^*
\!-\!\!\int\limits_0 ^s\! e^{-\mathcal{B}^0(s-\widetilde{s})}\mathcal{N}e^{-\mathcal{B}^0\widetilde{s}}\,d\widetilde{s}
\\
&\!=\!(\Psi ^\varepsilon \nabla +\widetilde{\Lambda}^\varepsilon )e^{-\mathcal{B}^0s}
\!+\!e^{-\mathcal{B}^0s}(\Psi ^\varepsilon \nabla \!+\!\widetilde{\Lambda}^\varepsilon )^*
\!-\!s\mathcal{N}e^{-\mathcal{B}^0s}.
\end{split}
\end{equation}
We have taken into account that, in the our case, the scalar differential operator 
$\mathcal{N}$ commutes with the exponential generated by the effective operator  $\mathcal{B}^0$ with constant coefficients.

From Theorem~\ref{Theorem exp B_eps final} we derive the following result.

\begin{theorem}
Let the operators $\mathcal{B}_\varepsilon$ and $\mathcal{B}^0$ be defined in Subsections~\textnormal{\ref{Subsection scalar elliptic operator}} and \textnormal{\ref{Subsection eff op for scalar case}}, respectively, and let $\mathcal{K}_\varepsilon ^0 (s)$ be the corrector \eqref{K_eps 0 for scalar ell op}. Then for $0<\varepsilon\leqslant 1$ and $s>0$ we have the approximation
\begin{equation*}
\Vert e^{-\mathcal{B}_\varepsilon s}-e^{-\mathcal{B}^0s}-\varepsilon \mathcal{K}_\varepsilon ^0 (s)\Vert _{L_2(\mathbb{R}^d)\rightarrow L_2(\mathbb{R}^d)}
\leqslant
C_{16}\varepsilon ^2 s^{-1}e^{-\check{c}_* s/2}.
\end{equation*}
The constants $C_{16}$ and $\check{c}_*$ depend only on the initial data \eqref{problem data for scalar elliptic operator}.
\end{theorem}

\end{document}